\documentclass[a4paper, twoside, 11pt, english]{article}

\usepackage[T1]{fontenc}
\usepackage[utf8]{inputenc}
\usepackage{mathtools, amsthm, amssymb, stmaryrd}
\usepackage{newtxtext, courier, euler}
\usepackage[scaled=0.95]{helvet}
\usepackage[cal=euler, scr=boondoxo, scrscaled=1.05]{mathalfa}
\usepackage{fancyhdr, titlesec, url, enumerate, microtype, setspace}
\usepackage[all, cmtip]{xy}
\usepackage{tikz}
\usepackage{babel}
\usepackage{hyperref}

\allowdisplaybreaks

\titleformat{\section}[block]{\scshape\filcenter\LARGE}{\thesection.}{.5em}{}
\titleformat{\subsection}[block]{\bfseries\filcenter\large}{\thesubsection.}{.5em}{\medskip}
\titleformat{\subsubsection}[runin]{\bfseries}{\thesubsubsection.}{.5em}{}[.]
\titlespacing{\subsubsection}{0pt}{\topsep}{.5em}

\newtheoremstyle{ntheorem}%
	{\topsep}{\topsep}{\itshape}{0pt}{\bfseries}{.}{.5em}%
	{\thmnumber{#2.\hspace{.5em}}\thmname{#1}\thmnote{ (#3)}}
	
\newtheoremstyle{ndefinition}%
	{\topsep}{\topsep}{\normalfont}{0pt}{\bfseries}{.}{.5em}%
	{\thmnumber{#2.\hspace{.5em}}\thmname{#1}\thmnote{ (#3)}}
	
\newtheoremstyle{nremark}%
	{\topsep}{\topsep}{\normalfont}{0pt}{\itshape}{.}{.5em}%
	{\thmnumber{}\thmname{#1}\thmnote{ (#3)}}

\theoremstyle{ntheorem}
  	\newtheorem{theorem}[subsubsection]{Theorem}
  	\newtheorem{proposition}[subsubsection]{Proposition}
	\newtheorem{lemma}[subsubsection]{Lemma}
  	\newtheorem{corollary}[subsubsection]{Corollary}

\theoremstyle{ndefinition}

	\newtheorem{example}[subsubsection]{Example}
	\newtheorem{remark}[subsubsection]{Remark}

\fancypagestyle{plain}{
\fancyhf{}\fancyfoot[LC,RC]{}
\fancyfoot[LE,RO]{$\thepage$}

}

\newcommand{\auteur}[3]{
\noindent
\begin{minipage}[t]{.45\textwidth}
\begin{flushright}
\textsc{#1} \\
{\footnotesize\textsf{#2}}
\end{flushright} 
\end{minipage}
\qquad
\begin{minipage}[t]{.45\textwidth}
#3
\end{minipage}
}

\newcommand{\pdf}[1]{\texorpdfstring{$#1$}{1}}

\newcommand{\fl}{\rightarrow}
\newcommand{\fll}{\longrightarrow}
\newcommand{\ofl}[1]{\overset{\displaystyle #1}{\fll}}
\newcommand{\ifl}{\rightarrowtail}
\newcommand{\pfl}{\twoheadrightarrow}
\newcommand{\fllg}{\longleftarrow}
\newcommand{\oflg}[1]{\overset{\displaystyle #1}{\fllg}}

\newcommand{\dfl}{\Rightarrow}

\newcommand{\tfl}{\Rrightarrow}
\newcommand{\tfll}{\xymatrix@1@C=20pt{\ar@3 [r] &}}
\newcommand{\otfl}[1]{\xymatrix@1@C=20pt{\ar@3 [r] ^-*+{#1} &}}

\newcommand{\pullbackcorner}[1][dr]{\save*!/#1-1.5pc/#1:(-1,1)@^{|-}\restore}

\newcommand{\ens}[1]{\{ #1\}}
\newcommand{\enspres}[2]{\{\, #1 \,\vert\, #2 \,\}}

\newcommand{\pres}[2]{\langle\, #1 \,\vert\, #2 \,\rangle}
\newcommand{\bigpres}[2]{\big\langle\; #1 \;\big\vert\; #2 \;\big\rangle}
\newcommand{\cohpres}[3]{\langle\, #1 \,\vert\, #2 \,\vert\, #3 \,\rangle}

\newcommand{\op}[1]{{#1}^{o}}
\newcommand{\env}[1]{{#1}^{e}}
\newcommand{\cl}[1]{\overline{#1}}
\newcommand{\rep}[1]{\widehat{#1}}
\renewcommand{\tilde}[1]{\widetilde{#1}}

\newcommand{\sm}{\scriptstyle}
\newcommand{\Item}[1]{{\upshape\textbf{(#1)}}}
\newcommand{\lin}[1]{\mathscr{A}(#1)}

\DeclareMathOperator{\id}{Id}

\DeclareMathOperator{\Tor}{\mathrm{Tor}}

\DeclareMathOperator{\FDT}{FDT}

\DeclareMathOperator{\FP}{FP}

\DeclareMathOperator{\Red}{Red}
\DeclareMathOperator{\red}{Red_m}
\DeclareMathOperator{\lm}{lm}
\DeclareMathOperator{\lc}{lc}
\DeclareMathOperator{\lt}{lt}
\DeclareMathOperator{\Lead}{Lead}

\DeclareMathOperator{\supp}{supp}
\DeclareMathOperator{\cell}{cell}
\DeclareMathOperator{\dom}{dom}
\DeclareMathOperator{\im}{im}

\DeclareMathOperator{\Sph}{Sph}
\DeclareMathOperator{\Sq}{Sq}
\DeclareMathOperator{\Std}{Std}
\DeclareMathOperator{\Sym}{Sym}

\newcommand{\tens}{\otimes}
\renewcommand{\leq}{\leqslant}
\renewcommand{\geq}{\geqslant}
\newcommand{\tri}{\vartriangleright}
\newcommand{\trieq}{\trianglerighteqslant}

\renewcommand{\phi}{\varphi}
\renewcommand{\epsilon}{\varepsilon}

\newcommand{\K}{\mathbb{K}}
\newcommand{\Nb}{\mathbb{N}}

\newcommand{\dr}{\partial}

\newcommand{\Br}{\mathcal{B}}
\newcommand{\Cr}{\mathcal{C}}

\newcommand{\Gr}{\mathcal{G}}

\newcommand{\catego}[1]{\text{\small $\mathsf{#1}$}}

\newcommand{\Cat}{\catego{Cat}}
\newcommand{\Gpd}{\catego{Gpd}}

\newcommand{\Bimod}{\catego{Bimod}}

\newcommand{\Set}{\catego{Set}}

\newcommand{\Vect}{\catego{Vect}}
\newcommand{\Vectg}{\catego{GrVect}}
\newcommand{\Alg}{\catego{Alg}}
\newcommand{\Algg}{\catego{GrAlg}}

\newcommand{\Gph}{\catego{Gph}}
\newcommand{\Glob}{\catego{Glob}}
\newcommand{\Ind}{\catego{Ind}}
\newcommand{\Pol}{\catego{Pol}}

\begin{document}

\thispagestyle{empty}
 
\begin{center}

\begin{doublespace}
\begin{huge}
{\scshape Convergent presentations and polygraphic resolutions of associative algebras}
\end{huge}

\bigskip
\hrule height 1.5pt 
\bigskip

\begin{Large}
{\scshape Yves Guiraud \qquad Eric Hoffbeck \qquad Philippe Malbos}
\end{Large}
\end{doublespace}

\vfill

\begin{small}\begin{minipage}{14cm}
\noindent\textbf{Abstract --}
Several constructive homological methods based on noncommutative Gröbner bases are known to compute free resolutions of associative algebras. In particular, these methods relate the Koszul property for an associative algebra to the existence of a quadratic Gröbner basis of its ideal of relations. In this article, using a higher-dimensional rewriting theory approach, we give several improvements of these methods. We define polygraphs for associative algebras as higher-dimensional linear rewriting systems that  generalise the notion of noncommutative Gröbner bases, and allow more possibilities of termination orders than those associated to monomial orders. We introduce polygraphic resolutions of associative algebras, giving a categorical description of higher-dimensional syzygies for presentations of algebras. We show how to compute polygraphic resolutions starting from a convergent presentation, and how these resolutions can be linked with the Koszul property.


\smallskip
\noindent\textbf{M.S.C. 2010 --} 
18C10, 18D05, 18G10, 16S37, 68Q42.

\bigskip
\hfill December 21, 2017
\end{minipage}\end{small}

\begin{small}
\renewcommand{\contentsname}{}
\setcounter{tocdepth}{2}
\tableofcontents
\end{small}

\end{center}

\section*{Introduction}

\subsection*{An overview on rewriting in algebras}

\subsubsection*{Linear rewriting}

Throughout the twentieth century, several rewriting approaches have been developed for computations in an algebraic setting. Rewriting consists in orienting the relations of a presentation, and  computational properties are deduced from the overlaps that may appear in the application of the oriented relations. A \emph{rewriting} is the application of a sequence of oriented relations, and the \emph{confluence} property says that all the rewritings reducing a given term end on a common term. Newman's lemma,~\cite{Newman42}, states that under the termination hypothesis, meaning that every rewriting ends, confluence is equivalent to confluence of the one-step rewritings, a property called \emph{local confluence}. The critical branching lemma,~\cite{Nivat73, KnuthBendix70, Huet80, GuiraudMalbos16}, states that local confluence can in turn be deduced from the confluence of \emph{critical branchings}, which are some minimal overlaps of two oriented relations. These two results allow to deduce confluence from a local analysis of branchings, and in this way the confluence property can be used to compute normal forms for algebras presented by generators and relations, with applications to the decision of the word problem (or of ideal membership) and to the construction of bases, such as Poincaré-Birkhoff-Witt bases. 

This is the main principle applied in numerous works in algebra. For example, Shirshov introduced in~\cite{Shirshov62}  an algorithm to compute a linear basis of a Lie algebra presented through generators and relations. He used the notion of \emph{composition} of elements in a free Lie algebra to describe the critical branchings. He gave an algorithm to compute bases in free algebras having the confluence property, and he proved the \emph{composition lemma}, which is the analogue of Newman's lemma for Lie algebras. In particular, he deduced a constructive proof of the Poincaré-Birkhoff-Witt theorem. Also, rewriting methods to compute with ideals of commutative polynomial rings were introduced by Buchberger with Gröbner basis theory. He described critical branchings with the notion of \emph{$S$-polynomial} and gave an algorithm for the computation of Gröbner bases, sets of relations having the confluence property,~\cite{Buchberger65, Buchberger06, Buchberger87}. Buchberger's algorithm, applied to a finite family of generators of an ideal of a commutative polynomial ring, always terminates and returns a Gröbner basis of the ideal. The central theorem for Gröbner basis theory is the counterpart of Newman's lemma, and Buchberger's algorithm is in essence the analogue of Knuth-Bendix's completion procedure in a linear setting, \cite{Buchberger87}. Note that, in the same period, ideas in the spirit of the Gröbner basis approach appear in several others works, by Hironaka in~\cite{Hironaka64} and Grauert in~\cite{Grauert72}, with \emph{standard bases} for power series rings, or for applications of Newman's lemma for universal algebras by Cohn, \cite{Cohn65}. The domain took foundation in several works on algorithmic methods in elimination theory by Macaulay,~\cite{Macaulay16}, with \emph{$H$-bases}, by Janet,~\cite{Janet20} with \emph{involutive bases}, or Gunther,~\cite{Gunther41}, with notions similar to Gröbner bases.

\subsubsection*{Noncommutative Gröbner bases}

Bokut and Bergman have independently extended completion methods to associative algebras, \cite{Bokut76, Bergman78}. They obtained Newman's lemma for rewriting systems in free associative algebras compatible with a monomial order, called respectively the \emph{composition lemma} and the \emph{diamond lemma}. 
Given an associative algebra~$A$ presented by a set of generators~$X$ and a set of relations~$R$,
the monomials of the free algebra~$\lin{X}$ on~$X$ form a linear basis of~$\lin{X}$.
One important application of noncommutative Gröbner bases is to explicitly find a basis of the algebra~$A$ in the form of a subset of the set of monomials on~$X$. This is based on a monomial order, that is a well-founded total order on the monomials, and the idea is to change the presentation of the ideal generated by~$R$ with respect to this order.
The property that the new presentation has to satisfy is the algebraic counterpart of the confluence property. 
The set of reduced monomials with respect to a Gröbner basis~$\mathcal{G}$ forms a linear basis of the quotient of the free algebra~$\lin{X}$ by the ideal generated by~$\mathcal{G}$.
In general, completion procedures do not terminate for ideals in a noncommutative polynomial ring. Even if the ideal is finitely generated, it may not have a finite Gröbner basis, but an infinite Gröbner basis can be computed over a field, see~\cite{Mora94, Ufnarovski98}.
Subsequently, rewriting methods were developed for a wide range of algebraic structures, such as Weyl algebras~\cite{SaitoSturmfelsTakayama00}, or operads~\cite{DotsenkoKhoroshkin10}. For a comprehensive treatment on noncommutative Gröbner bases we refer the reader to~\cite{BremnerDotsenko16, Mora94, Ufnarovski95}.

\subsubsection*{Free resolutions of associative algebras from confluence}

At the end of the eighties, through Anick's and Green's works~\cite{Anick85, Anick86, Anick87, Green99}, noncommutative Gröbner bases have found new applications giving constructive methods to compute free resolutions of associative algebras.
Their constructions provide small explicit resolutions to compute homological invariants (homology groups, Hilbert and Poincaré series) of algebras presented by generators and relations given by Gröbner bases. Anick's resolution consists in a complex generated by \emph{Anick's chains}, that are certain iterated overlaps of the leading terms of the relations, and whose differential is obtained by deforming the differential of a complex for an associated monomial algebra. This construction has many applications, such as an algorithm for the computation of Hilbert series,~\cite{Ufnarovski95}. The chains and the differential of the resolution are constructed recursively, making possible  its implementation,~\cite{BackelinCojocaruUfnarovski03}, but the differential is complicated to make explicit in general. Sköldberg introduced in~\cite{Skoldberg06} a homotopical method based on discrete Morse theory to derive Anick's resolution from the bar resolution.

\subsubsection*{Confluence and Koszulness}

Anick's resolution also gives a relation between the Koszul property of an algebra and the existence of a quadratic noncommutative Gröbner basis for its ideal of relations.
Recall that a connected graded algebra~$A$ is \emph{Koszul} if the Tor groups $\Tor^A_{k,(i)}(\K,\K)$ vanish for~$i\neq k$, where~$k$ is the homological degree and~$i$ corresponds to the grading of the algebra. 
Koszul algebras were introduced by Priddy, \cite{Priddy70}, and he proved that quadratic algebras having a Poincaré-Birkhoff-Witt basis are Koszul.
This notion was generalised by Berger in~\cite{Berger01} to the case of $N$-homogeneous algebras, asking that $\Tor^A_{k,(i)}(\K,\K)$ vanish for $i\neq \ell_N(k)$, where $\ell_N:\Nb\to\Nb$ is the function defined by $\ell_N(k) = lN$ if~$k=2l$, and~$lN+1$ if~$k=2l+1$. For a graded algebra~$A$ with $N$-homogeneous relations, the groups $\textrm{Tor}_{k,(i)}^{A} (\K, \K)$ always vanish for~$i<\ell_N(k)$, \cite{BergerMarconnet06}, so that the Koszul property corresponds to the limit case.
Koszulness of a quadratic algebra~$A$ can be obtained by showing the existence of a quadratic Gröbner basis, because it implies the existence of a Poincaré-Birkhoff-Witt basis of~$A$,~\cite{Green99}, and thus the Koszulness of~$A$ by Priddy's result.
For the $N$-homogeneous case, a Gröbner basis concentrated in weight~$N$ is not enough to imply Koszulness: an extra condition has to be checked as shown by Berger,~\cite{Berger01}.
When the algebra is monomial, this extra condition corresponds to the \emph{overlap property} defined by Berger, which consists in a combinatorial condition based on the overlaps of the monomials of the relations.
Anick's resolution generalises Priddy's results on Koszulness of quadratic algebras. Indeed, if an algebra~$A$ has a quadratic Gröbner basis, then Anick's resolution is concentrated in the right degrees, and thus~$A$ is Koszul,~\cite{Anick86, GreenHuang95}.
Note that Backelin gave a characterisation of the Koszul property for quadratic algebras in terms of lattices,~\cite{Backelin83, BackelinFroberg85}, and this condition was later interpreted in term of \emph{$X$-confluence} by Berger, \cite{Berger98}.
Berger extended this construction to the case of $N$-homogeneous algebras, where an extra condition has to be checked~\cite{Berger01}. Finally, the quadratic Gröbner basis method to prove Koszulness has been extended to the case of operads by Dotsenko and Khoroshkin in~\cite{DotsenkoKhoroshkin10}.

\subsection*{Polygraphic resolutions of associative algebras}

Known constructions of free resolutions using confluence, such as Anick's resolution, are not explicit and constructed inductively with respect to a monomial ordering. One of the objectives of this paper is to give such a construction in a higher-dimensional rewriting framework, in order to make explicit the contracting homotopy in terms of rewriting properties of the presentation, and using a non-monomial orientation of relations.
Higher-dimensional rewriting theory is the theory of presentations by generators and relations of higher-dimensional categories,~\cite{Street76,Street87,Burroni93}. The notion of \emph{polygraph} is the main concept of the theory, extending to higher dimensions the notion of presentation by generators and relations for categories, and giving an unified paradigm of rewriting,~\cite{GuiraudMalbos09,GuiraudMalbos16}. In this work, we introduce a variation of the notion of polygraph for higher-dimensional associative algebras. 

\subsubsection*{Higher-dimensional associative algebras}

The first section of the paper introduces the main categorical notions used throughout this work.
We define in~\ref{SSS:HigherDimensionalAssociativeAlgebras} the category~$\infty\Alg$ of \emph{$\infty$-algebras} and \emph{morphisms of $\infty$-algebras} as the category of $\infty$-categories and $\infty$-functors internal to the category~$\Alg$ of associative algebras. Theorem~\ref{T:nAlg} gives several interpretations of the structure of $\infty$-algebra. In particular,  the category~$\infty\Alg$ is isomorphic to the category of internal $\infty$-groupoids in the category~$\Alg$. Theorem~\ref{T:nAlg} also makes explicit isomorphisms between these categories and full subcategories of the category of globular algebras and the category of globular bimodules. The latter interpretation is used to construct free $\infty$-algebras in~\ref{SSS:FreeHigherAlgebras}.

\subsubsection*{Polygraphs for associative algebras}

In Section~\ref{S:nPolygraphs}, we adapt the set-theoretical notion of polygraph to presentations of higher-dimensional algebras. As in the set-theoretical case, the category $n\Pol(\Alg)$ of \emph{$n$-polygraphs for associative algebras}, called \emph{$n$-polygraphs} for short, is constructed by induction on~$n$ in~\ref{SSS:PolygraphFiniteDimension}. The category $0\Pol(\Alg)$ is the category of sets. Then \emph{$1$-polygraphs} are pairs $\langle X_0\,|\,X_1\rangle$ made of a set~$X_0$ and a \emph{cellular extension}~$X_1$ of the free algebra~$\lin{X_0}$ over~$X_0$, that is a set~$X_1$ with two maps
\[
\xymatrix @=1.5em {
\lin{X_0}
& X_1.
	\ar@<-0.5ex> [l] _-{s}
	\ar@<0.5ex> [l] ^-{t}
}
\]
The elements of $X_k$ are called the \emph{$k$-cells} of~$X$, and for a $1$-cell~$\alpha$, the $0$-cells~$s(\alpha)$ and~$t(\alpha)$ are its \emph{source} and \emph{target}. For instance the $1$-polygraph
\begin{equation}
\label{Equation:X:xyz=x3+y3+z3}
X = \bigpres{x,y,z}{xyz \ofl{\gamma} x^3 + y^3 + z^3}
\end{equation}
studied in~\ref{X:xyz=x3+y3+z3} is generated by three $0$-cells~$x$, $y$, and~$z$, and a $1$-cell~$\gamma$ that reduces the monomial~$xyz$ into $x^3+y^3+z^3$.
For~$n\geq 2$, assume that the category $(n-1)\Pol(\Alg)$ of $(n-1)$-polygraphs is defined, together with the free $(n-1)$-algebra functor $(n-1)\Pol(\Alg) \fl (n-1)\Alg$. The category~$n\Pol(\Alg)$ of $n$-polygraphs is defined as the pullback
\[
\vcenter{\xymatrix @=1.5em{
n\Pol(\Alg)
  \ar[r]
  \ar[d] \pullbackcorner
&
(n-1)\Alg^+
  \ar[d] 
\\
(n-1)\Pol(\Alg)
  \ar[r] 
&
(n-1)\Alg
}}
\]
where $(n-1)\Alg^+$ is the category of $(n-1)$-algebras $A_{n-1}$ with a cellular extension, i.e.\ a set~$X$ with two maps
\[
\xymatrix @=1.5em {
A_{n-1} 
& X
	\ar@<-0.5ex> [l] _-{s}
	\ar@<0.5ex> [l] ^-{t}
}
\]
such that~$s(\alpha)$ and~$t(\alpha)$ have the same source and the same target.
An $n$-polygraph is thus a sequence $\cohpres{X_0}{\cdots}{X_n}$, made of a set~$X_0$ and, for every $0\leq k<n$, a cellular extension~$X_{k+1}$ of the free $k$-algebra over $\cohpres{X_0}{\cdots}{X_k}$.
If~$X$ is an $n$-polygraph, we denote by~$\lin{X}$ the free $n$-algebra over~$X$.
The \emph{algebra presented} by a $1$-polygraph~$X$ is the quotient algebra~$\cl{X}$ of the free algebra~$\lin{X_0}$ by the congruence generated by~$X_1$. A \emph{coherent presentation} of an algebra~$A$ is a $2$-polygraph $\cohpres{X_0}{X_1}{X_2}$ whose underlying $1$-polygraph $\pres{X_0}{X_1}$ is a presentation of~$A$, and such that the cellular extension~$X_2$ is \emph{acyclic}, meaning that every $1$-sphere in the free $1$-algebra~$\lin{X_1}$ is trivial with respect to the congruence generated by~$X_2$.

\subsection*{Convergent presentation of associative algebras}

In Section~\ref{S:LinearRewriting}, we study the rewriting properties of $1$-polygraphs, whose  $1$-cells are not necessarily oriented with respect to a monomial ordering.  Our approach is thus less restrictive than those known for associative algebras, which rely on a monomial order. For instance, there is no monomial order compatible with the rule~$\gamma$ of the $1$-polygraph~\eqref{Equation:X:xyz=x3+y3+z3}, see~\ref{X:xyz=x3+y3+z3}, whereas this orientation makes proving its confluence trivial.

\subsubsection*{One-dimensional polygraphic rewriting}

A $1$-polygraph $X$ is \emph{left-monomial} if, for every $1$-cell~$\alpha$ of~$X$, the $0$-cell~$s(\alpha)$ is a monomial of~$\lin{X}$ that does not belong to the support of~$t(\alpha)$, as defined in~\ref{SSS:LeftMonomialPolygraph}. 
We define in~\ref{SSS:RewritingSteps} a \emph{rewriting step} of a left-monomial $1$-polygraph~$X$, as a $1$-cell $\lambda f + 1_a$ of the free $1$-algebra~$\lin{X}$ such that~$\lambda \neq 0$ and~$s(f)$ is not in the support of~$a$. A composition of rewriting steps in~$\lin{X}$ is called a \emph{positive} $1$-cell of~$\lin{X}$. A $0$-cell~$a$ of~$\lin{X}$ is \emph{reduced} if there is no rewriting step of source~$a$. The reduced $0$-cells of~$\lin{X}$ form a linear subspace~$\Red(X)$ of the free algebra~$\lin{X_0}$, admitting as a basis the set~$\red(X)$ of reduced monomials of~$\lin{X_0}$. 
We define in~\ref{SSS:RewriteOrderTermination} the \emph{termination} of~$X$ as the wellfoundedness of a binary relation on monomials induced by rewriting steps.
A \emph{branching} of~$X$ is a pair~$(f,g)$ of positive $1$-cells of~$\lin{X}$ with the same source. Such a branching is called \emph{confluent} if there exist positive $1$-cells~$h$ and~$k$ of~$\lin{X}$ as in
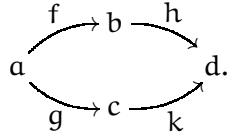
\[
\vcenter{\xymatrix @R=0.25em {
& b
	\ar@/^/ [dr] ^-{h}
\\
a
	\ar@/^/ [ur] ^-{f}
	\ar@/_/ [dr] _-{g}
&& d.
\\
& c
	\ar@/_/ [ur] _-{k}
}}
\]
One says that $X$ is \emph{confluent} if all its branchings are confluent. One defines a less restrictive notion of \emph{critical confluence} by requiring the confluence of \emph{critical branchings} only, which are branchings involving generating $1$-cells whose source overlap. A $1$-polygraph is \emph{convergent} if it terminates and it is confluent. In that case, every $0$-cell~$a$ of~$\lin{X}$ has a unique \emph{normal form}, denoted by $\rep{a}$, which is a reduced $0$-cell of~$\lin{X}$ in which~$a$ can be rewritten. Under the termination hypothesis, the confluence property is equivalent to saying that the vector space $\lin{X_0}$ admits the direct decomposition $\lin{X_0} = \Red(X) \oplus I(X)$, as proved by Proposition~\ref{P:CharacterisationConfluence}. As an immediate consequence of this decomposition, we deduce the following basis theorem.
\begin{quote}
{\bfseries Theorem~\ref{T:StandardBasis}.}
\emph{
Let~$A$ be an algebra and~$X$ be a convergent presentation of~$A$. Then~$\red(X)$ is a linear basis of~$A$. As a consequence, the vector space~$\Red(X)$, equipped with the product defined by $a\cdot b=\rep{ab}$, is an algebra that is isomorphic to~$A$.
}
\end{quote}
Finally, Propositions~\ref{P:ConvergentGröbner} and~\ref{P:ConvergentPBW} show that convergent $1$-polygraphs generalise noncommutative Gröbner bases and, in the case of $N$-homogeneous algebras, Poincaré-Birkhoff-Witt bases. In particular, for quadratic algebras and the deglex order, we recover Priddy's concept of Poincaré-Birkhoff-Witt basis as introduced in~\cite[Section 5.1]{Priddy70}.

\subsubsection*{The coherent critical branchings theorem}

The goal of Section~\ref{S:Squier} is to give coherent formulations of confluence results on $1$-polygraphs.
 A branching~$(f,g)$ of a left-monomial $1$-polygraph~$X$ is \emph{$Y$-confluent} with respect to a cellular extension~$Y$ of the free $1$-algebra~$\lin{X}$, if there exist positive $1$-cells~$h$ and~$k$ in~$\lin{X}$, and a $2$-cell~$F$ in the free $2$-algebra generated by the $2$-polygraph $\langle X\, |\, Y\rangle$ of the form
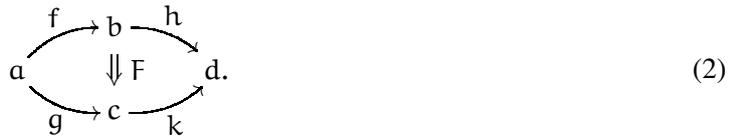
\begin{equation}
\label{E:YConfluence}
\vcenter{\xymatrix @R=0.25em {
& b
	\ar@/^/ [dr] ^-{h}
	\ar@2 []!<0pt,-10pt>;[dd]!<0pt,10pt> ^-*+{F}
\\
a
	\ar@/^/ [ur] ^-{f}
	\ar@/_/ [dr] _-{g}
&& d.
\\
& c
	\ar@/_/ [ur] _-{k}
}}
\end{equation}
When~$Y$ is the set $\Sph(\lin{X})$ of all the spheres created by the $1$-cells of~$\lin{X}$, we recover the notions of confluence, local confluence, critical confluence.
If~$X$ is terminating, we obtain two important properties: the \emph{coherent Newman's lemma}, Proposition~\ref{P:HNewman}, stating that local $Y$-confluence implies $Y$-confluence, and the \emph{coherent critical branching theorem}, Theorem~\ref{T:HCritical}, stating that critical $Y$-confluence implies $Y$-confluence.
When $Y=\Sph(\lin{X})$, we recover respectively Newman's lemma and the critical branching theorem. For instance, the $1$-polygraph $X$ given by (\ref{Equation:X:xyz=x3+y3+z3}) is terminating and has no critical branching, hence it is trivially confluent.
As a consequence of the critical branching theorem, Corollary~\ref{C:BuchbergerCriterion} gives a polygraphic interpretation of Buchberger's criterion for noncommutative Gröbner bases.
As it requires termination, the critical branching theorem in this linear setting differs from its set-theoretic counterpart. Indeed, nonoverlapping branchings are always confluent in the set-theoretic case, but we show that it is not the case for linear rewriting systems, for which confluence of non-overlapping branchings can depend on critical confluence. Finally, Section~\ref{S:Squier} gives a linearised version of Squier's theorem, intially stated for presentations of monoids~\cite{Squier94}, see also~\cite{GuiraudMalbos16}.
\begin{quote}
{\bfseries Theorem~\ref{T:Squier}.}
\emph{
Let~$X$ be a convergent left-monomial $1$-polygraph. A cellular extension~$Y$ of~$\lin{X}$ that contains one $2$-cell of the form~\eqref{E:YConfluence} for every critical branching~$(f,g)$ of~$X$, with~$h$ and~$k$ positive $1$-cells of~$\lin{X}$, is acyclic.
}
\end{quote}
This result is then extended in Section~\ref{S:SquierResolution} into a polygraphic resolution for the algebra~$\cl{X}$ presented by~$X$, involving critical branchings in every dimension.

\subsection*{Construction and reduction of polygraphic resolutions of associative algebras}

A \emph{polygraphic resolution} of an algebra~$A$ is an $\infty$-polygraph~$X$ whose underlying polygraph is a presentation of~$A$, and such that, for every~$k\geq 1$, the cellular extension~$X_{k+1}$ of the free $k$-algebra~$\lin{A_k}$ is acyclic. In Sections~\ref{S:PolygraphicResolutions} and~\ref{S:SquierResolution} we construct such a resolution starting with a convergent presentation of the algebra~$A$. Finally, in Section~\ref{S:FreeResolutions}, we deduce a resolution of~$A$ by free $A$-bimodules from a polygraphic resolution of~$A$.

\subsubsection*{Contractions of polygraphs}

A method to construct a polygraphic resolution is to consider a contraction inducing a notion of normal form in every dimension, together with a homotopically coherent reduction of every cell to its normal form. This notion was introduced in~\cite{GuiraudMalbos12advances} for presentations of categories, where it was called a normalisation strategy, and provides a constructive characterisation of the acyclicity of an  $(\infty,1)$-polygraph. 
In Subsection~\ref{SS:Contractions}, we introduce contractions for polygraphic resolutions of algebras.
We prove that a polygraphic resolution of an algebra~$A$ is equivalent to the data of an $\infty$-polygraph whose underlying $1$-polygraph is a presentation of~$A$ and equipped with a contraction. Explicitely, a \emph{unital section} of an $\infty$-polygraph~$X$ is a section of the canonical projection $\pi:\lin{X}\pfl\cl{X}$ mapping~$1$ to~$A$. Given such a unital section~$\iota$, an \emph{$\iota$-contraction of~$X$} is a homotopy~$\sigma: \id_{\lin{X}}\dfl\iota\pi$ that satisfies $\sigma_a = 1_a$ on the images of~$\iota$ and of~$\sigma$, and an $\iota$-contraction is called \emph{right} if it satisfies $\sigma_{ab}=s_0(a)\sigma_b \star_0 \sigma_{a\rep{s_0(b)}}$. The main result of Section~\ref{S:PolygraphicResolutions} relates the property for an $\infty$-polygraph to be a polygraphic resolution to the existence of a right $\iota$-contraction.
\begin{quote}
{\bfseries Theorem~\ref{T:ResolutionContraction}.}
\emph{
Let~$X$ be an $\infty$-polygraph with a unital section~$\iota$. Then~$X$ is a polygraphic resolution of~$\cl{X}$ if, and only if, it admits a right $\iota$-contraction.
}
\end{quote}

\subsubsection*{The standard polygraphic resolution}

Given an augmented algebra~$A$ and a linear basis~$\Br$ of its positive part, Theorem~\ref{T:StandardPolygraphicResolution} makes explicit a polygraphic resolution~$\Std(\Br)$, which is a cubical analogue of the standard resolution of the algebra $A$. Subsection~\ref{SS:StandardPolygraphicResolution} is devoted to the proof of the acyclicity of the resolution by exhibiting a right contraction for the resolution~$\Std(\Br)$. We can apply this construction when~$X$ is a convergent presentation of~$A$: Corollary~\ref{C:ConvergentStandardPolygraphicResolution} makes explicit a polygraphic resolution $\Std(\red(X))$ for~$A$, whose generating cells are finite families of nontrivial reduced monomials of the free algebra~$\lin{X}$. This resolution extends the coherent presentation given by Squier's theorem, Theorem~\ref{T:Squier}.

\subsubsection*{Collapsing, and Squier's polygraphic resolution}

The polygraphic resolution $\Std(\red(X))$ is too large in general to be used in practice. Following a construction given by Brown in~\cite{Brown92}, in Subsection~\ref{SS:Collapsing} we define polygraphic \emph{collapsing schemes}, which give, by a process similar to algebraic Morse theory for chain complexes,~\cite{Skoldberg06}, a method to contract a polygraphic resolution of an algebra into a smaller one, see Theorem~\ref{T:Collapsing}. We then use this result to contract the standard resolution $\Std(\red(X))$ into a smaller one, denoted by $\Sq(X)$ and called \emph{Squier's polygraphic resolution}, containing only the critical $n$-branchings of~$X$, which are the overlaps of~$n$ generating relations. Omitting the technical details, we obtain:

\begin{quote}
{\bfseries Theorem~\ref{T:SquierResolution}.}
\emph{
If~$X$ is a reduced convergent left-monomial presentation of an augmented algebra~$A$, then~$\Sq(X)$ is a polygraphic resolution of~$A$.
}
\end{quote}
By this result, any reduced convergent left-monomial presentation of an augmented algebra~$A$ extends to a polygraphic resolution of~$A$ whose generating cells correspond to the iterated overlaps of leading terms of relations: in this spirit, Squier's resolution is a categorical analogue of Anick's resolution of an augmented algebra presented by a Gröbner basis. In Subsection~\ref{SS:ExamplesSymmetricAlgebra}, we compute Squier's polygraphic resolutions for the symmetric algebra, the quantum deformation of the symmetric algebra, and the exterior algebra.

\subsubsection*{Free resolutions from polygraphic resolutions}

In Section~\ref{S:FreeResolutions}, we show that a polygraphic resolution of an algebra~$A$ induces free resolutions in categories of modules over~$A$.
Given an $\infty$-polygraph~$X$ whose underlying $1$-polygraph is presentation of~$A$, we construct in Subsection~\ref{SS:FreeBimodulesResolutions} a complex of $A$-bimodules, denoted by~$\env{A}[X]$, whose boundary maps are induced by the source and target maps of the polygraph. 
We prove that if the polygraph~$X$ is acyclic, then the induced complex is acyclic, giving the following result.
\begin{quote}
{\bfseries Theorem~\ref{T:PolRes->AbRes}.}
\emph{
If~$X$ is a polygraphic resolution of an algebra~$A$, then the complex~$\env{A}[X]$ is a free resolution of the $A$-bimodule~$A$. Moreover, if~$X$ is of finite type, then so is~$\env{A}[X]$.
}
\end{quote}
This bimodule resolution can be used to compute Hochschild homology, as in \cite[Section 5]{Berger01}. Using these constructions, we deduce finiteness homological properties of an associative algebra~$A$ given by a convergent presentation.
In~\ref{SSS:AlgebrasFiniteDerivationType},  we introduce the property of \emph{finite $n$-derivation type} for an associative algebra, that corresponds to admitting a polygraphic resolution with finitely many generating $k$-cells for $k<n$. We relate this property to a homological finiteness condition, \emph{type~$\FP_n$}, and we prove that associative algebras admitting a finite convergent presentation are of finite $n$-derivation type for any natural number~$n$, Proposition~\ref{Proposition:FDTnormalisationStrategy}.

\subsubsection*{Confluence and Koszulness}

In Section~\ref{S:ConvergenceKoszulness}, we apply our constructions to study Koszulness of associative algebras.
Given a map $\omega: \Nb \to \Nb$, we call a polygraphic resolution \emph{$\omega$-concentrated} 
if, for any integer~$k$, all its $k$-cells are concentrated in degree~$\omega(k)$.
Similarly, a free resolution~$F_\ast$ of bimodules is \emph{$\omega$-concentrated} if, for any integer~$k$, 
the bimodule~$F_k$ is generated in degree~$\omega(k)$. 
As a consequence of Theorem~\ref{T:PolRes->AbRes}, we deduce that if~$X$ is an $\omega$-concentrated polygraphic resolution of a graded algebra~$A$, then the free $A$-bimodule~$\env{A}[X_n]$ is generated by its component of degree $\omega(n+1)$, Proposition~\ref{P:Koszul}. In particular, if~$\omega = \ell_N$ for some~$N\geq 2$, this proves that the algebra~$A$ is Koszul. 
Consequently, an algebra admitting a quadratic convergent presentation is Koszul. Finally, we discuss several examples applying our rewriting methods to prove or disprove Koszulness. For instance, the algebra presented by the $1$-polygraph~$X$ given by~\eqref{Equation:X:xyz=x3+y3+z3} is trivially Koszul, whereas a proof of this property using Anick's resolution involves the computation of a resolution of infinite length.

\subsection*{Organisation of the article and conventions}

Section~\ref{S:nAlg} presents some categorical background, and constructions on the structure of higher-dimensional associative algebra. Section~\ref{S:nPolygraphs} deals with the notion of polygraph for associative algebras, to define presentations, coherent presentations and resolutions for these algebras. Section~\ref{S:LinearRewriting} is devoted to the study of termination and confluence properties of polygraphs, making comparisons with Poincaré-Birkhoff-Witt bases and Gröbner bases. Section~\ref{S:Squier} deals with coherent presentations of associative algebras and the construction of such presentations using convergence. This is the first step of the construction achieved in Sections~\ref{S:PolygraphicResolutions} and~\ref{S:SquierResolution} of a polygraphic resolution of an associative algebra starting with a convergent presentation. Section~\ref{S:PolygraphicResolutions} presents homotopical operations on polygraphic resolutions and a method to contract polygraphic resolutions into smaller ones. 
In Section~\ref{S:SquierResolution}, these operations are applied to construct Squier's polygraphic resolution of an associative algebra starting with a reduced convergent presentation. Finally, Section~\ref{S:FreeResolutions} is devoted to the application of polygraphic resolutions to the construction of resolutions of associative algebras by free bimodules, leading to finiteness conditions and several necessary or sufficient conditions for an associative algebra to be Koszul. 

\medskip
We fix a field~$\K$ for the whole article, and denote by~$\Vect$ the category of vector spaces over~$\K$. The category of unital and associative algebras over~$\K$ is denoted by~$\Alg$. In this article, we will say algebra for unital associative algebra if no confusion may arise. 

\section{Higher-dimensional associative algebras}
\label{S:nAlg}

This section introduces the higher-dimensional objects used throughout the paper: higher-dimensional vector spaces and higher-dimensional associative algebras, defined as (globular, strict) higher-dimensional categories internal to vector spaces and to associative algebras, respectively. Our notion of higher-dimensional vector space extends the $2$-vector spaces defined by Baez and Crans in~\cite{BaezCrans04}, but with a shift by~$1$ in the dimension: our $n$-vector spaces are 
$n$-categories in~$\Vect$, instead of $(n-1)$-categories. The main result of the section, Theorem~\ref{T:nAlg}, explores the structure of $\infty$-algebras, giving equivalences with other, simpler structures: this is used in the next section to build free $\infty$-algebras.

\subsection{Internal higher-dimensional categories}
\label{SS:nCat}

Let~$\Cr$ be a fixed category. The definitions of globular objects of~$\Cr$ and $\infty$-categories of~$\Cr$ can be given in a more abstract setting, but we assume here that~$\Cr$ is concrete over the category of sets, and that the corresponding forgetful functor admits a left adjoint.

\subsubsection{Indexed objects and morphisms}

For~$I\subseteq\Nb$, an \emph{$I$-indexed object of~$\Cr$} is a sequence $X=(X_n)_{n\in I}$ of objects of~$\Cr$. In what follows, we just say \emph{indexed object} when~$I=\Nb$. 
For~$I,J\subseteq\Nb$, for~$X$ an $I$-indexed object of~$\Cr$, for~$Y$ a $J$-indexed object of~$\Cr$, and for~$p$ an integer, an \emph{indexed morphism of degree~$p$ from~$X$ to~$Y$} is a sequence 
\[
F \:=\: \big( X_i \fl Y_{i+p} \big)_{
\begin{subarray}{l}
i \,\in\, I \\[0.2ex] i+p \,\in\, J
\end{subarray}
}
\]
of morphisms of~$\Cr$. In particular, if~$I=J=\Nb$ and~$p=-1$, then the index~$i$ ranges over~$\Nb\setminus\ens{0}$. 
If~$X$ is an indexed object of~$\Cr$, and~$n\geq 0$, we abusively denote by~$X_n$ the indexed object of~$\Cr$ that is constantly equal to~$X_n$.
Indexed objects (over any possible sets of indices) and indexed morphisms (of all possible degrees) of~$\Cr$ form a category, denoted by~$\Ind(\Cr)$. If~$\Cr$ has limits (resp.\ colimits), then so does~$\Ind(\Cr)$, and those limits (resp.\ colimits) are computed pointwise.

\subsubsection{Internal globular objects}

A \emph{globular object of~$\Cr$} is an indexed object $X=(X_n)_{n\geq 0}$ of~$\Cr$ equipped with indexed morphisms 
\[
X \ofl{s} X,
\qquad
X \ofl{t} X
\qquad\text{and}\qquad
X \ofl{i} X,
\]
of respective degrees~$-1$, $-1$ and~$1$, called the \emph{source map}, the \emph{target map} and the \emph{identity map of~$X$}, that satisfy the following equalities, collectively referred to as the \emph{globular relations of~$X$}:
\[
ss = st, 
\qquad 
ts = tt
\qquad\text{and}\qquad
 si = ti = \id_{X}.
\]
With the identity map and the last two relations removed, one gets a \emph{semiglobular object of~$\Cr$}.
Given two globular objects~$X$ and~$Y$ of~$\Cr$, a \emph{globular morphism from~$X$ to~$Y$} is an indexed morphism $F:X\fl Y$ of degree~$0$ that commutes with the source, target and identity maps:
\[
\vcenter{\xymatrix @=1.5em {
X
	\ar [r] ^-{s}
	\ar [d] _-{F}
	\ar@{} [dr] |{\sm =}
& X
	\ar [d] ^-{F}
\\
Y
	\ar [r] _-{s}
& Y
}}
\qquad\qquad
\vcenter{\xymatrix @=1.5em{
X
	\ar [r] ^-{t}
	\ar [d] _-{F}
	\ar@{} [dr] |{\sm =}
& X
	\ar [d] ^-{F}
\\
Y
	\ar [r] _-{t}
& Y
}}
\qquad\qquad
\vcenter{\xymatrix @=1.5em{
X
	\ar [r] ^-{i}
	\ar [d] _-{F}
	\ar@{} [dr] |{\sm =}
& X
	\ar [d] ^-{F}
\\
Y
	\ar [r] _-{i}
& Y
}}
\]
We denote by~$\Glob(\Cr)$ the category of globular objects and globular morphisms of~$\Cr$.

For~$n\geq 0$, an \emph{$n$-globular object of~$\Cr$} is defined in the same way as a globular object of~$\Cr$, but starting with a $\ens{0,\dots,n}$-indexed object of~$\Cr$. We denote by $n\Glob(\Cr)$ the category they form with the corresponding $n$-globular morphisms of~$\Cr$. 

\subsubsection{Sources, targets and spheres}

Let~$X$ be globular object of~$\Cr$, and fix~$n\geq 0$. An element~$x$ of~$X_n$ is called an \emph{$n$-cell of~$X$}, and, if~$n\geq 1$, the $(n-1)$-cells~$s(x)$ and~$t(x)$ are called the \emph{source of~$x$} and the \emph{target of~$x$}. We write $z:x\fl y$ if~$z$ is an $n$-cell of~$X$ of source~$x$ and target~$y$, and we use the more specific notation $z:x\dfl y$ (resp.\ $z:x\tfl y$) when~$n=2$ (resp.~$n=3$). Since the globular relations imply that~$i$ is injective, and when no confusion occurs, we just write~$i(x)$ or~$1_x$, or even just~$x$, instead of any iterate image~$i^p(x)$ of~$x$ through~$i$. 

Two $n$-cells~$x$ and~$y$ of~$X$ are called \emph{parallel} if either~$n=0$, or $s(x) = s(y)$ and $t(x)=t(y)$. An \emph{$n$-sphere of~$X$} is a pair $\gamma=(x,y)$ of parallel $n$-cells of~$X$, in which case we call~$x$ the \emph{source of~$\gamma$} and~$y$ the \emph{target of~$\gamma$}. 

\subsubsection{Composable cells}

Let~$X$ be a globular object of~$\Cr$, and fix~$k\geq 0$. Define the \emph{$k$-source map of~$X$} and the \emph{$k$-target map of~$X$} as the indexed morphisms
\[
X \ofl{s_k} X_k
\qquad\text{and}\qquad
X \ofl{t_k} X_k
\]
of degree~$0$, given, on an $n$-cell~$x$ of~$X$, by
\[
s_k(x) = 
\begin{cases}
s^{n-k} (x) &\text{if $n\geq k$} \\
i^{n-k} (x) &\text{if $n\leq k$},
\end{cases}
\qquad\text{and}\qquad
t_k(x) = 
\begin{cases}
t^{n-k} (x) &\text{if $n\geq k$} \\
i^{n-k} (x) &\text{if $n\leq k$}.
\end{cases}
\]
The globular relations generalise, for~$j<k$, to
\[
s_j s_k = s_j t_k = s_j
\qquad\text{and}\qquad
t_j s_k = t_j t_k = t_j.
\]

For~$k\geq 0$, denote by~$X\star_k X$ the pullback
\[
\xymatrix @=1.5em {
X \star_k X
	\ar[r]
	\ar[d]
	\pullbackcorner
& X
	\ar[d] ^-{s_k}
\\
X
	\ar[r] _-{t_k}
& X_k
}
\]
in~$\Ind(\Cr)$. Explicitly, the $n$-cells of~$X\star_k X$ are the pairs~$(x,y)$ of $n$-cells of~$X$ such that~$t_k(x)=s_k(y)$ holds. An $n$-cell~$(x,y)$ of~$X\star_k X$ is called a \emph{$k$-composable pair of~$n$-cells of~$X$}. Note that, by definition of~$s_k$ and~$t_k$, if~$n\leq k$, then the $n$-cells of~$X\star_k X$ are all pairs~$(x,x)$ for~$x$ an $n$-cell of~$X$.

\subsubsection{Internal \pdf{\infty}-categories}

An \emph{$\infty$-category of~$\Cr$} is a globular object~$X$ of~$\Cr$ equipped, for every~$k\geq 0$, with an indexed morphism
\[
X \star_k X \ofl{c_k} X
\]
of degree~$0$, called the \emph{$k$-composition of~$X$}, whose value at~$(x,y)$ is denoted by $x\star_k y$, and such that the following relations are satisfied for all $0\leq k<n$:
\begin{enumerate}
\item (\emph{compatibility with the source and target maps}) for every $n$-cell~$(x, y)$ of~$V \star_k V$,
\[
s(x\star_k y) = 
\begin{cases} 
s(x) 
	&\text{if $k=n-1$} \\ 
s(x)\star_k s(y) 
	&\text{otherwise} 
\end{cases}
\qquad\text{and}\qquad
t(x\star_k y) = 
\begin{cases} 
t(y) 
	&\text{if $k=n-1$} \\ 
t(x)\star_k t(y) 
	&\text{otherwise},
\end{cases}
\]
\item (\emph{compatibility with the identity map}) for every $n$-cell~$(x, y)$ of~$V\star_k V$, 
\[
1_{x\star_k y} = 1_x \star_k 1_y,
\]
\item (\emph{associativity}) for all $n$-cells~$x$, $y$ and~$z$ of~$V$ such that~$(x, y)$ and~$(y, z)$ are $n$-cells of~$V\star_k V$, 
\[
(x\star_k y) \star_k z = x \star_k (y \star_k z),
\]
\item (\emph{neutrality}) for every $n$-cell~$x$ of~$V$, 
\[
s_k(x) \star_k x = x = x \star_k t_k(x),
\]
\item (\emph{exchange}) for every~$j<k$, and all $n$-cells~$(x, x')$ and~$(y, y')$ of~$V\star_k V$ such that~$(x, y)$ and~$(x', y')$ are $n$-cells of~$V\star_j V$, 
\[
(x\star_k x') \star_j (y \star_k y') = (x\star_j y) \star_k (x'\star_j y').
\]
\end{enumerate}

Note that the compatibility of the compositions with the source and target maps ensures that the associativity axiom makes sense: if~$(x,y)$ and~$(y,z)$ are $n$-cells of~$V\star_k V$, then so do $(x\star_k y, z)$ and $(x, y\star_k z)$. The compatibility of compositions with identities implies that we can still write~$x$ for~$1_x$ with no ambiguity.

Given $\infty$-categories~$X$ and~$Y$ of~$\Cr$, an \emph{$\infty$-functor from~$X$ to~$Y$} is a globular morphism~$F:X\fl Y$ that commutes with all $k$-compositions:
\[
\xymatrix @=1.5em {
{X\star_k X}
	\ar [r] ^-{c_k}
	\ar [d] _-{F \times F}
	\ar@{} [dr] |-{\sm =}
& {X}
	\ar [d] ^-{F}
\\
{Y\star_k Y}
	\ar [r] _-{c_k}
& {Y}
}
\]

We denote by $\infty\Cat(\Cr)$ the category of $\infty$-categories and $\infty$-functors of~$\Cr$.
For $n\geq 0$, an \emph{$n$-category of~$\Cr$} is defined like an $\infty$-category of~$\Cr$, but starting with an $n$-globular object of~$\Cr$, equipped with compositions~$c_0$, \dots, $c_{n-1}$. 
We denote by $n\Cat(\Cr)$ the category of $n$-categories of~$\Cr$ and the corresponding $n$-functors.

\subsubsection{Internal \pdf{\infty}-groupoids}

In an $\infty$-category~$X$ of~$\Cr$, for~$n\geq 1$, an $n$-cell~$x$ is called \emph{invertible} if there exists an $n$-cell~$x^-$ in~$X$, of source~$t(x)$ and target~$s(x)$, such that the relations
\[
x\star_{n-1} x^- = s(x)
\qquad\text{and}\qquad
x^-\star_{n-1} x = t(x)
\]
are satisfied. An \emph{$\infty$-groupoid of~$\Cr$} is an $\infty$-category of~$\Cr$ in which all $k$-cells are invertible, for every~$k\geq 1$. 
Similarly, for~$n\geq 0$, an \emph{$n$-groupoid of~$\Cr$} is an $n$-category in which all $k$-cells are invertible, for every~$k \geq 1$. 
We denote by $\infty\Gpd(\Cr)$ the category of $\infty$-groupoids of~$\Cr$ and of $\infty$-functors between them, and by $n\Gpd(\Cr)$ its full subcategory whose objects are $n$-groupoids.

\subsection{Higher-dimensional vector spaces}
\label{SS:nVect}

\subsubsection{Globular vector spaces}

The objects and morphisms of the category $\Glob(\Vect)$ are called \emph{globular vector spaces} and \emph{globular linear maps}. Explicitly, a globular vector space is a diagram 
\[
\xymatrix @!C @C=3em @W=2em {
V_0 
	\ar [r] |-*+{i}
& V_1 
	\ar@<-1.5ex> [l] _-{s}
	\ar@<1.5ex> [l] ^-{t}
	\ar [r] |-*+{i}
& \cdots
	\ar@<-1.5ex> [l] _-{s}
	\ar@<1.5ex> [l] ^-{t}
	\ar [r] |-*+{i}
& V_n
	\ar@<-1.5ex> [l] _-{s}
	\ar@<1.5ex> [l] ^-{t}
	\ar [r] |-*+{i}
& V_{n+1}
	\ar@<-1.5ex> [l] _-{s}
	\ar@<1.5ex> [l] ^-{t}
	\ar [r] |-*+{i}
& \cdots
	\ar@<-1.5ex> [l] _-{s}
	\ar@<1.5ex> [l] ^-{t}
}
\]
of vector spaces and linear maps that satisfy the globular relations. 
If~$V$ is a globular vector space, and if~$a$ is an $n$-cell of~$V$, for~$n\geq 1$, then the \emph{boundary of~$a$} is the $(n-1)$-cell of~$V$ denoted by ~$\dr(a)$ and defined by
\[
\dr(a) = s(a) - t(a).
\]

\subsubsection{Higher-dimensional vector spaces}

The categories $\infty\Cat(\Vect)$ and $n\Cat(\Vect)$ are denoted by $\infty\Vect$ and $n\Vect$, and their objects and morphisms are called \emph{$\infty$-vector spaces}, \emph{linear $\infty$-functors}, \emph{$n$-vector spaces} and \emph{linear $n$-functors}, respectively. In particular, we have $0\Vect=\Vect$, and $1\Vect$ is the category of $2$-vector spaces of~\cite{BaezCrans04}.

If~$V$ is a globular vector space, the pullback~$V\star_k V$ is equipped with the vector space structure given by, for all $n$-cells~$(a,a')$ and~$(b,b')$ of~$V\star_k V$ and all scalars~$\lambda$ and~$\mu$,
\[
\lambda (a, a') + \mu (b, b') = (\lambda a + \mu b, \lambda a' + \mu b').
\]
Hence, if~$V$ is an $\infty$-vector space, the linearity of its composition~$c_k$ is equivalent, in the same context, to
\begin{equation}
\label{E:CompLinear}
(\lambda a + \mu b) \star_k (\lambda a'+ \mu b') \:=\: \lambda (a\star_k a') + \mu (b\star_k b').
\end{equation}

The following result states that globular vector spaces, $\infty$-vector spaces and $\infty$-groupoids of~$\Vect$ are the same notions. This also holds for $n$-globular vector spaces, $n$-vector spaces and $n$-groupoids of~$\Vect$, which is essentially~\cite[Prop.~2.5]{KhudaverdyanMandalPoncin11}. 

\begin{proposition}
\label{P:nVect}
The forgetful functors
\[
\infty\Gpd(\Vect) \fll \infty\Vect \fll \Glob(\Vect)
\]
are isomorphisms. In particular, a globular vector space~$V$ can be uniquely extended into an $\infty$-vector space, by putting
\begin{equation}
\label{E:nVectComp}
a \star_k b = a - t_k (a) + b,
\end{equation}
and, in an $\infty$-vector space, every $n$-cell~$a$ is invertible, with inverse
\begin{equation}
\label{E:nVectInv}
a^- = s(a) - a + t(a).
\end{equation}
\end{proposition}

\begin{proof}
Let~$V$ be an $\infty$-vector space. For~$(a,b)$ in~$V_n\star_k V_n$, we have
\[
a\star_k b = (a - s_k(b) + s_k(b)) \star_k (t_k(a) - t_k(a) + b),
\]
so that, using $t_k(a)=s_k(b)$, the linearity of the $k$-composition and the neutrality axioms, we obtain~\eqref{E:nVectComp}:
\[
a\star_k b = (a \star_k t_k(a)) - (s_k(b) \star_k t_k(a)) + (s_k(b) \star_k b) = a - t_k(a) + b.
\]
As a consequence, given a globular vector space~$V$, there exists at most one $\infty$-vector space with~$V$ as underlying vector space, whose compositions are given by~\eqref{E:nVectComp}. 

Conversely, let~$V$ be a globular vector space and define the $k$-composition~$c_k$ by~\eqref{E:nVectComp}. Let us check that the axioms of an $\infty$-vector space are fulfilled. First, $c_k$ is linear, which is obtained by replacing $k$-compositions by their definitions in both sides of the relations~\eqref{E:CompLinear} and by using the linearity of~$t$.
Next, $c_k$ is compatible with the source map:
\[
s (a\star_k b) 
	= s(a) - s(1_{s_k(b)}) + s(b) 
	= s(a) - s_k(b) + s(b)
	=
\begin{cases}
s(a) &\text{if $k=n-1$,} \\
s(a) \star_k s(b) &\text{otherwise.}
\end{cases}
\]
We proceed symmetrically to get the compatibility of~$c_k$ with the target map, and its compatibility with the identity map comes from the linearity of the latter. For associativity, we use the compatibility of~$c_k$ with the source and target maps to get, by induction on~$n\geq k$, that $t_k(a\star_k b)=t_k(b)$ and $s_k(b\star_k c)=s_k(b)$, and, then, we replace $k$-compositions in both sides of the associativity axiom to get the result. The neutrality relations are immediate consequences of the globular relations $si=ti=\id_{V}$. Finally, the exchange relations are obtained by using~$t_jt_k=t_j$ when~$j<k$ to prove that both sides are equal to
\[
a - t_k(a) + a' - t_j(a) + b - t_k(b) + b'.
\]
We conclude by observing that any globular linear map~$F$ is automatically a linear $\infty$-functor:
\[
F(a\star_k b) = F(a - t_k(a) + b) = F(a) - t_k(F(a)) + F(b) = F(a) \star_k F(b).
\]

Now, let~$V$ be an $\infty$-vector space and~$a$ be an $n$-cell of~$V$. With~$a^-$ defined by~\eqref{E:nVectInv}, we check $s(a^-)=t(a)$ and $t(a^-)=s(a)$ and, then, we apply~\eqref{E:nVectComp} to obtain $a\star_{n-1} a^-=s(a)$ and $a^-\star_{n-1} a = t(a)$.
\end{proof}

\subsubsection{The graded case}

Replacing~$\Vect$ by the category~$\Vectg$ of (non-negatively) graded vector spaces over~$\K$ and graded linear maps of degree~$0$, one obtains the category $\Glob(\Vectg)$ of \emph{globular graded vector spaces}: these are globular vector spaces~$V$ such that each vector space~$V_n$ is graded, i.e.\ each~$V_n$ admits a decomposition $V_n = \bigoplus_{i\in\Nb} V_n^{(i)}$, and the components of the source, target and identity maps are graded linear maps of degree~$0$. In such a graded globular vector space~$V$, the $n$-cells of~$V_n^{(i)}$ are said to be \emph{homogeneous of degree~$i$}. Similarly, the categories $\infty\Vectg$ and $\infty\Gpd(\Vectg)$ are obtained by replacing of~$\Vect$ by~$\Vectg$, and Proposition~\ref{P:nVect} extends in a straightforward way.

\subsection{Higher-dimensional associative algebras}
\label{SS:nAlg}

\subsubsection{Globular associative algebras and bimodules} 
The objects and morphisms of the category $\Glob(\Alg)$ are called \emph{globular algebras} and \emph{morphisms of globular algebras}. 
If~$A$ is an algebra, and $\Bimod(A)$ is the category of $A$-bimodules, a \emph{globular $A$-bimodule} is an object of $\Glob(\Bimod(A))$. In view of Theorem~\ref{T:nAlg}, the category $\Glob(\Bimod)$ of \emph{globular bimodules} is the one whose objects are pairs~$(A,M)$ formed by an algebra~$A$ and a globular $A$-bimodule~$M$, and whose morphisms from~$(A,M)$ to~$(B,N)$ are pairs~$(F,G)$ made of a morphism $F:A\fl B$ of algebras and a morphism $G:M\fl N$ of bimodules (in the sense that $G(ama')=F(a)G(m)F(a')$ holds for all~$a$ and~$a'$ in~$A$ and~$m$ in~$M$).

\subsubsection{Higher-dimensional associative algebras}
\label{SSS:HigherDimensionalAssociativeAlgebras}
The categories $\infty\Cat(\Alg)$ and $n\Cat(\Alg)$ are denoted by $\infty\Alg$ and $n\Alg$, and their objects and morphisms are called \emph{$\infty$-algebras}, \emph{morphisms of $\infty$-algebras}, \emph{$n$-algebras} and \emph{morphisms of $n$-algebras}, respectively. For a globular algebra~$A$, the product of $n$-cells~$(a,a')$ and~$(b,b')$ of~$A\star_k A$ is given by
\[
(a, a')(b, b') = (ab, a'b').
\]
This implies that, if~$A$ is an $\infty$-algebra, the fact that the composition~$c_k$ is a morphism of algebras is equivalent, with the same notations, to 
\begin{equation}
\label{E:nAlgComp}
ab \star_k a'b' = (a\star_k a')(b\star_k b').
\end{equation}

The following result states that the structure of $\infty$-algebra boils down to the one of a globular bimodule that satisfies an extra condition, corresponding to the fact that the composition~$c_0$ of an $\infty$-algebra satisfies~\eqref{E:nAlgComp}. This also holds for $n$-algebras and $n$-globular bimodules with the same condition.

\begin{theorem}
\label{T:nAlg}
The following categories are isomorphic:
\begin{enumerate}
\item The category~$\infty\Alg$ of~$\infty$-algebras.
\item The category $\infty\Gpd(\Alg)$ of internal $\infty$-groupoids in~$\Alg$.
\item The full subcategory of $\Glob(\Alg)$ whose objects are the globular algebras~$A$ that satisfy, for all $n$-cells~$a$ and~$b$ of~$A$, the relations 
\begin{equation}
\label{E:nAlgExchangeProd}
ab = a s_0(b) + t_0(a) b - t_0(a) s_0(b) = s_0(a) b + a t_0(b) - s_0(a) t_0(b).
\end{equation}
\item The full subcategory of $\Glob(\Bimod)$ whose objects are the pairs~$(A,M)$ such that~$M_0$  is equal to~$A$, with its canonical $A$-bimodule structure, and that satisfy, for all $n$-cells~$a$ and~$b$ of~$M$, the relation
\begin{equation}
\label{E:nAlgExchange}
a s_0(b) + t_0(a) b - t_0(a) s_0(b) = s_0(a) b + a t_0(b) - s_0(a) t_0(b).
\end{equation}
\end{enumerate}
\end{theorem}

\begin{proof}
First, one checks that~\eqref{E:nVectInv} defines an inverse for every $n$-cell of an $\infty$-algebra, for~$n\geq 1$, so that~$\infty\Alg$ is isomorphic to $\infty\Gpd(\Alg)$. 

Now, let~$A$ be an $\infty$-algebra, and let us prove that~\eqref{E:nAlgExchangeProd} holds. For~$0$-cells~$a$ and~$b$ of~$A$, \eqref{E:nAlgExchangeProd} reads $ab = ab = ab$. If~$a$ and~$b$ are $n$-cells of~$A$, for~$n\geq 1$, Proposition~\ref{P:nVect} and~\eqref{E:nAlgComp} imply
\[
ab = (a \star_0 t_0(a)) (s_0(b) \star_0 b) = a s_0(b) \star_0 t_0(a) b = a s_0(b) - t_0(a) s_0(b) + t_0(a) b,
\]
and symmetrically for the other part of~\eqref{E:nAlgExchangeProd}. Moreover, forgetting the compatibility with the compositions of a morphism of $\infty$-algebras gives a morphism of globular algebras.

Conversely, let~$A$ be a globular algebra that satisfies~\eqref{E:nAlgExchangeProd}. Let us check that the unique possible composition~$\star_k$ on~$A$, given by~\eqref{E:nVectComp}, satisfies~\eqref{E:nAlgComp}. We start with~$k=0$, and pairs~$(a, a')$ and~$(b, b')$ in $A_n\star_0 A_n$. Writing~$c$ for~$t_0(a)$ and~$d$ for~$t_0(b)$, we obtain
\[
(a\star_0 a')(b\star_0 b') = ab + ab' + a'b + a'b' - ad - a'd - cb - cb' + cd.
\]
We use~\eqref{E:nAlgExchangeProd} on~$ab'$ and~$a'b$ to get
\[
ab' = ag + cb' - cd
\qquad\text{and}\qquad
a'b = cb + a'd - cd,
\]
so that we conclude
\[
(a\star_0 a')(b\star_0 b') = ab + a'b' - cd = ab \star_0 a'b'.
\]
Now, let us fix~$k\geq 1$, and pairs~$(a, a')$ and~$(b, b')$ in $A_n\star_k A_n$. We write~$c$ for~$t_0(a)$ and~$d$ for~$s_0(b)$, and we note that $c=t_0(a)=t_0t_k(a)=t_0s_k(a')=t_0(a')$ and, similarly, that $d=s_0(b')$. Then, we use the fact that~$\star_0$ satisfies~\eqref{E:nAlgComp} and the exchange relation between~$\star_0$ and~$\star_k$ to get
\[
ab \star_k a'b' = (ad \star_0 cb) \star_k (a'd \star_0 cb') = (ad \star_k a'd) \star_0 (cb \star_k cb').
\]
By definition of~$\star_k$, we obtain
\[
ad \star_k a'd = ad - t_k(a)d + a'd = (a\star_k a') d
\]
and
\[
cb \star_k cb' = cb - c t_k(b) + cb' = c (b\star_k b').
\]
Hence we get, using~\eqref{E:nAlgComp} for~$\star_0$ again,
\[
ab \star_k a'b' = (a\star_k a') d \star_0 c (b\star_k b') = (a\star_k a') (b\star_k b').
\]
Moreover, if~$F$ is a morphism between two globular algebras that satisfy~\eqref{E:nAlgComp}, then~\eqref{E:nVectComp} and the linearity of~$F$ imply
\[
F(x\star_k y) = F(x) - t_k(F(x)) + F(y) = F(x)\star_k F(y),
\]
which proves that~$F$ is a morphism between the corresponding $\infty$-algebras. This concludes the proof that~$\infty\Alg$ is isomorphic to the category of~\Item{iii}.

Let~$A$ be a globular algebra that satisfies~\eqref{E:nAlgExchangeProd}. By hypothesis, $A_0$ is an algebra and each~$A_n$, for~$n\geq 1$, is equipped with a structure of $A_0$-bimodule by using its algebra product with iterated identities of elements of~$A_0$. Moreover, the source, target and identity maps commute with these $A_0$-bimodule structures because they are morphisms of algebras, and~\eqref{E:nAlgExchange} is satisfied by hypothesis. Finally, any morphism of globular algebras induces a morphism between the underlying globular bimodules.

Conversely, let~$A$ be an algebra, $M$ be an $A$-bimodule with $M_0=A$ and that satisfies~\eqref{E:nAlgExchange}. On each~$M_n$, for~$n\geq 1$, we define the product~$ab$ by
\[
ab = a s_0(b) + t_0(a) b - t_0(a) s_0(b),
\]
and we check that this equips~$M$ with a structure of a globular algebra that satisfies the relation~\eqref{E:nAlgExchangeProd}. Finally, a morphism of globular bimodules that satisfies~\eqref{E:nAlgExchange} commutes with the product defined by~\eqref{E:nAlgExchangeProd}. Thus~$\infty\Alg$ is indeed isomorphic to the category described in~\Item{iv}.
\end{proof}

\subsubsection{The graded case}

As in the case of vector spaces, we obtain the category $\Glob(\Algg)$ of \emph{globular graded algebras} by replacing~$\Alg$ with the category~$\Algg$ of graded algebras: in such an object~$A$, the product of each~$A_n$ is graded so that, if~$f$ and~$g$ are $n$-cells of~$A$, with~$f$ homogeneous of degree~$i$ and~$g$ homogeneous of degree~$j$, then their product~$fg$ is a homogeneous $n$-cell of~$A$ of degree $i+j$. A globular graded algebra~$A$ is called \emph{connected} if it satisfies $A_n^{(0)}=\K$ for every~$n\geq 0$. The notions of graded $\infty$-algebra and globular graded bimodule are obtained similarly, and Theorem~\ref{T:nAlg} extends to this context.

\section{Polygraphs for associative algebras}
\label{S:nPolygraphs}

In~\cite{Street76, Street87,Burroni93}, Street and Burroni have introduced categorical objects known as computads or polygraphs, to describe generating families and presentations of higher-dimensional categories. Here, we adapt these set-theoretic objects to provide bases and presentations of higher-dimensional algebras. As for the original polygraphs, which can be used to define homotopical resolutions of monoids, categories and $n$-categories, the polygraphs we define here give a notion of polygraphic resolution of an associative algebra.

\subsection{Extended higher-dimensional associative algebras}
\label{SS:ExtendedNAlg}

\subsubsection{Cellular extensions}

Fix a natural number~$n$, and let~$A$ be an $n$-algebra. A \emph{cellular extension of~$A$} is a set~$X$ equipped with maps
\[
\xymatrix @=1.5em {
A_n 
& X
	\ar@<-0.5ex> [l] _-{s}
	\ar@<0.5ex> [l] ^-{t}
}
\]
such that, for every~$x$ in~$X$, the pair $(s(x),t(x))$ is an $n$-sphere of~$A$. 

Let~$X$ be a cellular extension of~$A$. Denote by~$\equiv_X$ the congruence relation on the parallel $n$-cells of~$A$ generated by $s(x)\sim t(x)$ for every~$x$ in~$X$ (that is, $\equiv_X$ is the smallest equivalence relation on the parallel $n$-cells of~$A$, compatible with all the compositions of~$A$ and relating~$s(x)$ and~$t(x)$ for every~$x$ in~$X$).
Call~$X$ \emph{acyclic} if, for every $n$-sphere~$(a,b)$ of~$A$, we have $a\equiv_X b$. Every $n$-algebra~$A$ has two canonical cellular extensions: the empty one, and the one denoted by $\Sph(A)$ that contains all the $n$-spheres of~$A$ and is, as a consequence, acyclic.

\subsubsection{Extended higher-dimensional algebras}

For~$n\geq 0$, the category $n\Alg^+$ of \emph{extended $n$-algebras} is defined as the pullback
\[
\xymatrix @=1.5em{
n\Alg^+
 	\ar[r]
	\ar[d] 	
	\pullbackcorner
&
(n+1)\Gph
  \ar[d] 
\\
n\Alg
  \ar[r] 
&
n\Gph
}
\]
of the forgetful functors $n\Alg\fl n\Gph$ (forgetting the algebra structures) and $(n+1)\Gph \fl n\Gph$ (forgetting the $(n+1)$-cells), where~$n\Gph$ is the category of set-theoretic $n$-graphs (or $n$-semiglobular objects in~$\Set$).
In an explicit way, an object of~$n\Alg^+$ is a pair~$(A,X)$ formed by an $n$-algebra~$A$ and a cellular extension~$X$ of~$A$, while a morphism from~$(A,X)$ to~$(B,Y)$ is a pair~$(F,\phi)$ made of a morphism of $n$-algebras $F:A\fl B$, and a map $\phi:X\fl Y$ that commutes with the source and target maps.

\subsubsection{Free higher-dimensional algebras}
\label{SSS:FreeHigherAlgebras}
Fix~$n\geq 1$. The forgetful functor from~$n\Alg$ to~$(n-1)\Alg^+$, discarding the compositions of $n$-cells,  admits a left adjoint, that maps an extended $(n-1)$-algebra~$(A,X)$ to the $n$-algebra denoted by~$A[X]$, called the \emph{free $n$-algebra over~$(A,X)$}, and defined as follows.
First, we consider the $A_0$-bimodule
\[
M = \big( A_0 \tens \K X \tens A_0 \big) \oplus A_{n-1}
\]
obtained by the direct sum of the free $A_0$-bimodule with basis~$X$ and of a copy of~$A_{n-1}$, equipped with its canonical $A_0$-bimodule structure. Thus~$M$ contains linear combinations of elements $axb$, for~$a$ and~$b$ in~$A_0$ and~$x$ in~$X$, and of an $(n-1)$-cell~$c$ of~$A$. We define the source, target and identity maps between~$M$ and~$A_{n-1}$ by
\[
s(axb) = a s(x) b,
\qquad
s(c) = c, 
\qquad
t(axb) = a t(x) b,
\qquad
t(c) = c
\qquad\text{and}\qquad
i(c) = c, 
\]
for all~$x$ in~$X$, $a$ and~$v$ in~$A_0$, and~$c$ in~$A_{n-1}$. 
Then we define the $A_0$-bimodule~$A[X]_n$ as the quotient of~$M$ by the $A_0$-bimodule ideal generated by all the elements
\[
\big( a s_0(b) + t_0(a) b - t_0(a) s_0(b) \big) 
	\: - \: 
	\big( s_0(a) b + a t_0(b) - s_0(a) t_0(b) \big),
\] 
where~$a$ and~$b$ range over $A_0 \tens \K X \tens A_0$. We check that the source and target maps are compatible with the quotient, so that, by Theorem~\ref{T:nAlg}, the $A_0$-bimodule $A[X]_n$ extends~$A$ into a uniquely defined $n$-algebra~$A[X]$. 
Note that, by construction of~$A[X]$, the cellular extension~$X$ is acyclic if, and only if, for every $(n-1)$-sphere~$(a,b)$ of~$A$, there exists an $n$-cell $f:a\fl b$ in the free $n$-algebra~$A[X]$.

\subsubsection{Quotient higher-dimensional algebras}

For~$n\geq 0$ and~$(A,X)$ an extended $n$-algebra, the \emph{quotient of~$A$ by~$X$} is the $n$-algebra denoted by~$A/X$, and obtained by quotient of the $n$-cells of~$A$ by the congruence~$\equiv_X$. 
Thus, two $n$-cells~$a$ and~$b$ of~$A$ are identified in~$A/X$ if, and only if, there exists an $(n+1)$-cell $f:a\fl b$ in the free $(n+1)$-algebra~$A[X]$. As a consequence, the cellular extension~$X$ is acyclic if, and only if, the canonical projection $A\pfl A/X$ identifies all the parallel $n$-cells of~$A$.

\subsection{Polygraphs for associative algebras}
\label{SS:Polygraphs}

\subsubsection{Polygraphs of finite dimension}
\label{SSS:PolygraphFiniteDimension}
The category $n\Pol(\Alg)$ of \emph{$n$-polygraphs for associative algebras}, simply called \emph{$n$-polygraphs} in this paper, and the \emph{free $n$-algebra functor} from $n\Pol(\Alg)$ to $n\Alg$ are defined by induction on~$n\geq 0$ as follows. 

For~$n=0$, we define $0\Pol(\Alg)$ as the category of sets. If~$X$ is a set, the \emph{free $0$-algebra functor} maps~$X$ to the free algebra~$\lin{X}$ over~$X$.
Assume that~$n\geq 1$ and that the category $(n-1)\Pol(\Alg)$ of $(n-1)$-polygraphs is defined, together with the free $(n-1)$-algebra functor $(n-1)\Pol(\Alg) \fl (n-1)\Alg$. The category~$n\Pol(\Alg)$ of $n$-polygraphs is defined as the pullback
\begin{equation}
\label{E:DefLinPol}
\vcenter{\xymatrix @=1.5em{
n\Pol(\Alg)
  \ar[r]
  \ar[d] \pullbackcorner
&
(n-1)\Alg^+
  \ar[d] 
\\
(n-1)\Pol(\Alg)
  \ar[r] 
&
(n-1)\Alg
}}
\end{equation}
of the free $(n-1)$-algebra functor and the functor forgetting the cellular extension. The free $n$-algebra functor is obtained as the composite
\[
n\Pol(\Alg) \fll (n-1)\Alg^+ \fll n\Alg
\]
of the functor $n\Pol(\Alg) \fll (n-1)\Alg^+$ of~\eqref{E:DefLinPol}, followed by the functor mapping~$(A,X)$ to~$A[X]$. 

If~$X$ is an $n$-polygraph, we denote by~$\lin{X}$ the free $n$-algebra over~$X$. Expanding the definition, an $n$-polygraph is a sequence $(X_0,\dots,X_n)$,  written~$\cohpres{X_0}{\cdots}{X_n}$, made of a set~$X_0$ and, for every $0\leq k<n$, a cellular extension~$X_{k+1}$ of the free $k$-algebra over $\cohpres{X_0}{\cdots}{X_k}$. The free $n$-algebra over~$X$ is given by
\[
\lin{X} = \lin{X_0}[X_1] \cdots [X_n].
\]

\subsubsection{Polygraphs of infinite dimension}

The category~$\infty\Pol(\Alg)$ of \emph{$\infty$-polygraphs} and the corresponding \emph{free $\infty$-algebra functor} are obtained as the limit of the functors $(n+1)\Pol(\Alg)\fl n\Pol(\Alg)$ of~\eqref{E:DefLinPol}. Thus, by construction, an $\infty$-polygraph $X$ is a sequence $\cohpres{X_0}{\cdots}{X_n \,\vert\, \cdots}$ such that $\cohpres{X_0}{\cdots}{X_n}$ is an $n$-polygraph for every~$n\geq 0$, and $n$-polygraphs are the $\infty$-polygraphs~$X$ such that $X_p=\emptyset$ for~$p>n$.

Let~$X$ be an $\infty$-polygraph. The elements of~$X_n$ are called the \emph{$n$-cells of~$X$}. We commit the abuse to also denote by~$X_n$ the underlying $n$-polygraph of~$X$. We say that an $\infty$-polygraph is of \emph{finite type} if it has finitely many $n$-cells for every~$n\geq 0$.

\subsubsection{Higher-dimensional monomials}

Let~$X$ be an $\infty$-polygraph. A \emph{monomial of~$\lin{X}$} is an element of the free monoid~$X_0^*$ over~$X_0$, that is, a (possibly empty) product $x_1\cdots x_p$ of elements of~$X_0$. The monomials of~$\lin{X}$ form a linear basis of the free algebra~$\lin{X_0}$, which means that every $0$-cell~$a$ of~$\lin{X}$ can be uniquely written as a (possibly empty) linear combination
\[
a = \sum_{i=1}^p \lambda_i u_i
\]
of pairwise distinct monomials~$u_1$, \dots, $u_p$ of~$\lin{X}$, with~$\lambda_1$, \dots, $\lambda_p$ nonzero scalars. This expression is called \emph{the canonical decomposition of~$a$} and we denote by~$\supp(a)$ the set~$\ens{u_1,\dots,u_p}$. 

For~$n\geq 1$, an \emph{$n$-monomial of~$\lin{X}$} is an $n$-cell of~$\lin{X}$ with shape~$u\alpha v$, where~$\alpha$ is an $n$-cell of~$X$, and~$u$ and~$v$ are monomials of~$\lin{X}$. By construction of the free $n$-algebra over~$(\lin{X_{n-1}},X_n)$, and by freeness of~$\lin{X_{n-1}}$, every $n$-cell~$a$ of~$\lin{X}$ can be written as a linear combination
\begin{equation}
\label{E:DecompositionNCell}
a = \sum_{i=1}^p \lambda_i a_i + 1_c
\end{equation}
of pairwise distinct $n$-monomials~$a_1$, \dots, $a_p$ and of an identity $n$-cell~$1_c$ of~$\lin{X}$, and this decomposition is unique up to the relations
\[
a s_0(b) + t_0(a) b - t_0(a) s_0(b) = s_0(a) b + a t_0(b) - s_0(a) t_0(b),
\]
where~$a$ and~$b$ range over the $n$-monomials of~$\lin{X}$. 

If~$a$ is an $n$-cell of~$\lin{X}$, the \emph{size of~$a$} is the minimum number of $n$-monomials of~$\lin{X}$ required to write~$a$ as in~\eqref{E:DecompositionNCell}, and we denote by~$\cell(a)$ the subset of~$X_n$ that consists of the $n$-cells of~$X$ that appear in the corresponding $n$-monomials. 

\begin{lemma}
\label{L:NCellDecomposition}
Let~$X$ be an $\infty$-polygraph, and fix~$n\geq 1$. Then, every nonidentity $n$-cell~$a$ of~$\lin{X}$ admits a decomposition
\[
a = a_1 \star_{n-1} \cdots \star_{n-1} a_p,
\]
where~$a_1$, \dots, $a_n$ are $n$-cells of size~$1$ in~$\lin{X}$.
\end{lemma}

\begin{proof}
Write~$a=\lambda_1 b_1 + \cdots + \lambda_p b_p + 1_c$ as in~\eqref{E:DecompositionNCell}.
If~$p=1$, then~$a$ is of size~$1$. Otherwise, put
\[
d_i = \sum_{j=1}^i \lambda_j t(b_j)
\qquad\text{and}\qquad
e_i = \sum_{j=i}^p \lambda_j s(b_j)
\]
for each~$i$ in~$\ens{1,\dots,p}$, and $d_0 = e_{p+1} =0$. Define, for each~$i$ in~$\ens{1,\dots,p}$, the $n$-cell of size~$1$
\[
a_i = \lambda_i b_i + 1_c + 1_{d_{i-1}} + 1_{e_{i+1}}.
\]
First, check $s(a_i) = c + d_{i-1} + e_i$ and $t(a_i) = c + d_i + e_{i+1}$, so that $a_1\star_{n-1}\cdots\star_{n-1}a_p$ is a well-defined $n$-cell of~$\lin{X}$. Then, \eqref{E:nVectComp} gives 
\[
a_1 \star_{n-1} \cdots \star_{n-1} a_p 
	= \sum_{i=1}^p \lambda_i b_i 
	+ \sum_{i=1}^p (1_c + 1_{d_{i-1}} + 1_{e_{i+1}}) 
	- \sum_{i=1}^{p-1} (\lambda_i 1_{t(b_i)} + 1_c + 1_{d_{i-1}} + 1_{e_{i+1}}).
\]
We conclude thanks to
\[
d_{p-1} = \sum_{i=1}^{p-1} \lambda_i t(b_i)
\qquad\text{and}\qquad
e_{p+1}=0.
\qedhere
\]
\end{proof}

\subsubsection{Polygraphs for graded associative algebras}
\label{SS:GradedLinearPolygraphs}

The category $n\Algg^+$ of \emph{graded extended $n$-algebras} is defined similarly to $n\Alg^+$: its objects are pairs~$(A,X)$, where~$A$ is a graded $n$-algebra, and~$X$ is a graded cellular extension of~$A$, meaning that $X = \amalg_{i\in\Nb} X^{(i)}$ and that the source and target of each~$x$ in~$X^{(i)}$ are homogeneous of degree~$i$. In that case, the free $(n+1)$-algebra~$A[X]$ and the quotient $n$-algebra~$A/X$, defined as in the nongraded case, are also graded.

A \emph{graded $\infty$-polygraph} is an $\infty$-polygraph~$X$ such that each set~$X_n$ is graded, for $n\geq 0$. This notion restricts to $n$-polygraphs, and, if~$N\geq 2$, a $1$-polygraph~$X$ is called \emph{$N$-homogeneous} if~$X_0$ is concentrated in degree~$1$ and~$X_1$ is concentrated in degree~$N$. We say \emph{quadratic} and \emph{cubical} instead of $2$-homogeneous and $3$-homogeneous, respectively. 
Given a fixed map $\omega : \Nb\setminus\ens{0} \fl \Nb\setminus\ens{0}$, we say that a graded $\infty$-polygraph~$X$ is \emph{$\omega$-concentrated} if each graded set~$X_n$, for~$n\geq 0$, is concentrated in degree~$\omega(n+1)$. In that case, because the source and target maps are graded, for~$n\geq 1$, the source and target of every $n$-cell of~$X$ are homogeneous $(n-1)$-cells of~$\lin{X}$ of degree~$\omega(n+1)$. 

\subsection{Presentations, coherent presentations and polygraphic resolutions}
\label{SS:PresentationsResolutions}

\subsubsection{Polygraphic presentations}

Let~$X$ be a $1$-polygraph. The \emph{algebra presented by~$X$} is the quotient algebra
\[
\cl{X} = \lin{X_0} / X_1.
\]
When the context is clear, we denote by~$\cl{a}$ the image of a $0$-cell~$a$ of~$\lin{X}$ through the canonical projection. 

Let~$A$ be an algebra. We say that~$A$ is \emph{presented by~$X$}, or that~$X$ is a \emph{presentation of~$A$}, if~$A$ is isomorphic to~$\cl{X}$. For example, if~$X_0$ is a generating set of~$A$, and if $X_1$ is the cellular extension of the free algebra over~$X_0$ that contains a $1$-cell $\alpha:a\fl b$ for all~$0$-cells~$a$ and~$b$ of~$\lin{X_0}$ that are equal in~$A$, then the $1$-polygraph $\pres{X_0}{X_1}$ is a presentation of~$A$.
We say that~$A$ is \emph{$N$-homogeneous} if it admits a presentation by an $N$-homogeneous graded $1$-polygraph. 

\subsubsection{Left-monomial \pdf{1}-polygraphs}
\label{SSS:LeftMonomialPolygraph}

Let~$X$ be a $1$-polygraph. We say that~$X$ is \emph{left-monomial} if, for every $1$-cell~$\alpha$ of~$X$, the source of~$\alpha$ is a monomial of~$\lin{X}$ that does not belong to~$\supp(t(\alpha))$. In that case, the source of any $1$-monomial of~$\lin{X}$ is a monomial of~$\lin{X}$. Note that, from~$X$, one obtains a left-monomial $1$-polygraph that presents the same algebra as~$X$ as follows. For every $1$-cell~$\alpha$ of~$X$, consider the boundary $\dr(\alpha) = s(\alpha) - t(\alpha)$ of~$\alpha$: if it is~$0$, discard~$\alpha$, otherwise, replace~$\alpha$ with 
\[
u \:\ofl{\alpha'}\: u - \frac{1}{\lambda} \dr(\alpha),
\]
where~$u$ is any chosen monomial in $\supp(\dr(\alpha))$ and~$\lambda$ is the coefficient of~$u$ in~$\dr(\alpha)$. 

\subsubsection{Presentations and ideals}

Let~$X$ be a $1$-polygraph. We denote by~$I(X)$ the ideal of the free algebra~$\lin{X_0}$ generated by the boundaries of the $1$-cells of~$X$. Since the algebra~$\lin{X_0}$ is free, the ideal~$I(X)$ is made of all the linear combinations
\[
\sum_{i=1}^p \lambda_i u_i \dr(\alpha_i) v_i,
\]
where $u_1\alpha_1 v_1$, \dots, $u_p\alpha_p v_p$ are pairwise distinct $1$-monomials of~$\lin{X}$, and~$\lambda_1$, \dots, $\lambda_p$ are nonzero scalars, so that the algebra~$\cl{X}$ presented by~$X$ is isomorphic to the quotient of~$\lin{X_0}$ by~$I(X)$.

\begin{lemma}
\label{L:Ideal}
Let~$X$ be a $1$-polygraph. For all $0$-cells~$a$ and~$b$ of~$\lin{X}$, the following are equivalent:
\begin{enumerate}
\item The $0$-cell~$a-b$ belongs to the ideal~$I(X)$.
\item There exists a $1$-cell $f:a\fl b$ in~$\lin{X}$.
\end{enumerate}
As a consequence, $I(X)$ exactly contains the $0$-cells~$a$ of~$\lin{X}$ such that~$\cl{a}=0$ holds in~$\cl{X}$.
\end{lemma}

\begin{proof}
Let us assume that~$a-b$ is in~$I(X)$, that is, 
$
a - b = \sum_{1\leq i\leq p} \lambda_i u_i\dr(\alpha_i) v_i.
$
Then the following $1$-cell~$f$ of $\lin{X}$ has source~$a$ and target~$b$:
\[
f \:=\: \sum_{i=1}^p \lambda_i u_i \alpha_i v_i \; + \; \big( a - \sum_{i=1}^p \lambda_i u_i s(\alpha_i) v_i \big).
\]

Conversely, let $f:a\fl b$ be a $1$-cell of~$\lin{X}$. Using Lemma~\ref{L:NCellDecomposition}, we decompose~$f$ into $1$-cells of size~$1$:
\[
f = f_1\star_0\cdots\star_0 f_p
\qquad\text{with}\qquad
f_i = \lambda_i u_i \alpha_i v_i + h_i.
\]
Since $t(f_i)=s(f_{i+1})$, we have $a - b = \dr(f_1) + \cdots + \dr(f_p)$. Moreover, $\dr(f_i) = \lambda_i u_i \dr(\alpha_i) v_i$ implies that each~$\dr(f_i)$ belongs to~$I(X)$, and thus so does~$a-b$.

Finally, if one applies the equivalence to the case~$b=0$, since~$\cl{0}=0$ holds in~$\cl{X}$, we get that~$a$ is in~$I(X)$ if, and only if, we have~$\cl{a}=0$ in~$\cl{X}$.
\end{proof}

\subsubsection{Coherent presentations and polygraphic resolutions}
\label{SSS:CoherentPresentationsPolygraphicResolutions}

Let~$A$ be an algebra. A \emph{coherent presentation of~$A$} is a $2$-polygraph~$X$ such that $\pres{X_0}{X_1}$ is a presentation of~$A$, and such that~$X_2$ is acyclic. A \emph{polygraphic resolution of~$A$} is an $\infty$-poly\-graph~$X$ such that the $1$-polygraph $\pres{X_0}{X_1}$ is a presentation of~$A$, and, for every~$n\geq 1$, the cellular extension~$X_{n+1}$ of~$\lin{X_n}$ is acyclic.
For example, if $X=\pres{X_0}{X_1}$ is a presentation of~$A$, then the $2$-polygraph
$\cohpres{X_0}{X_1}{\Sph(\lin{X_1})}$
is a coherent presentation of~$A$, and, if $X=\cohpres{X_0}{X_1}{X_2}$ is a coherent presentation of~$A$, then the $\infty$-polygraph 
\[
\cohpres{X_0}{X_1}{X_2 \,\vert\, \Sph(\lin{X_2}) \,\vert\, \Sph(\lin{\Sph(\lin{X_2})}) \,\vert\, \cdots}
\]
is a polygraphic resolution of~$A$.

\section{Convergent presentations of associative algebras}
\label{S:LinearRewriting}

This section develops a rewriting theory for associative algebras, based on the polygraphs introduced in Section~\ref{S:nPolygraphs}. For that, the usual rewriting notions are introduced: rewriting steps, normal forms, termination, confluence and convergent presentations. As usual in rewriting theories, Theorem~\ref{T:StandardBasis} asserts that a convergent presentation of an associative algebra yields a basis of that algebra, together with a mechanism to compute its product. In the end of the section, Propositions~\ref{P:ConvergentGröbner} and~\ref{P:ConvergentPBW} detail how convergent presentations of associative algebras generalise other rewriting-like objects, namely Gröbner bases and Poincaré-Birkhott-Witt bases.

\subsection{Rewriting steps and normal forms}
\label{SS:RewritingSteps}

\subsubsection{Rewriting steps and positive \pdf{1}-cells}
\label{SSS:RewritingSteps}
Assume that~$X$ is a left-monomial $1$-polygraph. 
A \emph{rewriting step of~$X$} is a $1$-cell $\lambda f + 1_a$ of size~$1$ of the free $1$-algebra~$\lin{X}$ that satisfies the condition
\[
\supp(\lambda s(f) + a) \;=\; \ens{s(f)} \,\amalg\, \supp(a),
\]
that is, such that $\lambda\neq 0$ and $s(f)\notin\supp(a)$. A $1$-cell of~$\lin{X}$ is called \emph{positive} if it is a (possibly empty) $0$-composite $f_1\star_0\cdots\star_0 f_p$ of rewriting steps of~$\lin{X}$.

\begin{lemma}
\label{L:FactElem1Cell}
Let~$X$ be a left-monomial $1$-polygraph. Every $1$-cell~$f$ of size~$1$ of~$\lin{X}$ can be decomposed into $f=g\star_0 h^-$, where each of~$g$ and~$h$ is either an identity or a rewriting step of~$\lin{X}$.
\end{lemma}

\begin{proof}
Write $f=\lambda f' + 1_b$, where $f':u\fl a$ is a $1$-monomial of~$\lin{X}$. Let~$\mu$ be the coefficient of~$u$ in~$b$, possibly zero, so that $b=\mu u + c$ with ~$c$ such that $\supp(c)$ does not contain~$u$. Put
\[
g = (\lambda + \mu) f' + 1_c
\qquad\text{and}\qquad
h = \lambda 1_a + \mu f' + 1_c.
\]
The linearity of the $0$-composition of~$\lin{X}$ gives $f=g\star_0 h^-$.
Moreover, by hypothesis, $u$ does not belong to any of $\supp(a)$ or $\supp(c)$. As a consequence, each of the $1$-cells~$g$ and~$h$ is either an identity (if $\lambda+\mu=0$ for~$g$, if $\mu=0$ for~$h$) or a rewriting step.
\end{proof}

\subsubsection{Reduced cells and normal forms}

A $0$-cell~$a$ of~$\lin{X}$ is called \emph{reduced} if~$\lin{X}$ contains no rewriting step of source~$a$. As a direct consequence of the definition, we have that a $0$-cell is reduced if, and only if, it is a (possibly empty) linear combination of reduced monomials of~$\lin{X}$. The reduced $0$-cells of~$\lin{X}$ form a linear subspace of the free algebra~$\lin{X_0}$ that is denoted by $\Red(X)$. Because~$X$ is left-monomial, the set of reduced monomials of~$\lin{X}$, denoted by $\red(X)$, forms a basis of $\Red(X)$.

If~$a$ is a $0$-cell of~$\lin{X}$, a \emph{normal form of~$a$} is a reduced $0$-cell~$b$ of~$\lin{X}$ such that there exists a positive $1$-cell of source~$a$ and target~$b$ in $\lin{X}$. 

\subsection{Termination}
\label{SS:Termination}

\subsubsection{Binary relations on free algebras}

Assume that~$X$ is a set and that~$\vdash$ is a binary relation on the free monoid~$X^*$. We say that~$\vdash$ is \emph{stable by product} if $u\vdash u'$ implies $vuw\vdash vu'w$ for all~$u$, $u'$, $v$ and~$w$ in~$X^*$. If~$Y$ is a left-monomial cellular extension of~$\lin{X}$, we say that~$\vdash$ is \emph{compatible with~$Y$} if $u\vdash v$ holds for every $1$-cell $y:u\fl a$ of~$Y$ and every monomial~$v$ in~$\supp(a)$.

The relation~$\vdash$ is extended to the $0$-cells of the free algebra~$\lin{X}$ by putting~$a\vdash b$ if the following two conditions hold:
\begin{enumerate}
\item $\supp(a)\setminus\supp(b) \neq \emptyset$,
\item for every~$v$ in $\supp(b)\setminus\supp(a)$, there exists~$u$ in $\supp(a)\setminus\supp(b)$, such that $u\vdash v$.
\end{enumerate}
As a consequence, if~$u$ is a monomial and~$a$ is a $0$-cell of~$\lin{X}$, then $u\vdash a$ holds if, and only if, $u\vdash v$ holds for every~$v$ in~$\supp(a)$. For this reason, we use the same notation for the relation on~$X^*$ and for its extension to the $0$-cells of~$\lin{X}$.

The relation~$\vdash$ on the $0$-cells of~$\lin{X}$ corresponds to the restriction to finite subsets of~$X^*$ of the so-called multiset relation generated by~$\vdash$. We refer to~\cite[Section 2.5]{BaaderNipkow98} for the general definition and the main properties of multiset relations, and, in particular, the fact that~$\vdash$ is wellfounded on the $0$-cells if, and only if, it is wellfounded on the monomials, the proof being based on König's lemma.

\subsubsection{The rewrite order and termination}
\label{SSS:RewriteOrderTermination}

Let~$X$ be a left-monomial $1$-polygraph. Define~$\succ_X$ as the smallest transitive binary relation on~$X_0^*$ that is stable by product and compatible with~$X_1$. We say that~$X$ \emph{terminates} if the relation~$\succ_X$ is wellfounded. In that case, the reflexive closure~$\succcurlyeq_X$ of the relation~$\succ_X$ is a wellfounded order, called the \emph{rewrite order of~$X$}.

Assume that the $1$-polygraph~$X$ terminates. Then the minimal $0$-cells for the rewrite order of~$X$ are the reduced ones. Moreover, for every nonidentity positive $1$-cell~$a$ of~$\lin{X}$, we have~$s(a)\succ_X t(a)$. This implies that the $1$-algebra~$\lin{X}$ contains no infinite sequence of $0$-composable rewriting steps
\[
a_0 \ofl{f_1} a_1 \ofl{f_1} \cdots \ofl{f_{n-1}} a_{n-1} \ofl{f_n} a_n \ofl{f_{n+1}}\: \cdots
\]
As a consequence, every $0$-cell of~$\lin{X}$ admits at least one normal form.
If~$X$ terminates, induction on the wellfounded order~$\succ_X$ is called \emph{noetherian induction}.

\subsubsection{Monomial orders}

Let~$X$ be a set. A wellfounded total order on~$X^*$, whose strict part is stable by product, is called a \emph{monomial order on~$\lin{X}$}. A classical example of a monomial order is given, for any wellfounded total order relation~$>$ on~$X$, by the \emph{deglex order generated by~$>$}, defined by
\begin{enumerate}
\item $u>_{\text{deglex}}v$ for all monomials~$u$ and~$v$ of~$\lin{X}$ such that~$u$ has greater length than~$v$, and
\item $uxv>_{\text{deglex}}uyw$ for all~$x>y$ of~$X$, and monomials~$u$, $v$ and~$w$ of~$\lin{X}$ such that~$v$ and~$w$ have the same length.
\end{enumerate}

Now, assume that~$X$ is a left-monomial $1$-polygraph. If there exists a monomial order~$\succ$ on~$\lin{X_0}$ that is compatible with~$X_1$, then~$X$ terminates: the order~$\succ$ is wellfounded, and~$a\succ_X b$ implies~$a\succ b$ for all $0$-cells~$a$ and~$b$. However, the converse implication does not hold, as illustrated by the following example.

\subsubsection{Example}
\label{X:xyz=x3+y3+z3}

The following left-monomial $1$-polygraph terminates:
\[
X = \bigpres{x,y,z}{xyz \ofl{\gamma} x^3 + y^3 + z^3}
\]
Indeed, for every monomial~$u$ of~$\lin{X}$, denote by~$A(u)$ the number of factors~$xyz$ that occur in~$u$, by~$B(u)$ the number of~$y$ that~$u$ contains, and put $\Phi(u)=3A(u)+B(u)$. It is sufficient to check that $\Phi(uxyzv)$ is strictly greater than each of $\Phi(ux^3v)$, $\Phi(uy^3v)$ and $\Phi(uz^3v)$, for all monomials~$u$ and~$v$ of~$\lin{X}$:
\[
\begin{array}{r c l c r c l}
\Phi(uxyzv) &= &\Phi(u)+\Phi(v)+4,
&\qquad
&\Phi(ux^3v) &= 
&\begin{cases} 
\Phi(u)+\Phi(v) + 3 
	&\text{if $v=yzv'$} \\ 
\Phi(u)+\Phi(v) 
	&\text{otherwise}, 
\end{cases}
\\ \strut \\
\Phi(uy^3v) &= &\Phi(u)+\Phi(v)+3,
&\qquad
&\Phi(uz^3v) &= 
&\begin{cases} 
\Phi(u)+\Phi(v) + 3 
	&\text{if $u=u'xy$} \\ 
\Phi(u)+\Phi(v) 
	&\text{otherwise}. 
\end{cases}
\end{array}
\]
However, no monomial order on~$\lin{X_0}$ is compatible with~$X_1$, because, for such an order~$\succ$, one of the monomials~$x^3$, $y^3$, $z^3$ is always greater than~$xyz$. Indeed, since~$\succ$ is total, one of~$x$, $y$ or~$z$ is greater than the other two. If it is~$x$, the case of~$z$ being symmetric, $x\succ y$ implies~$x^2\succ yx$ and~$x\succ z$ implies~$yx\succ yz$, so that~$x^2\succ yz$, hence~$x^3\succ xyz$. Now, if~$y\succ x$ and~$y\succ z$, we get~$y^2\succ xy$, thus~$y^2z \succ xyz$ and~$y^3\succ y^2z$. 

\begin{lemma}
\label{L:DecompositionAlgebra}
If~$X$ is a terminating left-monomial $1$-polygraph, then, as a vector space, $\lin{X_0}$ admits the decomposition
\[
\lin{X_0} = \Red(X) + I(X).
\]
\end{lemma}

\begin{proof}
Since~$X$ terminates, every $0$-cell~$a$ of~$\lin{X}$ admits at least a normal form~$b$, i.e.\ a reduced $0$-cell~$b$ such that there exists a positive $1$-cell $f:a\fl b$ in~$\lin{X}$. We conclude by writing $a = b + (a-b)$, and by observing that~$b$ belongs to~$\Red(X)$, by hypothesis, and that~$a-b$ is in~$I(X)$, by Lemma~\ref{L:Ideal}. 
\end{proof}

\subsection{Branchings and confluence}
\label{SS:Confluence}

\subsubsection{Branchings}

Assume that~$X$ is a left-monomial $1$-polygraph.
A \emph{branching of~$X$} is a pair~$(f,g)$ of positive $1$-cells of~$\lin{X}$ with the same source, this $0$-cell being called the \emph{source of~$(f,g)$}. We do not distinguish the branchings~$(f,g)$ and~$(g,f)$. A branching~$(f,g)$ of~$X$ is called \emph{local} if both~$f$ and~$g$ are rewriting steps of~$\lin{X}$. For a branching $(f,g)$ of~$X$ of source~$a$, define the branching 
\[
\lambda u(f,g)v + b = (\lambda ufv + b,\, \lambda ugv + b),
\]
of~$X$ of source $\lambda uav + b$, for all scalar~$\lambda$, monomials~$u$ and~$v$ and $0$-cell~$b$ of~$\lin{X}$. Note that, if~$(f,g)$ is local and~$\lambda\neq 0$, then $\lambda u(f,g) + b$ is also local.

\subsubsection{Classification of local branchings}

Assume that~$X$ is a left-monomial $1$-polygraph.
Given a local branching~$(\lambda u_1 \alpha u_2 + a,\, \mu v_1 \beta v_2 + b)$ of~$X$, we have two main possibilities, depending if $u_1 s(\alpha) u_2=v_1 s(\beta) v_2$ or not. Moreover, in the equality case, there are three different situations, depending on the respective positions of~$s(\alpha)$ and~$s(\beta)$ in this common monomial. This analysis leads to a partition of the local branchings of~$X$ into the following four families.

\begin{enumerate}

\item \emph{Aspherical} branchings: $\lambda (f,f) + b$, for all $1$-monomial $f:u\fl a$ of~$\lin{X}$, nonzero scalar~$\lambda$, and $0$-cell~$b$ of~$\lin{X}$, with $u\notin\supp(b)$.

\item \emph{Additive} branchings: $(\lambda f+\mu v+c,\,\lambda u+\mu g+c)$, for all $1$-monomials $f:u\fl a$ and $g:v\fl b$ of~$\lin{X}$, nonzero scalars~$\lambda$ and~$\mu$, and $0$-cell~$c$ of~$\lin{X}$, with $u\neq v$ and $u,v\notin\supp(c)$.

\item \emph{Peiffer} branchings: $\lambda(fv,ug)+c$, for all $1$-monomials $f:u\fl a$ and $g:v\fl b$ of~$\lin{X}$, nonzero scalar~$\lambda$, and $0$-cell~$c$ of~$\lin{X}$, with $uv\notin\supp(c)$.

\item \emph{Overlapping} branchings: $\lambda(f,g)+c$, for all $1$-monomials $f:u\fl a$ and $g:u\fl b$ of~$\lin{X}$ such that~$(f,g)$ is neither aspherical nor Peiffer, every nonzero scalar~$\lambda$, and every $0$-cell~$c$ of~$\lin{X}$, with $u\notin\supp(c)$.

\end{enumerate}

The \emph{critical branchings of~$X$} are the overlapping branchings of~$X$ such that~$\lambda=1$ and~$c=0$, and that cannot be factored $(f,g)=u(f',g')v$ in a nontrivial way. Note that an overlapping branching has a unique decomposition $\lambda u(f,g)v + c$, with~$(f_0,g_0)$ critical.

\subsubsection{Confluence}

Assume that~$X$ is a left-monomial $1$-polygraph.
A branching~$(f,g)$ of~$X$ is called \emph{confluent} if there exist positive $1$-cells~$h$ and~$k$ of~$\lin{X}$ as in
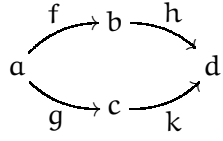
\[
\vcenter{\xymatrix @R=0.25em {
& b
	\ar@/^/ [dr] ^-{h}
\\
a
	\ar@/^/ [ur] ^-{f}
	\ar@/_/ [dr] _-{g}
&& d
\\
& c
	\ar@/_/ [ur] _-{k}
}}
\]
If~$a$ is a $0$-cell of~$\lin{X}$, we say that~$X$ is \emph{confluent at~$a$} (resp.\ \emph{locally confluent at~$a$}, resp.\ \emph{critically confluent}) if every branching (resp.\ local branching, resp.\ critical branching) of~$X$ of source~$a$ is confluent. We say that~$X$ is \emph{confluent} (resp.\ \emph{locally confluent}, resp.\ \emph{critically confluent}) if it is so at every $0$-cell of~$\lin{X}$. Observe that confluence implies that every $0$-cell of~$\lin{X}$ admits at most one normal form. 

\begin{proposition}
\label{P:CharacterisationConfluence}
Let~$X$ be a terminating left-monomial $1$-polygraph. The following assertions are equivalent:
\begin{enumerate}
\item $X$ is confluent.
\item Every $0$-cell of~$I(X)$ admits~$0$ as a normal form.
\item As a vector space, $\lin{X_0}$ admits the direct decomposition $\lin{X_0} = \Red(X) \oplus I(X)$.
\end{enumerate}
\end{proposition}

\begin{proof}
\Item{i}$\:\dfl\:$\Item{ii}. By Lemma~\ref{L:Ideal}, if~$a$ is in~$I(X)$, then there exists a $1$-cell $f:a\fl 0$ in~$\lin{X}$. Since~$X$ is confluent, this implies that~$a$ and~$0$ have the same normal form, if any. And, since~$0$ is reduced, this implies that~$0$ is a normal form of~$a$.

\Item{ii}$\:\dfl\:$\Item{iii}. By Lemma~\ref{L:DecompositionAlgebra}, it is sufficient to prove that $\Red(X)\cap I(X)$ is reduced to~$0$. On the one hand, if~$a$ is in~$\Red(X)$, then~$a$ is reduced and, thus, admits itself as only normal form. On the other hand, if~$a$ is in~$I(X)$, then~$a$ admits~$0$ as a normal form by hypothesis.

\Item{iii}$\:\dfl\:$\Item{i}. Consider a branching~$(f,g)$ of~$X$, with $f:a\fl b$ and $g:a\fl c$. Since~$X$ terminates, each of~$b$ and~$c$ admits at least one normal form, say~$b'$ and~$c'$ respectively. Hence, there exist positive $1$-cells $h:b\fl b'$ and $k:c\fl c'$ in~$\lin{X}$. Note that the difference~$b'-c'$ is also reduced. Moreover, the $1$-cell $(f\star_0 h)^- \star_0 (g\star_0 k)$ has~$b'$ as source and~$c'$ as target. This implies, by Lemma~\ref{L:Ideal}, that~$b'-c'$ also belongs to~$I(X)$. The hypothesis gives~$b'-c'=0$, so that~$(f,g)$ is confluent.
\end{proof}

\subsection{Convergence}
\label{SS:Convergence}

\subsubsection{Convergence}

Let~$X$ be a left-monomial $1$-polygraph.
We say that~$X$ is \emph{convergent} if it is both terminating and confluent. In that case, every $0$-cell~$a$ of~$\lin{X}$ has a unique normal form, denoted by~$\rep{a}$, such that~$\cl{a}=\cl{b}$ holds in~$\cl{X}$ if, and only if, $\rep{a}=\rep{b}$ holds in~$\lin{X}$. 
Hence, if~$X$ is a convergent presentation of an algebra~$A$, the assignment of each element~$a$ of~$A$ to the normal form of any representative of~$a$ in~$\lin{X}$, written~$\rep{a}$ by extension, defines a section~$A\ifl\lin{X}$ of the canonical projection, where~$A$ is seen as a $1$-algebra with identity $1$-cells only. The section is linear (the normal form of $\lambda a + \mu b$ is $\lambda\rep{a}+\mu\rep{b}$), it preserves the unit (termination implies $\rep{1}=1$), but $\rep{ab}\neq \rep{a}\rep{b}$ in general.

\begin{theorem}
\label{T:StandardBasis}
Let~$A$ be an algebra and~$X$ be a convergent presentation of~$A$. Then the set~$\red(X)$ is a linear basis of~$A$. As a consequence, the vector space~$\Red(X)$, equipped with the product defined by $a\cdot b=\rep{ab}$, is an algebra that is isomorphic to~$A$.
\end{theorem}

\begin{proof}
If~$X$ is convergent, Proposition~\ref{P:CharacterisationConfluence} induces that the following sequence of vector spaces is exact:
\[
\xymatrix @M=7.5pt @C=1.5em {
0
	\ar[r] 
& I(X)
	\ar@{>->}[r] 
& \lin{X_0}
	\ar@{->>}[r]
& \Red(X)
	\ar[r] 
& 0
}
\]
Thus, since the algebra~$\lin{X_0}/I(X)$ is isomorphic to~$\cl{X}$, convergence implies that~$\red(X)$ is a linear basis of~$\cl{X}$. We deduce that~$\Red(X)$ and~$\cl{X}$ are isomorphic as vector spaces. There remains to transport the product of~$\cl{X}$ to~$\Red(X)$ to get the result.
\end{proof}

\begin{example}
\label{X:BasisReduced}

Let~$A$ be the algebra presented by the $1$-polygraph $X=\pres{x,y}{\alpha: xy \fl x^2}$, which terminates, because $xy>x^2$ holds for the deglex order generated by~$y>x$. This presentation is also confluent, because it has no critical branching (see Corollary~\ref{C:Critical}). Hence, the set 
\[
\red(X) = \enspres{y^i x^j} {i,j\in\Nb}
\]
is a linear basis of the algebra~$A$. Moreover, the product defined by
\[
y^i x^j \cdot y^k x^l
\:=\:
\begin{cases}
y^i x^{j + k + l}
	&\text{if $j\geq k$}
\\
y^{i - j + k} x^{2j + l}
	&\text{if $j\leq k$}
\end{cases}
\]
turns~$\Red(X)$ into an algebra that is isomorphic to~$A$.

Now, consider the presentation~$Y=\pres{x,y}{\beta : x^2 \fl xy}$ of~$A$. Termination of~$Y$ follows from the deglex order generated by~$x>y$, but~$Y$ is not confluent, since it has a nonconfluent critical branching:
\[
\xymatrix @R=0em {
& xyx
\\
x^3
  \ar@/^/ [ur] ^-{\beta x}
  \ar@/_/ [dr] _-{x\beta}
\\
& x^2y 
  \ar [r] _-{\beta y}
& xy^2
}
\]
Thus the $0$-cell $xyx-xy^2$ is both in~$\Red(Y)$ and~$I(Y)$, proving that the sum $\Red(Y)+I(Y)$ is not direct. As a consequence, $\red(Y)$ is not a linear basis of~$A$. 
\end{example}

\subsubsection{Reduced convergent presentations}

Let~$X$ be a left-monomial $1$-polygraph.
We say that~$X$ is \emph{left-reduced} if, for every $1$-cell~$\alpha$ of~$X$, the only rewriting step of~$\lin{X}$ of source~$s(\alpha)$ is~$\alpha$ itself. We say that~$X$ is \emph{right-reduced} if, for every $1$-cell~$\alpha$ of~$X$, the $0$-cell~$t(\alpha)$ is reduced. We say that~$X$ is \emph{reduced} if it is both left-reduced and right-reduced.

Assume that~$X$ is convergent. Then one obtains a reduced convergent left-monomial $1$-polygraph~$Y$, such that~$\cl{X}\simeq\cl{Y}$, through the following operations. First, replace every $1$-cell $a\fl b$ with $a\fl\rep{a}$. Then, if there exist parallel $1$-cells, discard all of them but one. Finally, eliminate all the remaining $1$-cells of $X$ whose source contains the source of another $1$-cell. 

\subsubsection{Completion of presentations}

The completion procedure, developed by Buchberger for commutative algebras~\cite{Buchberger65} and by Knuth and Bendix for term rewriting systems~\cite{KnuthBendix70}, adapts to terminating left-monomial $1$-polygraphs as follows, to transform them into convergent ones.

Fix a left-monomial $1$-polygraph~$X$, and a well-founded strict order that is stable by product and compatible with~$X_1$. For each nonconfluent critical branching~$(f,g)$ of~$X$, consider $a=c-d$, where~$c$ and~$d$ are arbitrary normal forms of~$t(f)$ and~$t(g)$, respectively. If $\supp(a)$ contains a maximal element~$u$, add the $1$-cell $u\fl b$ to~$X$, where~$b$ is defined by $a=\lambda u + b$ and $u\notin\supp(b)$; otherwise, the procedure fails. After the exploration of all the critical branchings of~$X$, the procedure, if it has not failed, yields a terminating left-monomial $1$-polygraph~$Y$ such that~$\cl{X}\simeq\cl{Y}$. If~$Y$ is not confluent, restart with~$Y$. The procedure either stops when it reaches a convergent left-monomial $1$-polygraph, or runs forever. 

\subsection{Comparison with Gröbner bases and Poincaré-Birkhoff-Witt bases}
\label{SS:PBWGröbner}

\subsubsection{Gröbner bases}

Let~$X$ be a set and~$\preccurlyeq$ be a monomial order on the free algebra~$\lin{X}$. If~$a$ is a nonzero element of~$\lin{X}$, the \emph{leading monomial of~$a$} is the maximum element~$\lm(a)$ of~$\supp(a)$ for~$\preccurlyeq$ (or~$0$ if~$\supp(a)$ is empty), the \emph{leading coefficient of~$a$} is the coefficient~$\lc(a)$ of~$\lm(a)$ in~$a$, and the \emph{leading term of~$a$} is the element $\lt(a)=\lc(a)\lm(a)$ of~$\lin{X}$. 
Observe that, for~$a$ and~$b$ in~$\lin{X}$, we have~$a \prec b$ if, and only if, either $\lm(a) \prec \lm(b)$ or ($\lt(a)=\lt(b)$ and $a-\lt(a) \prec b-\lt(b)$). 

Let~$I$ be an ideal of~$\lin{X}$. A \emph{Gröbner basis for~$(I,\preccurlyeq)$} is a subset~$\Gr$ of~$I$ such that the ideals of~$\lin{X}$ generated by~$\lm(I)$ and by~$\lm(\Gr)$ coincide.

\begin{proposition}
\label{P:ConvergentGröbner}
If~$X$ is a convergent left-monomial $1$-polygraph, and~$\preccurlyeq$ is a monomial order on~$\lin{X_0}$ that is compatible with~$X_1$, then the set~$\dr(X_1)$ of boundaries of $1$-cells of~$X$ is a Gröbner basis for~$(I(X),\preccurlyeq)$. 

Conversely, let~$X$ be a set, let~$\preccurlyeq$ be a monomial order on~$\lin{X}$, let~$I$ be an ideal of~$\lin{X}$ and~$\Gr$ be a subset of~$I$. Define $\Lead(\Gr)$ as the $1$-polygraph with $0$-cells~$X$ and one $1$-cell
\[
\lm(a) \:\ofl{\raisebox{0.5ex}{$\alpha_a$}}\: \lm(a) - \frac{1}{\lc(a)} a
\]
for each~$a$ in~$\Gr$. If~$\Gr$ is a Gröbner basis for~$(I,\preccurlyeq)$, then~$\Lead(\Gr)$ is a convergent left-monomial presentation of~$\lin{X}/I$, such that $I(\Lead(\Gr))=I$, and~$\preccurlyeq$ is compatible with~$\Lead(\Gr)_1$.
\end{proposition}

\begin{proof}
If~$X$ is convergent, then~$\dr(\alpha)$ is in~$I(X)$ for every $1$-cell~$\alpha$ of~$X$. Since~$\preccurlyeq$ is compatible with~$X_1$, we have $\lm(\dr(\alpha))=s(\alpha)$ for every $1$-cell~$\alpha$ of~$X$. Now, if~$a$ is in~$I(X)$, it is a linear combination
\[
a = \sum_{i} \lambda_i u_i \dr(\alpha_i) v_i
\]
of $1$-cells $u_i\dr(\alpha_i) v_i$, where~$\alpha_i$ is a $1$-cell of~$X$, and~$u_i$ and~$v_i$ are monomials of~$\lin{X}$. This implies that
\[
\lm(a) = u_is(\alpha_i)v_i = u_i \lm(\dr(\alpha_i)) v_i
\]
hold for some~$i$. Thus~$\dr(X_1)$ is a Gröbner basis for~$(I(X),\preccurlyeq)$.

Conversely, assume that~$\Gr$ is a Gröbner basis for~$(I,\preccurlyeq)$. By definition, $\preccurlyeq$ is compatible with~$\Lead(\Gr)_1$, hence~$\Lead(\Gr)$ terminates, and $I(\Lead(\Gr))=I$ holds, so that the algebra presented by~$\Lead(\Gr)$ is indeed isomorphic to $\lin{X}/I$. Moreover, the reduced monomials of~$\lin{\Lead(\Gr)}$ are the monomials of~$\lin{X}$ that cannot be decomposed as~$u\lm(a)v$ with~$a$ in~$\Gr$, and~$u$ and~$v$ monomials of~$\lin{X}$. Thus, if a reduced $0$-cell~$a$ of~$\lin{\Lead(\Gr)}$ is in~$I$, its leading monomial must be~$0$, because~$\Gr$ is a Gröbner basis of~$(I,\preccurlyeq)$. As a consequence of Proposition~\ref{P:CharacterisationConfluence}, we get that~$\Lead(\Gr)$ is confluent. 
\end{proof}

\subsubsection{Poincaré-Birkhoff-Witt bases}

Let~$A$ be an $N$-homogeneous algebra, for~$N\geq 2$, let~$X$ be a generating set of~$A$, concentrated in degree~$1$, and let~$\preccurlyeq$ be a monomial order on~$\lin{X}$. A \emph{Poincaré-Birkhoff-Witt} (\emph{PBW}) \emph{basis for~$(A,X,\preccurlyeq)$} is a subset~$\Br$ of~$X^*$  such that:
\begin{enumerate} 
\item $\Br$ is a linear basis of~$A$, with~$[u]_{\Br}$ denoting the decomposition of an element~$u$ of~$X^*$ in the basis~$\Br$,
\item for all~$u$ and~$v$ in~$\Br$, we have $uv\succcurlyeq [uv]_{\Br}$, 
\item an element~$u$ of~$X^*$ belongs to~$\Br$ if, and only if, for every decomposition $u=vu'w$ of~$u$ in~$X^*$ such that~$u'$ has degree~$N$, then~$u'$ is in~$\Br$.
\end{enumerate}

\begin{proposition}
\label{P:ConvergentPBW}
If~$X$ is a convergent left-monomial $N$-homogeneous presentation of an algebra~$A$, and~$\preccurlyeq$ is a monomial order on~$\lin{X_0}$ that is compatible with~$X_1$, then the set~$\red(X)$ of reduced monomials of~$\lin{X}$ is a PBW basis for~$(A,X_0,\preccurlyeq)$.

Conversely, let~$A$ be an $N$-homogeneous algebra, let~$X$ be a generating set of~$A$ that is concentrated in degree~$1$, let~$\preccurlyeq$ a monomial order on~$\lin{X}$, and~$\Br$ be a PBW basis of~$(A,X,\preccurlyeq)$. Define~$\tilde{\Br}$ as the $1$-polygraph with $0$-cells~$X$ and with one $1$-cell
\[
uv \:\ofl{\raisebox{0.5ex}{$\alpha_{u,v}$}}\: [uv]_{\Br}
\]
for all~$u$ and~$v$ in~$\Br$ such that~$uv$ has degree~$N$ and~$uv\neq [uv]_{\Br}$. Then~$\tilde{\Br}$ is a convergent left-monomial $N$-homogeneous presentation of~$A$, such that $\red(\tilde{\Br})=\Br$, and~$\preccurlyeq$ is compatible with~$\tilde{\Br}_1$.
\end{proposition}

\begin{proof}
If~$X$ is a convergent left-monomial presentation of~$A$, Theorem~\ref{T:StandardBasis} implies that the set~$\red(X)$ of reduced monomials of~$\lin{X}$ is a linear basis of~$A$. The fact that~$\preccurlyeq$ is compatible with~$X_1$ implies Axiom~\Item{ii} of a PBW basis, and Axiom~\Item{iii} comes from the definition of a reduced monomial for an $N$-homogeneous left-monomial $1$-polygraph.

Conversely, assume that~$\Br$ is a PBW basis for~$(A,X,\preccurlyeq)$. By definition, $\tilde{\Br}$ is $N$-homogeneous and left-monomial, and Axiom~\Item{iii} of a PBW basis implies $\red(\tilde{\Br})=\Br$. Termination of~$\tilde{\Br}$ is given by Axiom~\Item{ii} of a PBW basis, because~$\preccurlyeq$ is wellfounded. By Proposition~\ref{P:CharacterisationConfluence}, it is sufficient to prove that $\Red(\tilde{\Br})\cap I(\tilde{\Br})=0$ to get confluence: on the one hand, a reduced $0$-cell~$a$ of~$\Red(\tilde{\Br})$ is a linear combination of $0$-cells of~$\Br$, so that~$a$ is its only normal form; and, on the other hand, if~$a$ belongs to~$I(\tilde{\Br})$, then~$a$ admits~$0$ as a normal form by Lemma~\ref{L:Ideal}. Finally, the algebra presented by~$\tilde{\Br}$ is isomorphic to~$\Red(\tilde{\Br})$, that is to~$\K\Br$, hence to~$A$, by Theorem~\ref{T:StandardBasis} and because~$\Br$ is a linear basis of~$A$.
\end{proof}

\section{Coherent presentations of associative algebras}
\label{S:Squier}

In this section, we define a coherent presentation of an associative algebra as a presentation by generators and relations, extended with a family of $2$-cells that generates all the ``relations among relations''. We prove Proposition~\ref{P:HNewman} and Theorem~\ref{T:HCritical}, which are coherent versions of two classical results of rewriting theory: Newman's lemma and the critical branchings theorem. From these first results, we derive Squier's theorem for associative algebras (Theorem~\ref{T:Squier}), stating that the critical branchings of a left-monomial convergent presentation generate a coherent presentation. The section ends with two examples of applications of Squier's theorem.

\subsection{Coherent confluence and the coherent Newman's lemma}
\label{SS:Newman}

\subsubsection{Coherent confluence and convergence}

Let~$X$ be a left-monomial $1$-polygraph, and let~$Y$ be a cellular extension of the free $1$-algebra~$\lin{X}$. A branching~$(f,g)$ of~$X$ is \emph{$Y$-confluent} if there exist positive $1$-cells~$h$ and~$k$ in~$\lin{X}$ and a $2$-cell~$F$ in~$\lin{X}[Y]$ as in
\[
\vcenter{\xymatrix @R=0.25em{
& b
	\ar@/^/ [dr] ^-{h}
	\ar@2 []!<0pt,-10pt>;[dd]!<0pt,10pt> ^-*+{F}
\\
a
	\ar@/^/ [ur] ^-{f}
	\ar@/_/ [dr] _-{g}
&& d
\\
& c
	\ar@/_/ [ur] _-{k}
}}
\]
If~$a$ is a $0$-cell of~$\lin{X}$, say that~$X$ is \emph{$Y$-confluent} (resp.\ \emph{locally $Y$-confluent}, resp.\ \emph{critically $Y$-confluent}) \emph{at~$a$} if every branching (resp.\ local branching, resp.\ critical branching) of~$X$ of source~$a$ is $Y$-confluent. Say that~$X$ is \emph{$Y$-confluent} (resp.\ \emph{locally $Y$-confluent}, resp.\ \emph{critically $Y$-confluent}) if it is so at every $0$-cell of~$\lin{X}$, and that~$X$ is \emph{$Y$-convergent} if it is terminating and $Y$-confluent. 
Note that, if $Y=\Sph(\lin{X})$, the $Y$-confluence and $Y$-convergence properties boil down to the confluence and convergence of Section~\ref{S:LinearRewriting}.

\begin{lemma}
\label{L:Confluence2Cell}
Let~$X$ be a left-monomial $1$-polygraph, and $Y$ be a cellular extension of~$\lin{X}$, such that~$X$ is $Y$-confluent at every $0$-cell~$b\prec a$ for some fixed $0$-cell~$a$ of~$\lin{X}$. Let~$f$ be a $1$-cell of~$\lin{X}$ that admits a decomposition
\[
a_0 \:\ofl{f_1}\: a_1 \:\ofl{f_2}\: \cdots \:\ofl{f_p}\: a_p
\]
into $1$-cells of size~$1$. If~$a_i\prec a$ holds for every~$0<i<p$, then there exist positive $1$-cells~$g$ and~$h$ in~$\lin{X}$ and a $2$-cell~$F$ in~$\lin{X}[Y]$ as in
\[
\xymatrix @R=0.5em {
& a_p
	\ar@/^/ [dr] ^-{h}
	\ar@2 []!<0pt,-12pt>;[d]!<0pt,3pt> ^-*+{F}
\\	
a_0 
	\ar@/^/ [ur] ^-{f} 
	\ar@/_/ [rr] _-{g}
&& a'
}
\]
\end{lemma}

\begin{proof}
Proceed by induction on~$p$. If~$p=0$, then~$f$ is an identity, so taking~$g=h=1_{a_0}$ and~$F=1_f$ proves the result. Otherwise, construct 
\[
\xymatrix @R=1em {
&& a_p
	\ar@/^/ [dr] ^-{h_2}
	\ar@2 []!<0pt,-10pt>;[d]!<0pt,10pt> ^-*+{F}
\\
& a_1 
	\ar@/^/ [ur] ^(0.4){f_2\star_0\cdots\star_0 f_p}
	\ar [rr] |-*+{g_2}
	\ar [dr] |-{h_1}
	\ar@{} []!<-5pt,0pt>;[d]!<-5pt,0pt> |(0.66){\sm =}
&& b_2
	\ar@/^/ [dr] ^-{k_2}
	\ar@2 []!<-15pt,-12.5pt>;[d]!<-15pt,7.5pt> ^-*+{G}
\\	
a_0
	\ar@/^/ [ur] ^-{f_1}
	\ar@/_/ [rr] _-{g_1}
&& b_1 
	\ar@/_/ [rr] _-{k_1}
&& a'
}
\]
Apply Lemma~\ref{L:FactElem1Cell} to the $1$-cell~$f_1$ of size~$1$ to get the positive $1$-cells~$g_1$ and~$h_1$ such that $f_1=g_1\star_0 h_1^-$. We have~$a_i\prec a$ for every~$1<i<p$, so the induction hypothesis applies to $f_2\star_0\cdots\star_0 f_p$, providing the positive $1$-cells~$g_2$ and~$h_2$, and the $2$-cell~$F$. Then, consider the branching~$(h_1, g_2)$, whose source~$a_1$ satisfies~$a_1\prec a$: by hypothesis, the branching~$(h_1, g_2)$ is $Y$-confluent, giving the positive $1$-cells~$k_1$ and~$k_2$, and the $2$-cell~$G$.
\end{proof}

\begin{proposition}[Coherent Newman's lemma]
\label{P:HNewman}
Let~$X$ be a terminating left-monomial $1$-polygraph, and~$Y$ be a cellular extension of~$\lin{X}$. If~$X$ is locally $Y$-confluent then it is $Y$-confluent.
\end{proposition}

\begin{proof}
Prove that~$X$ is $Y$-confluent at every $0$-cell~$a$ of~$\lin{X}$ by noetherian induction on~$a$. 
If~$a$ is reduced, then~$(1_a,1_a)$ is the only branching of source~$a$, and it is $Y$-confluent, taking~$f'=g'=1_a$ and~$F=1_{1_a}$. 
Now, let~$a$ be a nonreduced $0$-cell of~$\lin{X}$ such that~$X$ is $Y$-confluent at every $0$-cell~$b\prec a$, and let~$(f,g)$ be a branching of~$X$ of source~$a$. If one of~$f$ or~$g$ is an identity, say~$f$, then~$(f,g)$ is $Y$-confluent, taking $f'=g$ and $g'=1_{t(g)}$, with~$F=1_g$. Otherwise, prove that~$(f,g)$ is $Y$-confluent thanks to the construction
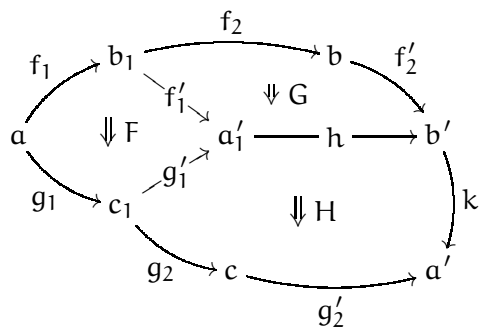
\[
\xymatrix @R=1em {
& b_1
	\ar@/^/ [rr] ^-{f_2} 
	\ar [dr] |-{f'_1}
	\ar@2 []!<-5pt,-22.5pt>;[dd]!<-5pt,22.5pt> ^-*+{F}
& \strut
	\ar@2 []!<15pt,-2.5pt>;[d]!<15pt,12.5pt> ^-*+{G}
& b
	\ar@/^/ [dr] ^-{f'_2}
\\
a 
	\ar@/^/ [ur] ^-{f_1}
	\ar@/_/ [dr] _-{g_1}
&& a'_1
	\ar [rr] |-*+{h}
	\ar@2 []!<25pt,-22.5pt>;[dd]!<25pt,17.5pt> ^-*+{H}
&& b'
	\ar@/^/ [dd] ^-{k}
\\
& c_1
	\ar [ur] |-{g'_1}
	\ar@/_/ [dr] _-{g_2}
\\
&& c
	\ar@/_/ [rr] _-{g'_2}
&& a'
}
\]
Since~$f$ and~$g$ are not identities, they admit decompositions $f=f_1\star_0 f_2$ and $g=g_1\star_0 g_2$ where~$f_1$ and~$g_1$ are rewriting steps, and~$f_2$ and~$g_2$ are positive $1$-cells. By hypothesis, the local branching $(f_1,g_1)$ is $Y$-confluent, yielding the positive $1$-cells~$f'_1$ and~$g'_1$ and the $2$-cell~$F$. Since both~$a\succ b_1$ and~$a\succ c_1$ hold, the induction hypothesis applies to the branching~$(f_2, f'_1)$ to give~$f'_2$, $h$ and~$G$, and, then, to the branching~$(g'_1\star_0 h, g_2)$ to give~$h$, $g'_2$ and~$H$. 
\end{proof}

Taking $Y=\Sph(\lin{X})$ in Proposition~\ref{P:HNewman}, we deduce Newman's lemma~\cite{Newman42}:

\begin{corollary}
\label{C:Newman}
For terminating left-monomial $1$-polygraphs, confluence and local confluence are equivalent properties.
\end{corollary}

\subsection{The coherent critical branchings theorem}
\label{SS:CriticalBranchings}

\begin{theorem}[The coherent critical branchings theorem]
\label{T:HCritical}
Assume that~$X$ is a terminating left-monomial $1$-polygraph, and that~$Y$ is a cellular extension of~$\lin{X}$. If~$X$ is critically $Y$-confluent, then~$X$ is locally $Y$-confluent.
\end{theorem}

\begin{proof}
Proceed by noetherian induction on the sources of the local branchings to prove that~$X$ is locally $Y$-confluent at every $0$-cell of~$\lin{X}$. 
First, note that a reduced $0$-cell cannot be the source of a local branching, so~$X$ is locally $Y$-confluent at reduced $0$-cells.
Next, fix a nonreduced $0$-cell~$a$ of~$\lin{X}$, and assume that~$X$ is locally $Y$-confluent at every~$b\prec a$. With a termination-based argument similar to that of the coherent Newman's lemma, we deduce that~$X$ is $Y$-confluent at every~$b\prec a$. Then proceed by case analysis on the type of the local branchings.

\smallskip
\Item{i} An aspherical branching $\lambda(f,f)+b$ is always $Y$-confluent.

\smallskip
\Item{ii} For an additive branching, construct
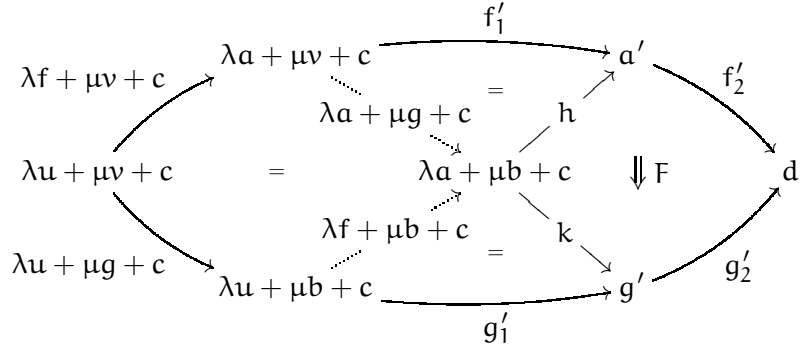
\[
\xymatrix @R=2.5em @C=1em {
& {\lambda a + \mu v + c}
	\ar@/^/ [rr] ^-{f'_1}
	\ar@{.>} [dr] |-*+{\lambda a + \mu g + c}
&& a'
	\ar@/^/ [drrr] ^-{f'_2}
	\ar@2 []!<2.5pt,-37.5pt>;[dd]!<2.5pt,37.5pt> ^-*+{F}
\\
{\lambda u + \mu v + c}
	\ar@/^2ex/ [ur] ^-*+{\lambda f + \mu v + c}
	\ar@/_2ex/ [dr] _-*+{\lambda u + \mu g + c}
	\ar@{} [rr] |(0.45){\sm =}
&& {\lambda a + \mu b + c}
	\ar [ur] |-*+{h}
	\ar [dr] |-*+{k}
	\ar@{} [u] |(0.7){\sm =}
	\ar@{} [d] |(0.7){\sm =}
&&&& d
\\
& {\lambda u + \mu b +c}
	\ar@{.>} [ur] |-*+{\lambda f + \mu b + c}
	\ar@/_/ [rr] _-{g'_1}
&& {g'}
	\ar@/_/ [urrr] _-{g'_2}
}
\]
By linearity of the $0$-composition, we have
\[
(\lambda f + \mu v + c) \star_0 (\lambda a + \mu g + c) 
	= \lambda f + \mu g + c
	= (\lambda u + \mu g + c) \star_0 (\lambda f + \mu b + c).
\]
Note that the dotted $1$-cells $\lambda a + \mu g + c$ and $\lambda f + \mu b + c$ are not positive in general, since~$u$ can be in~$\supp(b)$ or~$v$ in~$\supp(a)$. However, those $1$-cells are of size~$1$, and Lemma~\ref{L:FactElem1Cell} applies to both of them, to give positive $1$-cells~$f'_1$, $g'_1$, $h$ and~$k$ that satisfy
\[
f'_1 = (\lambda a + \mu g + c) \star_0 h
\qquad\text{and}\qquad
g'_1 = (\lambda f + \mu b + c) \star_0 k.
\]
Now, $u\succ a$, $v\succ b$, $\lambda\neq 0$ and~$\mu\neq 0$ imply $\lambda u + \mu v + c \succ \lambda a + \mu b + c$. Thus, the branching $(h,k)$ is $Y$-confluent by hypothesis, yielding the positive $1$-cells~$f'_2$ and~$g'_2$ and the $2$-cell~$F$.

\medskip
\textbf{(iii)} 
In the case of a Peiffer branching, construct
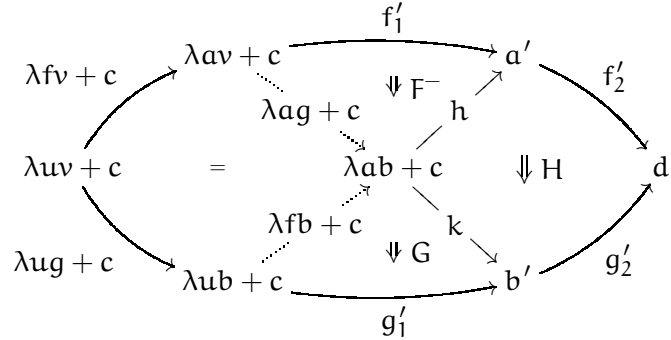
\[
\xymatrix @R=2.5em @C=1.5em {
& \lambda av + c
	\ar@/^/ [rr] ^-{f'_1}
	\ar@{.>} [dr] |-*+{\lambda ag + c}
&& a'
	\ar@/^/ [drr] ^-{f'_2}
	\ar@2 []!<2.5pt,-37.5pt>;[dd]!<2.5pt,37.5pt> ^-*+{H}
\\
\lambda uv + c
	\ar@/^2ex/ [ur] ^-*+{\lambda fv + c}
	\ar@/_2ex/ [dr] _-*+{\lambda ug + c}
	\ar@{} [rr] |(0.45){\sm =}
&& \lambda ab + c
	\ar [ur] |-*+{h}
	\ar [dr] |-*+{k}
	\ar@2 [u]!<0pt,-7.5pt>;[]!<0pt,27.5pt> ^-*+{F^-}
	\ar@2 []!<0pt,-27.5pt>;[d]!<0pt,7.5pt> ^-*+{G}
&&& d
\\
& \lambda ub + c
	\ar@{.>} [ur] |-*+{\lambda fb + c}
	\ar@/_/ [rr] _-{g'_1}
&&	b'
	\ar@/_/ [urr] _-{g'_2}
}
\]
Use the linearity of the $0$-composition to obtain
\[
(\lambda fv + c) \star_0 (\lambda ag + c) = \lambda fg + c = (\lambda ug + c) \star_0 (\lambda fb + c).
\]
Again, the dotted $1$-cells $\lambda fb + c$ and $\lambda ag + c$ are not positive in general: this is the case, for example, if either $\supp(ub)\cap\supp(c)$ or $\supp(av)\cap\supp(c)$ is not empty. Let 
$a = \sum_{i=1}^p \mu_i u_i$
be the canonical decomposition of~$a$. By linearity of the $0$-composition, the $1$-cell $\lambda ag + c$ admits the following decomposition in $1$-cells of size~$1$:
\[
\lambda ag + c = g_1\star_0\cdots\star_0 g_p,
\qquad \text{with} \quad 
g_j \:=\: 
	\sum_{1\leq i < j} \lambda\mu_i u_i b 
	\:+\: \lambda\mu_j u_j g 
	\:+\: \sum_{j<i\leq p} \lambda\mu_i u_i v + c.
\]
We have~$u\succ u_i$ for every~$i$, and~$v\succ b$, giving $\lambda uv + v \succ s(g_j)$ for every~$j$. Hence~$\lambda ag + c$ is eligible to Lemma~\ref{L:Confluence2Cell}, yielding~$f'_1$, $h$ and~$F$. The cells~$g'_1$, $k$ and~$G$ are obtained similarly from $\lambda fb + c$. Finally, $\lambda uv + c \succ \lambda ab + c$ implies, by induction hypothesis, that~$(h,k)$ is $Y$-confluent, giving~$f'_2$, $g'_2$ and~$H$. 

\medskip
\textbf{(iv)} 
For an overlapping branching $(\lambda f + c, \lambda g + c)$, construct
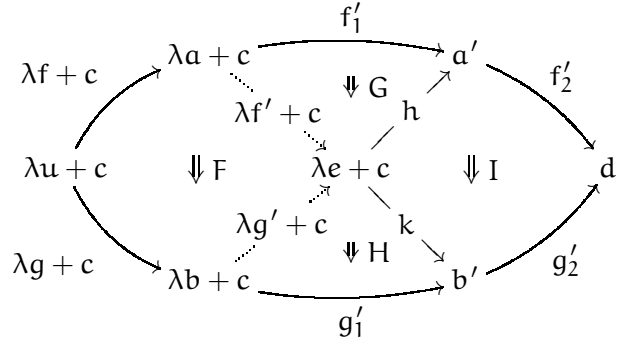
\[
\xymatrix @R=2.5em @C=1.5em {
& \lambda a + c
	\ar@/^/ [rr] ^-{f'_1}
	\ar@{.>} [dr] |-*+{\lambda f' + c}
	\ar@2 []!<-5pt,-37.5pt>;[dd]!<-5pt,37.5pt> ^-*+{F}
&& a'
	\ar@/^/ [drr] ^-{f'_2}
	\ar@2 []!<2.5pt,-37.5pt>;[dd]!<2.5pt,37.5pt> ^-*+{I}
\\
\lambda u + c
	\ar@/^2ex/ [ur] ^-*+{\lambda f + c}
	\ar@/_2ex/ [dr] _-*+{\lambda g + c}
&& \lambda e + c
	\ar [ur] |-*+{h}
	\ar [dr] |-*+{k}
	\ar@2 [u]!<0pt,-7.5pt>;[]!<0pt,27.5pt> ^-*+{G}
	\ar@2 []!<0pt,-27.5pt>;[d]!<0pt,7.5pt> ^-*+{H}
&&& d
\\
& \lambda b + c
	\ar@{.>} [ur] |-*+{\lambda g' + c}
	\ar@/_/ [rr] _-{g'_1}
&& b'
	\ar@/_/ [urr] _-{g'_2}
}
\]
Consider the unique decomposition $(f,g)=v(f_0,g_0)w$, with~$(f_0,g_0)$ critical. Since~$(f_0,g_0)$ is $Y$-confluent by hypothesis, one obtains
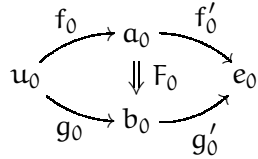
\[
\xymatrix @R=0.25em {
& a_0
	\ar@/^/ [dr] ^-{f_0'}
	\ar@2 []!<0pt,-10pt>;[dd]!<0pt,10pt> ^-*+{F_0}
\\
u_0
	\ar@/^/ [ur] ^-{f_0}
	\ar@/_/ [dr] _-{g_0}
&& e_0
\\
& b_0
	\ar@/_1.4ex/ [ur] _-{g_0'}
}
\]
Define the positive $1$-cells $f'=v f'_0 w$ and $g'=v g'_0 w$, and the $2$-cell~$F=v F_0 w$. As previously, the dotted $1$-cells are not positive in general, if~$\supp(c)$ intersects~$\supp(a)$ or~$\supp(b)$ for example. However, the $1$-cell~$f'$ is positive, so that it is a $0$-composite $f' = l_1\star_0\cdots\star_0 l_p$ of rewriting steps. As a consequence, we have the chain of inequalities
\[
u \:\succ\: a \:=\: s(l_1) \:\succ\: (\cdots) \:\succ\: s(l_p) \:\succ\: e.
\]
Since we have~$\lambda \neq 0$ and $u\notin\supp(c)$ by hypothesis, the inequality $\lambda u + c \succ \lambda s(l_i) + c$ holds for every~$i$, so that the following decomposition of the $1$-cell $\lambda f'+c$ satisfies the hypotheses of Lemma~\ref{L:Confluence2Cell}:
\[
\lambda f' + c \:=\: \big(\lambda l_1 + c \big) \star_1 \cdots \star_1 \big(\lambda l_p + c\big).
\]
This gives~$f'_1$, $h$ and~$G$. Proceed similarly with the $1$-cell $\lambda g' + c$ to obtain~$g'_1$, $k$ and~$H$. Finally, apply the induction hypothesis on~$(h,k)$, since $\lambda u + c \succ \lambda e + c$, to get~$f'_2$, $g'_2$ and~$I$.
\end{proof}

Taking $Y=\Sph(\lin{X})$ in Theorem~\ref{T:HCritical}, we deduce the critical branching theorem~\cite{Nivat73, KnuthBendix70, Huet80}:

\begin{corollary}
\label{C:Critical}
For terminating left-monomial $1$-polygraphs, local confluence and critical confluence are equivalent properties. In particular, a terminating left-monomial $1$-polygraph with no critical branching is convergent.
\end{corollary}

Proposition~\ref{P:ConvergentGröbner} and Corollaries~\ref{C:Newman} and~\ref{C:Critical} imply

\begin{corollary}[Buchberger's criterion]
\label{C:BuchbergerCriterion}
Let~$X$ be a set, $\preccurlyeq$ be a monomial order on the free algebra~$\lin{X}$, and~$I$ be an ideal of~$\lin{X}$. A subset~$\Gr$ of~$I$ is a Gröbner basis for~$(I,\preccurlyeq)$ if, and only if, the $1$-polygraph $\Lead(\Gr)$ of Proposition~\ref{P:ConvergentGröbner} is critically confluent.
\end{corollary}

\begin{remark} 
The critical branching theorem for $1$-polygraphs differs from its set-theoretic counterpart~\cite[3.1.5]{GuiraudMalbos16}. Indeed, in the set-theoretic case, the termination hypothesis is not required, and nonoverlapping branchings are always confluent, independently of critical confluence. The following two counterexamples show that the linear case is different.

On the one hand, some local branchings can be nonconfluent without termination, even if critical confluence holds. Indeed, the $1$-polygraph
\[
\bigpres {x, y, z, t}{xy\ofl{\alpha} xz, \, zt\ofl{\beta} 2yt}
\]
has no critical branching, but it has a nonconfluent additive branching:
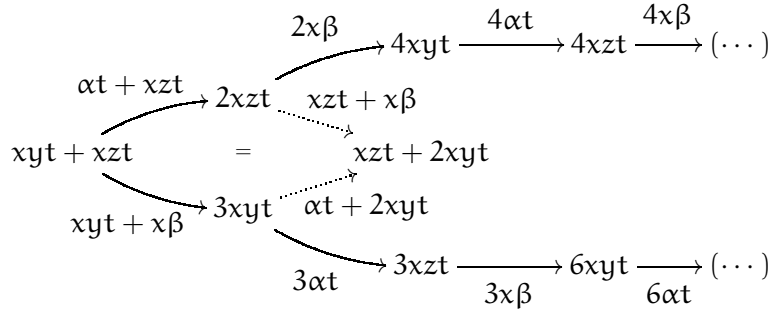
\[
\xymatrix @R=0.5em {
&& 4xyt 
	\ar [r] ^-*+{4\alpha t}
& 4xzt
	\ar [r] ^-*+{4x\beta}
& (\cdots)
\\
& 2xzt
	\ar@/^/ [ur] ^-*+{2x\beta}
	\ar@{.>} [dr] ^(0.55)*+{xzt + x\beta}
\\
xyt + xzt 
	\ar@/^/ [ur] ^-*+{\alpha t + xzt}
	\ar@/_/ [dr] _-*+{xyt + x\beta}
	\ar@{} [rr] |-{\sm =}
&& xzt + 2xyt
\\
& 3xyt
	\ar@{.>} [ur] _(0.55)*+{\alpha t + 2xyt}
	\ar@/_/ [dr] _-*+{3\alpha t}
\\
&& 3xzt
	\ar [r] _-*+{3x\beta}
& 6xyt
	\ar [r] _-*+{6\alpha t}
& (\cdots)
}
\]
The only positive $1$-cells of source~$2xzt$ are alternating $0$-compositions  of $2^p x\beta$ and $2^{p+1}\alpha t$, whose targets are all the $0$-cells $2^p xzt$ and $2^{p+1}xyt$, for~$p\geq 1$. Similarly, the only positive $1$-cells of source~$3xyt$ have the $0$-cells $3.2^p xyt$ and $3.2^p xzt$ as targets, for~$p\geq 0$. The other possible $1$-cells of source~$2xzt$ and~$3xyt$ are not positive, like the dotted ones. Here, it is the termination hypothesis that fails, as testified by the infinite sequences of rewriting steps in the previous diagram.

On the other hand, the lack of critical confluence may imply that some nonoverlapping local branchings are not confluent, even under the hypothesis of termination. For example, the $1$-polygraph
\[
\bigpres{x, y, z}{xy\ofl{\alpha} 2x, \, yz\ofl{\beta} z}
\]
terminates, but it has a nonconfluent Peiffer branching:
\[
\xymatrix @R=1em {
& 6xz
& 3xz
\\
& 3xyz
	\ar [u] ^-*+{3\alpha z}
	\ar [ur] ^-{3x\beta}
	\ar@{.>} [dr] ^(0.55)*+{2x\beta + xyz}
\\
xyyz + xyz
	\ar@/^/ [ur] ^(0.4)*+{\alpha yz + xyz}
	\ar@/_/ [dr] _(0.4)*+{xy\beta + xyz}
	\ar@{} [rr] |-{\sm =}
&& 2xz + xyz
\\
& 2xyz
	\ar@{.>} [ur] _(0.55)*+{\alpha z + xyz}
	\ar [dr] _-{2\alpha z}
	\ar [d] _-*+{2x\beta}
\\
& 4xz
& 2xz
}
\]
Here, it is the hypothesis on confluence of critical branchings that is not satisfied, since the critical branching~$(\alpha z, x\beta)$ of source~$xyz$ is not confluent. As a consequence, the only $1$-cells that would close the confluence diagram of the Peiffer branching are the dotted ones, which are not positive.
\end{remark}

\subsection{Squier's theorem}
\label{SS:Squier}

\begin{proposition}
\label{P:PreSquier}
Let~$X$ be a left-monomial $1$-polygraph, and~$Y$ be a cellular extension of~$\lin{X}$. If~$X$ is $Y$-convergent, then~$Y$ is acyclic.
\end{proposition}

\begin{proof}
Since~$X$ is $Y$-convergent, it is convergent, so every $0$-cell~$a$ of~$\lin{X}$ admits a unique normal form~$\rep{a}$, and~$\lin{X}$ contains a positive $1$-cell
\[
a \ofl{\eta_a} \rep{a}.
\]

Now, let $f:a\fl b$ be a positive $1$-cell of~$\lin{X}$. Since~$f$, $\eta_a$ and~$\eta_b$ are positive $1$-cells of~$\lin{X}$, the pair $(f\star_0\eta_b,\eta_a)$ is a branching of~$X$. By hypothesis, this branching is $Y$-confluent, so that, using the fact that~$\rep{a}$ and~$\rep{b}$ are reduced $0$-cells of~$\lin{X}$ that are necessarily equal, we get a $2$-cell
\[
\xymatrix @R=0.5em {
& b
	\ar@/^/ [dr] ^-{\eta_b}
	\ar@2 []!<0pt,-10pt>;[d]!<0pt,0pt> ^-*+{\eta_f}
\\
a
	\ar@/^/ [ur] ^-{f}
	\ar@/_/ [rr] _-{\eta_a}
&& {\rep{a}}
}
\]
in~$\lin{X}[Y]$. Put $\eta_{f^-}=f^-\star_0\eta_f^-$ to obtain the following $2$-cell of~$\lin{X}[Y]$:
\[
\xymatrix @R=0.5em {
& a
	\ar@/^/ [dr] ^-{\eta_a}
	\ar@2 []!<0pt,-10pt>;[d]!<0pt,0pt> ^-*+{\eta_{f^-}}
\\
b
	\ar@/^/ [ur] ^-{f^-}
	\ar@/_/ [rr] _-{\eta_b}
&& {\rep{a}}
}
\]

Next, let~$f:a\fl b$ be any $1$-cell of~$\lin{X}$. By Lemmas~\ref{L:NCellDecomposition} and~\ref{L:FactElem1Cell}, the $1$-cell~$f$ factorises into
\[
f = g_1 \star_0 h_1^- \star_0 \dots \star_0 g_p \star_0 h_p^-,
\]
where~$g_1$, \dots, $g_p$ and~$h_1$, \dots, $h_p$ are positive $1$-cells of~$\lin{X}$. Then define~$\eta_f$ as the following composite $2$-cell of~$\lin{X}[Y]$, with source $f\star_0\eta_b$ and target~$\eta_a$:
\[
\xymatrix @!C @C=3.5em @R=2.5em {
a
	\ar [r] ^-{g_1} 
	\ar [d] _-{\eta_{a}}
& b_1
	\ar [r] ^-{h_1^-}
	\ar [d] |-*+{\eta_{b_1}}
	\ar@2 []!<-30pt,-18pt>;[dl]!<20pt,18pt> ^-{\eta_{g_1}}
& a_2
	\ar [r] 
	\ar [d] |-*+{\eta_{a_2}}
	\ar@2 []!<-30pt,-18pt>;[dl]!<20pt,18pt> ^-{\eta_{h_1^-}}
& \cdots
	\ar [r] 
& a_p
	\ar [r] ^-{g_p}
	\ar [d] |-*+{\eta_{a_p}}
& b_p
	\ar [r] ^-{h_p^-}
	\ar [d] |-*+{\eta_{b_p}}
	\ar@2 []!<-30pt,-18pt>;[dl]!<20pt,18pt> ^-{\eta_{g_p}}
& b
	\ar [d] ^-{\eta_b}
	\ar@2 []!<-30pt,-18pt>;[dl]!<20pt,18pt> ^-{\eta_{h_p^-}}
\\
{\rep{a}}
	\ar@{=} [r]
& {\rep{a}}
	\ar@{=} [r]
& {\rep{a}}
	\ar@{=} [r]
& \cdots 
	\ar@{=} [r]
& {\rep{a}}
	\ar@{=} [r]
& {\rep{a}}
	\ar@{=} [r]
& {\rep{a}}
}
\]

Finally, for all parallel $1$-cells $f,g:a\fl b$ of~$\lin{X}$, the composite $2$-cell
\[
\xymatrix @R=1.5em @C=3em{
& b
	\ar [dr] |-*+{\eta_b}
	\ar@/^2.5ex/ [drr] _-{}="src1" ^(0.66){1_b}
	\ar@2 []!<-10pt,-12.5pt>;[d]!<-10pt,10pt> ^-*+{\eta_f}
	\ar@{} "src1"!<-5pt,0pt>;[dr] |-{\sm =}
\\
a
	\ar@/^/ [ur] ^-{f}
	\ar [rr] |-*+{\eta_a}
	\ar@/_/ [dr] _-{g}
&& {\rep{a}}
	\ar [r] |-*+{\eta_b^-}
& b
\\
& b
	\ar [ur] |-*+{\eta_b}
	\ar@/_2.5ex/ [urr] ^-{}="src2" _(0.66){1_b}
	\ar@2 [u]!<-10pt,-10pt>;[]!<-10pt,12.5pt> ^-*+{\eta_g^-}
	\ar@{} "src2"!<-5pt,0pt>;[ur] |(0.4){\sm =}
}
\]
of~$\lin{X}[Y]$ has source~$f$ and target~$g$, thus concluding the proof that~$Y$ is acyclic. 
\end{proof}

Composing Theorem~\ref{T:HCritical} and Propositions~\ref{P:HNewman} and~\ref{P:PreSquier} gives the analogue of Squier's theorem~\cite[Theorem 5.2]{Squier94} for convergent left-monomial $1$-polygraphs:

\begin{theorem}[Squier's theorem]
\label{T:Squier}
Let~$X$ be a convergent left-monomial $1$-polygraph. A cellular extension~$Y$ of~$\lin{X}$ that contains a $2$-cell
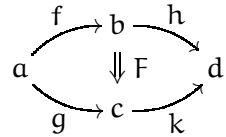
\[
\vcenter{\xymatrix @R=0.25em {
& b
	\ar@/^/ [dr] ^-{h}
	\ar@2 []!<0pt,-10pt>;[dd]!<0pt,10pt> ^-*+{F}
\\
a
	\ar@/^/ [ur] ^-{f}
	\ar@/_/ [dr] _-{g}
&& d
\\
& c
	\ar@/_/ [ur] _-{k}
}}
\]
for every critical branching~$(f,g)$ of~$X$, with~$f'$ and~$g'$ positive $1$-cells of~$\lin{X}$, is acyclic.
\end{theorem}

\begin{example}
\label{X:PP05}

From~\cite[4.3]{PP05}, we consider the quadratic algebra~$A$ presented by 
\[
\bigpres{x,y,z}{x^2 + yz=0,\; x^2 + \lambda zy=0},
\]
where~$\lambda$ is a fixed scalar different from $0$ and $1$. Put~$\mu=\lambda^{-1}$. The algebra~$A$ admits the presentation 
\[
X = 
	\bigpres
		{x, y, z}
		{yz \ofl{\alpha} -x^2 ,\; zy \ofl{\beta} - \mu x^2}.
\]
The deglex order generated by $z>y>x$ satisfies $yz>x^2$ and $zy>x^2$, proving that~$X$ terminates. However, $X$ is not confluent. Indeed, it has two critical branchings:
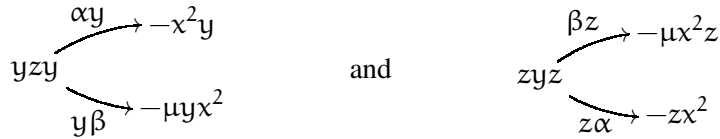
\[
\vcenter{\xymatrix @R=0em {
& -x^2y
\\
yzy
	\ar @/^/ [ur] ^-{\alpha y}
	\ar @/_/ [dr] _-{y\beta}
\\
& - \mu yx^2
}}
\qquad\qquad\text{and}\qquad\qquad
\vcenter{\xymatrix @R=0.5em @C=2em @R=0em {
& - \mu x^2z
\\
zyz
	\ar @/^/ [ur] ^-{\beta z}
	\ar @/_/ [dr] _-{z\alpha}
\\
& -zx^2
}}
\]
and neither of them is confluent, because the monomials~$x^2 y$, $yx^2$, $x^2z$ and~$zx^2$ are reduced.
The adjunction of the $1$-cells 
\[
yx^2 \:\ofl{\gamma}\: \lambda x^2y
\qquad\text{and}\qquad
zx^2 \:\ofl{\delta}\: \mu x^2z
\]
gives a left-monomial $1$-polygraph~$Y$ that also presents~$A$, since~$\gamma$ and~$\delta$ induce relations that already hold in~$\cl{X}$, and that also terminates, because of $yx^2>x^2y$ and $zx^2>x^2z$. Moreover, each one of the four critical branchings of~$Y$ is confluent:
\[
\vcenter{\xymatrix @R=0.75em @C=1.25em {
& -x^2y
	\ar@2 []!<-20pt,-20pt>;[dd]!<-20pt,20pt> ^-*+{F}
\\
yzy
	\ar @/^/ [ur] ^-{\alpha y}
	\ar @/_/ [dr] _-{y\beta}
\\
& - \mu yx^2
	\ar @/_/ [uu] _-{-\mu\gamma}
}}
\qquad
\vcenter{\xymatrix @R=0.75em @C=1.25em {
& -\mu x^2z
	\ar@2 []!<-20pt,-20pt>;[dd]!<-20pt,20pt> ^-*+{G}
\\
zyz
	\ar @/^/ [ur] ^-{\beta z}
	\ar @/_/ [dr] _-{z\alpha}
\\
& -zx^2
	\ar @/_/ [uu] _-{-\delta}
}}
\qquad
\vcenter{\xymatrix @R=0.75em @C=0.25em {
& -x^4
	\ar@2 []!<-5pt,-22.5pt>;[dd]!<-5pt,22.5pt> ^-*+{H}
\\
yzx^2
	\ar@/^/ [ur] ^-{\alpha x^2}
	\ar@/_/ [dr] _-{y\delta}
&& x^2yz
	\ar@/_/ [ul] _-{x^2\alpha}
\\
& \mu yx^2z
	\ar@/_/ [ur] _-{\mu \gamma z}
}}
\qquad
\vcenter{\xymatrix @R=0.75em @C=0.25em {
& - \mu x^4
	\ar@2 []!<-5pt,-22.5pt>;[dd]!<-5pt,22.5pt> ^-*+{I}
\\
zyx^2
	\ar@/^/ [ur] ^-{\beta x^2}
	\ar@/_/ [dr] _-{z\gamma}
&& x^2zy
	\ar@/_/ [ul] _-{x^2\beta}
\\
& \lambda zx^2y
	\ar@/_/ [ur] _-{\lambda \delta y}
}}
\]
Theorem~\ref{T:Squier} implies that the $2$-polygraph $\pres{Y}{F,G,H,I}$ is a coherent presentation of~$A$.

This coherent presentation can be reduced to a smaller one by a collapsing mechanism, that is formalised in Theorem~\ref{T:Collapsing}, but hinted at on this example. First, some $2$-cells may be removed without breaking acyclicity, because their boundary can also be filled by a composite of other $2$-cells. Here, the ``critical $3$-branchings'', where three rewriting steps overlap, reveal two relations between $2$-cells:
\[
\vcenter{\xymatrix @!C @C=3em {
& -x^2yz
	\ar@/^3ex/ [dr] ^-{-x^2\alpha}
	\ar@2 []!<-40pt,-20pt>;[d]!<-40pt,15pt> ^-*+{Fz}
	\ar@2 []!<32.5pt,-37.5pt>;[dd]!<32.5pt,37.5pt> ^-*+{-H}
\\
yzyz
	\ar@/^3ex/ [ur] ^-{\alpha yz}
	\ar [r] |-*+{y\beta z}
	\ar@/_3ex/ [dr] _-{yz\alpha}
& -\mu yx^2z
	\ar [u] _-{-\mu\gamma z}
	\ar@2 []!<-40pt,-15pt>;[d]!<-40pt,20pt> ^-*+{yG}
& x^4
\\
& -yzx^2
         \ar [u] _-{-y\delta}
	\ar@/_3ex/ [ur] _-{-\alpha x^2} 
}}
\qquad\tfl\qquad
\vcenter{\xymatrix @!C @C=1.5em @R=1.5em {
& -x^2yz
        \ar @/^2ex/ [dr] ^{-x^2\alpha}
\\
yzyz
	\ar @/^2ex/ [ur] ^{\alpha yz}
	\ar @/_2ex/ [dr] _{y z\alpha}
	\ar@{} [rr] |-{\sm =}
&& x^4
\\
& -yzx^2
	\ar @/_2ex/ [ur] _{-\alpha x^2}
}}
\]
\[
\vcenter{\xymatrix @!C @C=3em{
& -\mu x^2zy
	\ar@/^3ex/ [dr] ^-{-\mu x^2\beta}
	\ar@2 []!<-40pt,-20pt>;[d]!<-40pt,15pt> ^-*+{Gy}
	\ar@2 []!<32.5pt,-37.5pt>;[dd]!<32.5pt,37.5pt> ^-*+{-\mu I}
\\
zyzy
	\ar@/^3ex/ [ur] ^-{\beta zy}
	\ar [r] |-*+{z\alpha y}
	\ar@/_3ex/ [dr] _-{zy\beta}
& -zx^2y
	\ar [u] _-{-\delta y}
	\ar@2 []!<-40pt,-15pt>;[d]!<-40pt,20pt> ^-*+{zF}
& \mu^2 x^4
\\
& -\mu zyx^2
	\ar [u] _-{-\mu z\gamma}
	\ar@/_3ex/ [ur] _-{-\mu\beta x^2}
}}
\qquad\tfl\qquad
\vcenter{\xymatrix @!C @C=1.5em @R=1.5em {
& bx^2zy
	\ar @/^2ex/ [dr] ^{-\mu x^2\beta}
\\
zyzy
	\ar @/^2ex/ [ur] ^{\beta zy}
	\ar @/_2ex/ [dr] _{zy \beta}
	\ar@{} [rr] |-{\sm =}
&& \mu^2 x^4
\\
& bzyx^2
	\ar @/_2ex/ [ur] _{-\mu\beta x^2}
}}
\]
Since the boundaries of~$H$ and~$I$ can also be filled using~$F$ and~$G$ only, the $2$-polygraph $\pres{Y}{F,G}$ is also a coherent presentation of~$A$. Next, the $1$-cells~$\gamma$ and~$\delta$ are redundant, because the corresponding relations can be derived from~$\alpha$ and~$\beta$, as testified by the $2$-cells~$F$ and~$G$: removing~$\gamma$ with~$F$, and~$\delta$ with~$G$, yields the smaller coherent presentation $\pres{X}{\emptyset}$ of~$A$.
\end{example}

\begin{example}[The standard coherent presentation]
\label{X:StandardCoherentPresentation}
Assume that $A=\K\oplus A_+$ is an augmented algebra, and fix a linear basis~$\Br$ of~$A_+$. 
For~$u$ and~$v$ in~$\Br$, write~$u\tens v$ for the product of~$u$ and~$v$ in the free algebra over~$\Br$, and~$uv$ for their product in~$A$. 
Consider the $1$-polygraph~$\Std(\Br)_1$ whose $0$-cells are the elements of~$\Br$, and with a $1$-cell
\[
u\tens v \ofl{u\vert v} uv,
\]
for all~$u$ and~$v$ in~$\Br$. Note that~$uv$ belongs to the free algebra over~$\Br$ because~$A$ is augmented. By definition, $\Std(\Br)_1$ is a presentation of~$A$. Moreover, $\Std(\Br)_1$ terminates by a length argument: for all~$u$ and~$v$ in~$\Br$, the monomial $u\tens v$ is a word of length~$2$ in the free monoid over~$\Br$, while~$uv$ is a word of length~$1$. Finally, $\Std(\Br)_1$ has one critical branching $(u\vert v\tens w, u\tens v\vert w)$ for each triple~$(u,v,w)$ of elements of~$\Br$, and this critical branching is confluent. Thus, extending $\Std(\Br)_1$ with a $2$-cell
\[
\xymatrix @R=0.75em @!C @C=0em {
& uv\tens w
	\ar@/^2ex/ [dr] ^-{uv\vert w}
	\ar@2 []!<-10pt, -15pt>;[dd]!<-10pt, 15pt> ^-*+{u\vert v\vert w}
\\
u\tens v\tens w
	\ar@/^2ex/ [ur] ^-{u\vert v \tens w}
	\ar@/_2ex/ [dr] _-{u\tens v\vert w}
&& uvw
\\
& u\tens vw
	\ar@/_2ex/ [ur] _-{u\vert vw}
}
\]
for each triple~$(u,v,w)$ of elements of~$\Br$ produces, by Theorem~\ref{T:Squier}, a coherent presentation of~$A$, denoted by $\Std(\Br)_2$. Note that the free $2$-algebra over $\Std(\Br)_2$ does not depend (up to isomorphism) on the choice of the basis~$\Br$.

This coherent presentation of~$A$ is extended in every dimension in Subsection~\ref{SS:StandardPolygraphicResolution} to obtain a polygraphic version of the standard resolution of an algebra. As in the previous example, the next dimension contains the $3$-cells generated by the ``critical $3$-branchings'' of~$\Std_1(\Br)$: there is one such $3$-cell~$u\vert v\vert w\vert x$ for each quadruple~$(u,v,w,x)$ of elements of~$\Br$, with source
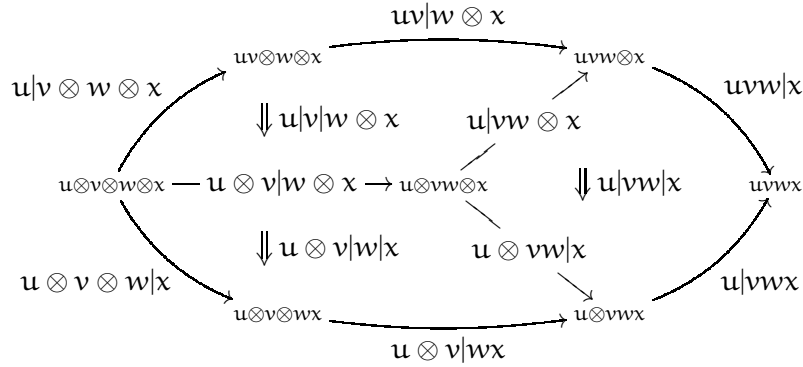
\[
\xymatrix @R=3em @C=1.5em @!C {
& {\sm uv \tens w\tens x}
	\ar@/^/ [rr] ^-{uv\vert w\tens x}
	\ar@2 []!<-5pt,-18pt>;[d]!<-5pt,18pt> ^-{\:u\vert v\vert w\tens x}
&& {\sm uvw\tens x}
	\ar@/^2ex/ [dr] ^-{uvw\vert x}
	\ar@2 []!<-10pt,-42.5pt>;[dd]!<-10pt,42.5pt> ^-{\:u\vert vw\vert x}
\\
{\sm u\tens v\tens w\tens x}
	\ar@/^2ex/ [ur] ^-{u\vert v\tens w\tens x}
	\ar [rr] |-*+{u\tens v\vert w\tens x}
	\ar@/_2ex/ [dr] _-{u\tens v\tens w\vert x}
& \strut
	\ar@2 []!<-5pt,-18pt>;[d]!<-5pt,18pt> ^-{\:u\tens v\vert w\vert x}
& {\sm u\tens vw\tens x}
	\ar [ur] |-*+{u\vert vw\tens x}
	\ar [dr] |-*+{u\tens vw\vert x}
&& {\sm uvwx}
\\
& {\sm u\tens v\tens wx}
	\ar@/_/ [rr] _-{u\tens v\vert wx}
&& {\sm u\tens vwx}
	\ar@/_2ex/ [ur] _-{u\vert vwx}
}
\]
and target
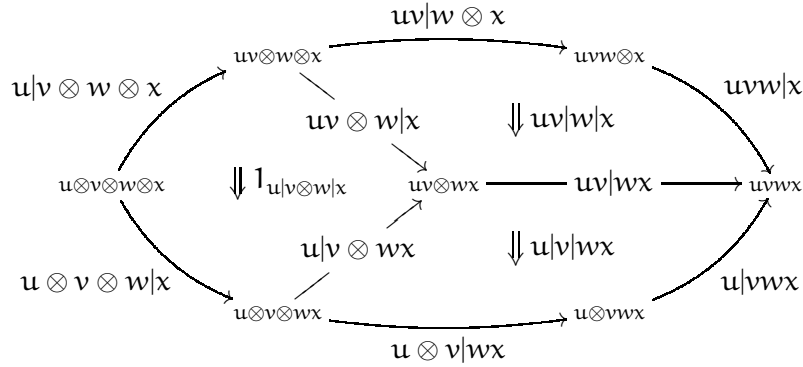
\[
\xymatrix @R=3em @C=1.5em @!C {
& {\sm uv\tens w\tens x}
	\ar@/^/ [rr] ^-{uv\vert w\tens x}
	\ar [dr] |-*+{uv\tens w\vert x}
	\ar@2 []!<-15pt,-42.5pt>;[dd]!<-15pt,42.5pt> ^-{\:1_{u\vert v\tens w\vert x}}
&& {\sm uvw\tens x}
	\ar@/^2ex/ [dr] ^-{uvw\vert x}
	\ar@2 []!<-35pt,-18pt>;[d]!<-35pt,18pt> ^-{\:uv\vert w\vert x}
\\
{\sm u\tens v\tens w\tens x}
	\ar@/^2ex/ [ur] ^-{u\vert v\tens w\tens x}
	\ar@/_2ex/ [dr] _-{u\tens v\tens w\vert x}
&& {\sm uv\tens wx}
	\ar [rr] |-*+{uv\vert wx}
& \strut
	\ar@2 []!<-35pt,-18pt>;[d]!<-35pt,18pt> ^-{\:u\vert v\vert wx}
& {\sm uvwx}
\\
& {\sm u\tens v\tens wx}
	\ar [ur] |-*+{u\vert v\tens wx}
	\ar@/_/ [rr] _-{u\tens v\vert wx}
&& {\sm u\tens vwx}
	\ar@/_2ex/ [ur] _-{u\vert vwx}
}
\]
\end{example}

\section{Construction and reduction of polygraphic resolutions}
\label{S:PolygraphicResolutions}

In this section, we introduce concepts and techniques that are used in Section~\ref{S:SquierResolution} to construct Squier's polygraphic resolution. 
First, we define contractions as sort of homotopy, and prove in Theorem~\ref{T:ResolutionContraction} that they characterise $\infty$-polygraphs that are polygraphic resolutions.
Then, following~\cite{Brown92}, we define collapsing schemes of polygraphs, and Theorem~\ref{T:Collapsing} shows that they give a method to contract polygraphic resolutions into smaller ones. This process, similar to algebraic Morse theory for chain complexes~\cite{Skoldberg06}, was already used to obtain minimal coherent presentations of Artin monoids in~\cite{GaussentGuiraudMalbos}.

\subsection{Linear homotopies}
\label{SS:Homotopies}

\subsubsection{Homotopies}

Let~$A$ and~$B$ be $\infty$-vector spaces (resp.\ $\infty$-algebras), and $F,G:A\fl B$ be linear $\infty$-functors (resp.\ morphisms of $\infty$-algebras). A \emph{homotopy from~$F$ to~$G$} is an indexed linear map (resp.\ indexed morphism of algebras)
\[
A \ofl{\eta} B
\]
of degree~$1$ that satisfies, writing~$\eta_x$ for~$\eta(x)$,
\begin{enumerate}
\item for every~$0$-cell~$x$ of~$A$, 
\begin{equation}
\label{E:SourceTargetHomotopy0}
s(\eta_x) = F(x)
\qquad\text{and}\qquad
t(\eta_x) = G(x),
\end{equation}
\item for every~$n\geq 1$ and every $n$-cell~$x$ of~$A$,
\begin{align}
\label{E:SourceHomotopy}
s(\eta_x) &= F(x) \star_0 \eta_{t_0(x)} \star_1 \cdots \star_{n-1} \eta_{t_{n-1}(x)}, 
\\[0.5ex]
\label{E:TargetHomotopy}
t(\eta_x) &= \eta_{s_{n-1}(x)} \star_{n-1} \cdots \star_1 \eta_{s_0(x)} \star_0 G(x),
\end{align}
with parentheses omitted according to the convention that~$\star_i$ binds more tightly than~$\star_j$ if~$i<j$, 
\item for every~$n\geq 0$ and every $n$-cell~$x$ of~$A$,
\begin{equation}
\label{E:HomotopyIdentity}
\eta_{1_x} = 1_{\eta_x}.
\end{equation}
 \end{enumerate}
 
\subsubsection{Remarks}

\Item{i}
Expanding~\eqref{E:SourceHomotopy} and~\eqref{E:TargetHomotopy}, a homotopy~$\eta$ from~$F$ to~$G$ maps a $1$-cell~$x:y\fl z$ to a $2$-cell
\[
\xymatrix @R=0.25em @C=1.5em {
& F(z)
	\ar@/^/ [dr] ^-{\eta_z}
	\ar@2 []!<0pt,-15pt>;[dd]!<0pt,15pt> ^-*+{\eta_x}
\\
F(y)
	\ar@/^/ [ur] ^-{F(x)}
	\ar@/_/ [dr] _-{\eta_y}
&& G(z)
\\
& G(y)
	\ar@/_/ [ur] _-{G(x)}
}
\]
and a $2$-cell $x:y\dfl y':z\fl z'$ of~$A$ to a $3$-cell
\[
\vcenter{\xymatrix @R=0.25em @C=4em {
& F(z')
	\ar@/^/ [dr] ^-{\eta_{z'}}
		\ar@2 []!<5pt,-15pt>;[dd]!<5pt,15pt> ^-*+{\eta_{y'}}
\\
F(z)
	\ar@/^5ex/ [ur] ^(0.2){F(y)} _{}="s"
	\ar@/_/ [ur] |-*+{F(y')} ^{}="t"
		\ar@2 "s"!<-6pt,-10pt>;"t"!<-9pt,10pt> ^-{F(x)}
	\ar@/_/ [dr] _{\eta_z}
&& {G(z')} 
\\
& G(z)
	\ar@/_/ [ur] _-{G(y')}
}}
\quad\otfl{\eta_x}\quad
\vcenter{\xymatrix @R=0.25em @C=4em  {
& F(z')
	\ar@/^/ [dr] ^-{\eta_{z'}}
		\ar@2 []!<-10pt,-15pt>;[dd]!<-10pt,15pt> ^-*+{\eta_{y}}
\\
F(z)
	\ar@/^/ [ur] ^-{F(y)}
	\ar@/_/ [dr] _{\eta_z}
&& {G(z')} 
\\
& G(z)
	\ar@/^/ [ur] |-*+{G(y)} _-{}="s"
	\ar@/_5ex/ [ur] _(0.8){G(y')} ^{}="t"
		\ar@2 "s"!<-7pt,-12pt>;"t"!<-10pt,7pt> ^-{G(x)}
}}
\]

\Item{ii}
By definition, a homotopy is linear. Together with the relation $x\star_k y = x - t_k(x) + y$, this implies that a homotopy also satisfies the same compatibility condition with respect to compositions as the set-theoretic homotopies of~\cite[B.8]{AraMaltsiniotis15}.

\begin{lemma}
\label{L:Homotopy}
Let~$A$ and~$B$ be $\infty$-vector spaces (resp.\ $\infty$-algebras) and $F,G:A\fl B$ be linear $\infty$-functors (resp.\ morphisms of $\infty$-algebras).
\begin{enumerate}
\item The definition of a homotopy~$\eta$ from~$F$ to~$G$ makes sense: the right-hand sides of~\eqref{E:SourceHomotopy} and~\eqref{E:TargetHomotopy} are well defined, the assignments of~$x$ to~$s(\eta_x)$ and~$t(\eta_x)$ are linear maps (resp.\ morphisms of algebras), $s(\eta_x)$ and~$t(\eta_x)$ are parallel, and $s(\eta_{1_x})=t(\eta_{1_x})=\eta_x$ hold.
\item If~$\eta$ is a homotopy from~$F$ to~$G$, then, for every~$n\geq 1$ and every $n$-cell~$x$ of~$A$,
\[
s(\eta_x) = F(x) - 1_{t(F(x))} + \eta_{t(F(x))}
\qquad\text{and}\qquad
t(\eta_x) = G(x) - 1_{t(G(x))} + \eta_{t(G(x))}.
\]
\end{enumerate}
\end{lemma}

\begin{proof}
\Item{i} Fix an $n$-cell~$x$ of~$A$, with~$n\geq 1$. The right-hand side of~\eqref{E:SourceHomotopy} is well defined, since, for every~$0\leq k<n$,
\[
\begin{array}{l l l}
t_k(F(x)\star_0\eta_{t_0(x)} \star_1 \cdots \star_{k-1} \eta_{t_{k-1}(x)}) 
	&\:=\: t_k(F(x)) \star_0\eta_{t_0(x)} \star_1 \cdots \star_{k-1} \eta_{t_{k-1}(x)} \\[1ex]
	&\:=\: F(t_k(x)) \star_0\eta_{t_0t_k(x)} \star_1 \cdots \star_{k-1} \eta_{t_{k-1}t_k(x)} 
	&\:=\: s_k(\eta_{t_k(x)}).
\end{array}
\]
Similarly, the compositions that appear in the right-hand of~\eqref{E:TargetHomotopy} are legitimate. Next, mapping~$x$ to~$s(\eta_x)$ or to~$t(\eta_x)$ is a linear map (resp.\ a morphism of algebras), because each one is a composite of such operations. Moreover, $s(\eta_x)$ and $t(\eta_x)$ are parallel, because the globular relations and the compatibility of the source and target maps with the compositions imply
\begin{align*}
& ss(\eta_x)
	= s(F(x)) \star_0 \eta_{t_0(x)} \star_1 \cdots \star_{n-2} \eta_{t_{n-2}(x)} 
	= s(\eta_{s(x)}) 
	= st(\eta_x)
\\[0.5ex]
\text{and}\quad 
&
ts(\eta_x) 
	= t(\eta_{t(x)})
    = \eta_{s_{n-2}(x)} \star_{n-2} \cdots \star_1 \eta_{s_0(x)} \star_0 t(G(x))
	= tt(\eta_x).
\end{align*}
Finally, replacing~$x$ by~$1_x$ in~\eqref{E:SourceHomotopy} and~\eqref{E:TargetHomotopy} yields $s(\eta_{1_x})=t(\eta_{1_x})=\eta_x$, so that asking $\eta_{1_x}=1_{\eta_x}$ is legal.

\smallskip
\Item{ii} By induction on~$k$, using the relation $x\star_k y=x-1_{t_k(x)}+y$.
\end{proof}

\begin{lemma}
\label{L:FreeHomotopy}
Assume that~$X$ is an $\infty$-polygraph, and that~$A$ is an $\infty$-algebra.

\begin{enumerate}

\item Let~$F,G:\lin{X}\fl A$ be linear $\infty$-functors. A homotopy~$\eta$ from~$F$ to~$G$ is uniquely and entirely determined by its values on the $n$-monomials of~$\lin{X}$, for~$n\geq 0$, provided the relation
\begin{equation}
\label{E:HomotopyExchange}
\eta_{as_0(b)} + \eta_{t_0(a)b} - \eta_{t_0(a)s_0(b)} 
= \eta_{s_0(a)b} + \eta_{at_0(b)} - \eta_{s_0(a)t_0(b)} 
\end{equation}
is satisfied for all $n$-monomials~$a$ and~$b$ of~$\lin{X}$.

\item Let~$F,G:\lin{X}\fl A$ be morphisms of $\infty$-algebras. A homotopy~$\eta$ from~$F$ to~$G$ is uniquely and entirely determined by its values on the cells of~$X$.

\end{enumerate}

\end{lemma}

\begin{proof}
\Item{i}
By induction on~$n\geq 0$. 
For~$n=0$, assume that $\eta_u:F(u)\fl G(u)$ is a fixed $1$-cell of~$A$ for every monomial~$u$ of~$\lin{X}$. Extend~$\eta$ to every $0$-cell~$a$ of~$\lin{X}$ by linearity: the resulting $1$-cell~$\eta_a$ has source~$F(a)$ and~$G(a)$, as required by~\eqref{E:SourceTargetHomotopy0} by linearity of~$F$ and~$G$.

Now, fix~$n\geq 1$, and assume that an $(n+1)$-cell~$\eta_a$ has been chosen in~$A$ for every $n$-monomial~$a$ of~$\lin{X}$, with source and target given by~\eqref{E:SourceHomotopy} and~\eqref{E:TargetHomotopy}, and in such a way that~\eqref{E:HomotopyExchange} holds. By construction, the $n$-cells of~$\lin{X}$ are the linear combinations of the $n$-monomials of~$\lin{X}$ and of identities of $(n-1)$-cells of~$\lin{X}$, up to the relation
\[
a s_0(b) + t_0(a) b - t_0(a) s_0(b) 
	= s_0(a) b + a t_0(b) - s_0(a) t_0(b),
\]
where~$a$ and~$b$ range over the $n$-monomials of~$\lin{X}$. Note that the values of~$\eta$ on identities of $(n-1)$-cells are constrained by~\eqref{E:HomotopyIdentity}. Thus, one can extend~$\eta$ to any $n$-cell~$a$ of~$\lin{X}$, provided a decomposition of~$a$ into a linear combination of $n$-monomials and of an identity is chosen, and~\eqref{E:HomotopyExchange} implies that the result is independent of this choice. Finally, the source and target of the obtained $(n+1)$-cell~$\eta_a$ satisfy~\eqref{E:SourceHomotopy} and~\eqref{E:TargetHomotopy} by linearity of~$F$, $G$, $\eta$ and all the compositions~$\star_0$, \dots, $\star_{n-1}$.

\smallskip\Item{ii}
By definition, a homotopy~$\eta$ between morphisms~$F$ and~$G$ of $\infty$-algebras is a homotopy between the underlying linear $\infty$-functors that commutes with identities and products. Thus, according to~\Item{i}, the values of~$\eta$ on the $0$-cells of~$\lin{X}$ can be uniquely reconstructed from its values on the $0$-cells of~$X$, and the corresponding $1$-cells have the required source and target because~$F$ and~$G$ are morphisms of $\infty$-algebras. 

Now, if~$n\geq 1$, we observe that the fact that~$\eta$ commutes with products ensures that~\eqref{E:HomotopyExchange} is automatically satisfied, so that~$\eta$ can be uniquely reconstructed from its values on the $n$-cells of~$\lin{X}$. Moreover, \eqref{E:SourceHomotopy} and~\eqref{E:TargetHomotopy} are satisfied because~$F$, $G$, $\eta$ and the compositions~$\star_0$, \dots, $\star_{n-1}$ are morphisms of $\infty$-algebras.
\end{proof}

\subsection{Contractions of polygraphs for associative algebras}
\label{SS:Contractions}

\subsubsection{Unital sections and contractions}
\label{SSS:UnitalSectionsContractions}

Let~$X$ be an $\infty$-polygraph. A \emph{unital section of~$X$} is a linear $\infty$-functor $\iota:\cl{X}\fl\lin{X}$ which is a section of the canonical projection $\pi:\lin{X}\pfl\cl{X}$ and that satisfies $\iota(1)=1$. Here the quotient algebra~$\cl{X}$ is seen as an $\infty$-algebra whose $n$-cells are identities for every~$n\geq 1$. Note that a unital section of~$\pi$ is not required to be a morphism of $\infty$-algebras. 

Let~$\iota$ be a unital section of~$X$. If~$a$ is an $n$-cell of~$\lin{X}$, write~$\rep{a}$ for $\iota\pi(a)$ when no confusion occurs (note that~$\rep{a}$ is an identity if~$n\geq 1$). An \emph{$\iota$-contraction of~$X$} is a homotopy~$\sigma$ from~$\id_{\lin{X}}$ to the composite linear $\infty$-functor~$\iota\pi$ that satisfies
\[
\sigma_a = 1_a
\]
for every $n$-cell~$a$ of~$\lin{X}$ that belongs to the image of~$\iota$ or of~$\sigma$. From the definition of~$\iota\pi$, we note that, for every~$n\geq 0$ and every $n$-cell~$a$ of~$\lin{X}$, 
\[
t(\sigma_a) = 
\begin{cases}
\rep{a} 
	&\text{if $n=0$,} \\
\sigma_{s(a)}
	&\text{if $n\geq 1$.}
\end{cases}
\]
An $\iota$-contraction~$\sigma$ of~$X$ is called \emph{right} if, for every~$n\geq 0$ and all $n$-cells~$f$ and~$g$ of~$\lin{X}$ of respective $0$-sources~$a$ and~$b$, it satisfies the relation
\begin{equation}
\label{E:RightContraction}
\sigma_{fg} = a\sigma_g \star_0 \sigma_{f\rep{b}}.
\end{equation}

\subsubsection{Remarks}

Expanding the definition, an $\iota$-contraction~$\sigma$ of~$X$ maps a $0$-cell~$a$ to a $1$-cell ${\sigma_a:a\fl\rep{a}}$, a $1$-cell~$f:a\fl b$ to a $2$-cell
\[
\xymatrix @R=0.5em {
& b
	\ar@/^/ [dr] ^-{\sigma_b}
	\ar@2 []!<0pt,-10pt>;[d] ^-*+{\sigma_f}
\\
a
	\ar@/^/ [ur] ^-{f}
	\ar@/_/ [rr] _{\sigma_a}
&& {\rep{a}}
}
\]
and a $2$-cell $F:f\dfl g:a\fl b$ to a $3$-cell
\[
\vcenter{\xymatrix @R=1em {
&& b
	\ar @/^/ [drr] ^-{\sigma_b}
\\
a
	\ar @/^5ex/ [urr] ^-{f} _{}="s"
	\ar [urr] |-*+{g} ^{}="t"
	\ar @/_/ [rrrr] _{\sigma_a} ^{}="tgt"
&&&& {\rep{a}} 
	\ar@2 "s"!<3pt,-8pt>;"t"!<-3pt,8pt> ^-{F}
	\ar@2 "1,3"!<0pt,-12.5pt>;"tgt"!<0pt,10pt> ^-*+{\sigma_g}
}}
\quad\otfl{\sigma_F}\quad
\vcenter{\xymatrix @R=0.5em {
& b
	\ar@/^/ [dr] ^-{\sigma_b}
	\ar@2 []!<0pt,-10pt>;[d] ^-*+{\sigma_f}
\\
a
	\ar@/^/ [ur] ^-{f}
	\ar@/_/ [rr] _{\sigma_a}
&& {\rep{a}}
}}
\]

If~$\sigma$ is a right $\iota$-contraction, then~\eqref{E:RightContraction} reads, on $0$-cells~$a$ and~$b$, 
\[
\xymatrix @R=0.5em {
& {a\rep{b}}
	\ar @/^/ [dr] ^-{\sigma_{a\rep{b}}}	
	\ar@{} [d] |(0.70){\sm =} 
\\
ab
	\ar @/^/ [ur] ^-{a\sigma_{b}} 
	\ar @/_/ [rr] _-{\sigma_{ab}} 
&& {\rep{ab}}
}
\]
and, on $1$-cells $f:a\fl b$ and $g:c\fl d$,
\[
\vcenter{\xymatrix @R=0.5em {
& bd
	\ar@/^/ [dr] ^-{\sigma_{bd}}
	\ar@2 []!<0pt,-10pt>;[d]!<0pt,2.5pt> ^-*+{\sigma_{fg}}
\\
ac
	\ar@/^/ [ur] ^-{fg}
	\ar@/_/ [rr] _-{\sigma_{ac}}
&& {\rep{ac}}
}}
\qquad = \qquad
\vcenter{\xymatrix @R=1.5em @C = 3.5em {
&& bd
	\ar @/^5ex/ [ddrr] ^-{\sigma_{bd}} _-{}="B"
	\ar [dr] |-*+{b\sigma_{d}}
	\ar@{} [dd] |-{\sm =}
\\
& ad
	\ar [ur] |-*+{fd}
	\ar [dr] |-*+{f\sigma_{d}}
	\ar@2 []!<-5pt,-12pt>;[d]!<-5pt,11pt> ^-*+{a\sigma_g}
&& {b\rep{c}}
	\ar [dr] |-*+{\sigma_{b\rep{c}}}
	\ar@2 []!<-5pt,-12pt>;[d]!<-5pt,11pt> ^-*+{\sigma_{f\rep{c}}}
	\ar@{} [];"B" |(0.6){\sm =}
\\
ac 
	\ar @/^5ex/ [uurr] ^-{fg} _-{}="A"
	\ar@{} [ur];"A" |(0.5){\sm =}
	\ar [ur] |-*+{ag}
	\ar [rr] |-*+{a\sigma_c}
	\ar @/_5ex/ [rrrr] _-{\sigma_{ac}} ^-{}="C"
&& {a\rep{c}}
	\ar [ur] |-*+{f\rep{c}}
	\ar [rr] |-*+{\sigma_{a\rep{c}}}
	\ar@{} [];"C" |-{\sm =}
&& {\rep{ac}}
}}
\]

\begin{lemma}
\label{L:Contraction}
Let~$X$ be an $\infty$-polygraph, and~$\iota$ be a unital section of~$X$.
\begin{enumerate}
\item The notion of (right) $\iota$-contraction of~$X$ makes sense: if~$a$ is in the image of~$\iota$ or  of~$\sigma$, then $s(\sigma_a)=t(\sigma_a)=a$, and the composition of the right-hand side of~\eqref{E:RightContraction} is well defined.
\item If~$\sigma$ is a right $\iota$-contraction of~$X$, then, for every~$n\geq 0$ and all $n$-cells~$f$ and~$g$ of~$\lin{X}$, 
\[
\sigma_{fg} = f\sigma_g - f\rep{b} + \sigma_{f\rep{b}} \,,
\qquad\text{with $b=s_0(g)$.}
\]
\end{enumerate}
\end{lemma}

\begin{proof}
\Item{i} Assume that~$a$ is a $0$-cell of~$\lin{X}$. By definition of an $\iota$-contraction, we have $s(\sigma_a)=a$ and $t(\sigma_a)=\rep{a}$. Thus, if~$a$ is in the image of~$\iota$, then~$\rep{a}=a$, so that~$\sigma_a$ is parallel to~$1_a$. 

Then, fix~$n\geq 1$ and assume that~$s(\sigma_a)=t(\sigma_a)=a$ hold for every $k$-cell~$a$ of~$\lin{X}$, with~$k<n$, that lies in the image of~$\iota$ or of~$\sigma$. Fix an $n$-cell~$a$ of~$\lin{X}$. If~$a$ is in the image of~$\iota$, then~$a$ is an identity~$1_b$, with~$b$ in the image of~$\iota$ by functoriality: because~$\sigma$ is a homotopy, it satisfies~$\sigma_a=1_{\sigma_b}$, so that, by induction hypothesis, $\sigma_a=1_{1_b}=1_a$. Now, assume that~$a=\sigma_b$, for some $(n-1)$-cell~$b$ of~$\lin{X}$. Lemma~\ref{L:Homotopy} gives
\[
s(\sigma_a) = \sigma_b - 1_{t(\sigma_b)} + \sigma_{t(\sigma_b)}.
\]
We distinguish two cases: if $n=1$, then $t(\sigma_b)=\rep{b}$, and, if~$n>1$, then $t(\sigma_b)=\sigma_{s(b)}$. Either way, $t(\sigma_b)$ is an $(n-1)$-cell of~$\lin{X}$ that lies in the image of~$\iota$ or of~$\sigma$, so that the induction hypothesis implies $\sigma_{t(\sigma_b)}=1_{t(\sigma_b)}$, and thus $s(\sigma_a)=\sigma_b=a$. Moreover, $t(\sigma_a)=\sigma_{s(a)}$, and Lemma~\ref{L:Homotopy} implies
\[
t(\sigma_a) 
	= \sigma_{s(\sigma_b)}
	= \sigma_b - \sigma_{1_{t(b)}} + \sigma_{\sigma_{t(b)}}.
\]
By definition of~$\sigma$, we have $\sigma_{1_{t(b)}} = 1_{\sigma_{t(b)}}$, and, by induction hypothesis,  $\sigma_{\sigma_{t(b)}} = 1_{\sigma_{t(b)}}$, giving $t(\sigma_a)=\sigma_b=a$.

Now, assume that~$f$ and~$g$ are $n$-cells of~$\lin{X}$, of respective $0$-sources~$a$ and~$b$. We have
\[
t_0(a\sigma_g) = a\rep{b} = s_0(\sigma_{f\rep{b}})
\]
so the composition of the right-hand side of~\eqref{E:RightContraction} is well defined. 

\smallskip
\Item{ii} On the one hand, \eqref{E:RightContraction} implies
\[
\sigma_{fg} = a\sigma_g\star_0\sigma_{f\rep{b}} = a\sigma_g - a\rep{b} + \sigma_{f\rep{b}} \,,
\]
and, on the other hand, the fact that~$\star_0$ is a morphism of algebras yields
\[
f\sigma_g = a\sigma_g \star_0 f\rep{b} = a\sigma_g - a\rep{b} + f\rep{b}.
\]
Combining these two expressions gives the result.
\end{proof}

\subsubsection{Reduced and essential monomials}

Fix a unital section~$\iota$ of~$X$. A monomial~$u$ of~$\lin{X}$ is \emph{$\iota$-reduced} if $u=\rep{u}$ holds. A non-$\iota$-reduced monomial~$u$ of~$\lin{X}$ is \emph{$\iota$-essential} if $u=xv$, with~$x$ a $0$-cell of~$X$ and~$v$ an $\iota$-reduced monomial of~$\lin{X}$.

If~$\sigma$ is an $\iota$-contraction of~$X$, for~$n\geq 1$, an $n$-cell~$a$ of~$\lin{X}$ is \emph{$\sigma$-reduced} if it is an identity or in the image of~$\sigma$. If~$\sigma$ is a right $\iota$-contraction of~$X$, and~$n\geq 1$, a non-$\sigma$-reduced $n$-monomial~$a$ of~$\lin{X}$ is \emph{$\sigma$-essential} if there exist an $n$-cell~$\alpha$ of~$X$ and an $\iota$-reduced monomial of~$\lin{X}$ such that $a=\alpha v$.

\begin{lemma}
\label{L:FreeRightContraction}
Let~$X$ be an $\infty$-polygraph, and~$\iota$ be a unital section of~$X$. A right $\iota$-contraction~$\sigma$ of~$X$ is uniquely and entirely determined by its values on the $\iota$-essential monomials and, for every~$n\geq 1$, on the $\sigma$-essential $n$-monomials of~$\lin{X}$. 
\end{lemma}

\begin{proof}
By Lemma~\ref{L:FreeHomotopy}~\Item{i}, we know that the homotopy underlying~$\sigma$ is uniquely and entirely determined by its values on the $n$-monomials of~$\lin{X}$, provided~\eqref{E:HomotopyExchange} is satisfied. There remains to check that the values of~$\sigma$ on the $\iota$-essential monomials and on the $\sigma$-essential $n$-monomials completely determines its values on the other monomials and $n$-monomials, and that~\eqref{E:HomotopyExchange} is automatically satisfied.

If~$u$ is a non-$\iota$-essential monomial, then either $u=1$, or $u=xv$ with $x$ a $0$-cell of~$X$ and~$v$ a non-$\iota$-reduced monomial. In the former case, $\sigma_1=1$ is forced because~$1$ is $\iota$-reduced. In the latter case, \eqref{E:RightContraction} imposes
\[
\sigma_{xv} = x\sigma_v \star_0 \sigma_{x\rep{v}}.
\]
Then proceed by induction on the length of~$v$ to define~$\sigma_v$ from the values of~$\sigma$ on $\iota$-reduced monomials.

Now, for every $n$-monomial~$u\alpha v$ of~$\lin{X}$, with~$\alpha$ an $n$-cell of~$X$ of $0$-source~$a$, and~$u$ and~$v$ monomials of~$\lin{X}$, 
\[
\sigma_{u\alpha v} = ua\sigma_v \star_0 u\sigma_{\alpha\rep{v}} \star_0 \sigma_{u \rep{av}}
\]
is imposed by~\eqref{E:RightContraction}. Then, if~$\alpha v$ is $\sigma$-reduced, then $\alpha v=\sigma_b$ for some $(n-1)$-cell~$b$ of~$\lin{X}$, which implies $\sigma_{\alpha v}=1_{\sigma_b}$.

To check~\eqref{E:HomotopyExchange}, fix $n$-monomials~$a$ and~$b$ of~$\lin{X}$, and put~$c=s_0(a)$, $c'=t_0(a)$, $d=s_0(b)$ and~$d'=t_0(b)$. On the one hand, 
\begin{align*}
\sigma_{ad} + \sigma_{c'b} - \sigma_{c'd} 
	&= c\sigma_{d} \star_0 \sigma_{a\rep{d}}
		+ c'\sigma_{b} \star_0 \sigma_{c'\rep{d}}
		- c'\sigma_{d} \star_0 \sigma_{c'\rep{d}} \\
	&= \sigma_{a\rep{d}} + c'\sigma_b
		+ c\sigma_d - c'\sigma_d - c\rep{d},
\end{align*}
and, on the other hand,
\begin{align*}
\sigma_{cb} + \sigma_{ad'} - \sigma_{cd'}
	&= c\sigma_b \star_0 \sigma_{c\rep{d}}
		+ c\sigma_{d'} \star_0 \sigma_{a\rep{d}}
		- c\sigma_{d'} \star_0 \sigma_{c\rep{d}} \\
	&= \sigma_{a\rep{d}} + c\sigma_b - c\rep{d}.
\end{align*}
So, we are left with proving 
\begin{equation}
\label{E:ExchangeC}
c'\sigma_b + c\sigma_d = c\sigma_b + c'\sigma_d
\end{equation}
under the hypothesis that there exists an $n$-cell~$a$ in~$\lin{X}$ with $0$-source~$c$ and $0$-target~$c'$. Since the $0$-composition of~$\lin{X}$ is a morphism of $\infty$-algebras, we have
\[
ad \star_0 c'\sigma_d = a\sigma_d = c\sigma_d \star_0 a\rep{d}.
\]
Replacing the $0$-composition by its linear expression in the leftmost and the rightmost terms yields
\[
ad - c'd + c'\sigma_d = c\sigma_d - c\rep{d} + a\rep{d}.
\]
Similarly, using $s_0(\sigma_b)=d$ and $t_0(\sigma_b)=\rep{d}$, and decomposing~$a\sigma_b$, we obtain
\[
ad - c'd + c'\sigma_b = c\sigma_b - c\rep{d} + a\rep{d}.
\]
The difference between the two equations so obtained gives~\eqref{E:ExchangeC}. 
\end{proof}

\begin{theorem}
\label{T:ResolutionContraction}
Let~$X$ be an $\infty$-polygraph with a fixed unital section~$\iota$. Then~$X$ is a polygraphic resolution of~$\cl{X}$ if, and only if, it admits a right $\iota$-contraction.
\end{theorem}

\begin{proof}
Assume that~$X$ is a polygraphic resolution of~$\cl{X}$, and define a right $\iota$-contraction~$\sigma$ of~$X$ thanks to Lemma~\ref{L:FreeRightContraction}. If~$xu$ is an essential monomial of~$\lin{X}$, then~$xu$ and~$\rep{xu}$ have the same image in~$\cl{X}$, so that, by definition of~$\cl{X}$, there exists a $1$-cell~$\sigma_{xu}:xu\fl\rep{xu}$ in~$\lin{X}$. Assume that~$\sigma$ is defined on the $n$-cells of~$\lin{X}$, for~$n\geq 1$, and let~$\alpha u$ be a $\sigma$-essential $n$-monomial of~$\lin{X}$. The $n$-cells defining $s(\sigma_{\alpha u})$ and $t(\sigma_{\alpha u})$ are parallel, so, by hypothesis, there exists an $(n+1)$-cell~$\sigma_{\alpha u}$ with these source and target in~$\lin{X}$.

Conversely, let~$\sigma$ be an $\iota$-contraction of~$X$, and~$a$ and~$b$ be parallel $n$-cells of~$\lin{X}$, for~$n\geq 1$. We have $t(\sigma_a)=\sigma_{s(a)}=\sigma_{s(b)}=t(\sigma_b)$ by hypothesis, so that the $(n+1)$-cell~$\sigma_a\star_n\sigma_b^-$ is well defined, with source $s(\sigma_a)$ and target $s(\sigma_b)$. The fact that $t_k(a)=t_k(b)$ holds for every $0\leq k<n$ implies that
\[
(\sigma_a \star_n \sigma_b)^- \star_{n-1} \sigma_{t_{n-1}(a)}^- \star_{n-2} \cdots \star_0 \sigma_{t_0(a)}^-
\]
is a well-defined $n$-cell of~$\lin{X}$, with source~$a$ and target~$b$, thus proving that~$X_{n+1}$ is acyclic.
\end{proof}

\subsection{Collapsing polygraphic resolutions}
\label{SS:Collapsing}

\subsubsection{Collapsing schemes}

Let~$X$ be an $\infty$-polygraph, and~$Y$ be an indexed subset of~$X$. A \emph{collapsing scheme of~$X$ onto~$Y$} is an injective indexed partial map
\[
X \ofl{\phi} X
\]
of degree~$1$ that satisfies, writing~$\phi_x$ for~$\phi(x)$,
\begin{enumerate}
\item as an indexed set, $X$ admits the partition
\begin{equation}
\label{E:CollapsingPartition}
X \;=\; \im(\phi) \;\amalg\; Y \;\amalg\; \dom(\phi),
\end{equation}
\item for every~$x$ in $\dom(\phi)$, the boundary of~$\phi_x$ satisfies
\begin{equation}
\label{E:CollapsingBoundary}
\dr(\phi_x)=\lambda x + a,
\end{equation}
where~$\lambda$ is a nonzero scalar, and~$a$ is an $n$-cell of~$\lin{X}$ such that $x\notin\cell(a)$,
\item putting $x\trieq_{\phi} y$, for all~$x$ in $\dom(\phi)$ and~$y$ in $\cell(\dr(\phi_x))$, defines a wellfounded order on the cells of~$X$.
\end{enumerate}

Let~$\phi$ be a fixed collapsing scheme of~$X$ onto~$Y$. 
If~$x$ is an $n$-cell of~$\dom(\phi)$, using the same notation as in~\eqref{E:CollapsingBoundary}, define the $n$-cell~$\tilde{x}$ and the $(n+1)$-cell~$\tilde{\phi}_x$ of~$\lin{X}$ by
\[
\tilde{x} = -\frac{1}{\lambda} a
\qquad\text{and}\qquad 
\tilde{\phi}_x = \frac{1}{\lambda} (\phi_x - t(\phi_x)) + \tilde{x},
\]
so that $s(\tilde{\phi}_x)=x$ and $t(\tilde{\phi}_x)=\tilde{x}$ are satisfied, and~$x\tri_{\phi} y$ holds for every $y\in\cell(\tilde{x})$. 

\begin{lemma}
\label{L:Collapsing}
Let~$X$ be an $\infty$-polygraph, $Y$ be an indexed subset of~$X$, and~$\phi$ be a collapsing scheme of~$X$ onto~$Y$. Setting, for every~$n\geq 1$ and every $n$-cell~$x$ of~$Y$, 
\begin{equation}
\label{E:CollapsingSourceTarget}
\cl{s}(x) = \pi(s(x)) 
\qquad\text{and}\qquad
\cl{t}(x) = \pi(t(x)),
\end{equation}
and, for every~$n\geq 0$ and every $n$-cell~$x$ of~$X$, 
\begin{equation}
\label{E:CollapsingProjection}
\pi (x) =
\begin{cases}
x
	&\text{if~$x\in Y$,} \\
\pi(\tilde{x})
	&\text{if~$x\in\dom(\phi)$,} \\
1_{\cl{s}(x)}
	&\text{if~$x\in\im(\phi)$,}
\end{cases}
\end{equation}
defines, by mutual induction, a structure of $\infty$-polygraph on~$Y$ and a morphism $\pi:\lin{X}\fl\lin{Y}$ of $\infty$-algebras.
\end{lemma}

\begin{proof}
By induction on~$n\geq 0$. For~$n=0$, we only have to check that~\eqref{E:CollapsingProjection} defines a map from the $0$-cells of~$X$ to the $0$-cells of~$\lin{Y}$: this holds by induction on the wellfounded order~$\trieq_{\phi}$, using the fact that~$x\tri_{\phi} y$ is satisfied for every~$x$ in~$\dom(\phi)$ and every $y$ in $ \cell(\tilde{x})$. 

Now, fix~$n\geq 1$, and assume that~$Y_{n-1}$ is an $(n-1)$-polygraph, with~$\cl{s}$ and~$\cl{t}$ as source and target maps, and suppose that~\eqref{E:CollapsingProjection} defines a morphism of $(n-1)$-algebras $\pi:\lin{X_{n-1}}\fl\lin{Y_{n-1}}$.

First, define~$\cl{s}$ and~$\cl{t}$ on the~$n$-cells of~$Y$ by~\eqref{E:CollapsingSourceTarget}, and check that the globular relations are satisfied. If~$n>1$ and~$x$ is an $n$-cell of~$Y$, the definition of~$\cl{s}$ and the fact that~$\pi$ commutes with the source map in dimension~$n-1$ give
\[
\cl{s}(\cl{s}(x)) = \cl{s}(\pi(s(x))) = \pi(s(s(x)).
\]
Similarly, we obtain $\cl{s}(\cl{t}(x)) = \pi(s(t(x)))$, so that $\cl{s}\cl{s}=\cl{s}\cl{t}$ is deduced from $ss=st$, and, with the same reasoning, $\cl{t}\cl{s}=\cl{t}\cl{t}$ also holds. 

Next, since~$\trieq_{\phi}$ is wellfounded and $x\tri_{\phi} y$ holds for every~$x$ in $\dom(\phi)$ and every $y$ in $\cell(\tilde{x})$, deduce that~\eqref{E:CollapsingProjection} defines~$\pi$ as a map from the set of $n$-cells of~$X$ to the set of $n$-cells of~$\lin{Y}$. Now, check that~$\pi$ commutes with the source and target maps, to prove that~$\pi$ is a morphism of $n$-algebras from~$\lin{X_n}$ to~$\lin{Y_n}$. Assume that~$x$ is an $n$-cell of~$X$, and that $\cl{s}(\pi(y))=\pi(s(y))$ and $\cl{t}(\pi(y))=\pi(t(y))$ hold for every $n$-cell~$y$ of~$X$ such that $x\tri_{\phi} y$. We distinguish three cases.

If~$x$ is in~$Y$, then the definitions of~$\pi$, $\cl{s}$ and~$\cl{t}$ give
\[
\cl{s}(\pi(x)) =\cl{s}(x) = \pi(s(x))
\qquad\text{and}\qquad
\cl{t}(\pi(x)) = \cl{t}(x) = \pi(t(x)).
\]

If~$x$ is in $\im(\phi)$ (and, thus, $n\geq 1$), then the definition of~$\pi(x)$, the compatibility of~$\cl{s}$ with identities, and the definition of~$\cl{s}$ produce
\[
\cl{s}(\pi(x)) =\cl{s}(1_{\cl{s}(x)}) = \cl{s}(x) = \pi(s(x)).
\]
For the same reasons, we deduce $\cl{t}(\pi(x))=\pi(s(x))$, leaving $\pi(s(x))=\pi(t(x))$ to prove. Let~$y$ be the $(n-1)$-cell of~$X$ such that $\phi_y=x$. The linearity of~$\pi$ on the $(n-1)$-cells, and the definitions of~$\tilde{y}$ and~$\pi(y)$ imply
\[
\pi(s(x)) - \pi(t(x)) = \pi(\dr(x)) = \pi(y) - \pi(\tilde{y}) = 0.
\]
 
If~$x$ is in $\dom(\phi)$, 
\[
\cl{s}(\pi(x)) = \cl{s}(\pi(\tilde{x})) = \pi(s(\tilde{x}))
\qquad\text{and}\qquad
\cl{t}(\pi(x)) = \cl{t}(\pi(\tilde{x})) = \pi(t(\tilde{x}))
\]
come from the definition of~$\pi(x)$ and the induction hypothesis, since every $y\in\cell(\tilde{x})$ satisfies $x\tri_{\phi}y$. Moreover, the globular relations, and $s(\tilde{\phi}_x)=x$ and $t(\tilde{\phi}_x)=\tilde{x}$, imply $s(x)=s(\tilde{x})$ and $t(x)=t(\tilde{x})$.
\end{proof}

\begin{lemma}
\label{L:SectionContraction}
Let~$X$ be an $\infty$-polygraph, $Y$ be an indexed subset of~$X$, and~$\phi$ be a collapsing scheme of~$X$ onto~$Y$. There exist a unique morphism of $\infty$-algebras $\iota:\lin{Y}\fl\lin{X}$, and a unique homotopy~$\eta$ from~$\id_{\lin{X}}$ to the composite morphism of $\infty$-algebras~$\iota\pi$, that satisfy, for every $n$-cell~$x$ of~$Y$,
\begin{equation}
\label{E:DefinitionSection}
\eta_{s_{n-1}(x)} \star_{n-1} \cdots \star_1 \eta_{s_0(x)} \star_0 \iota(x)
\:=\:
x \star_0 \eta_{t_0(x)} \star_1 \cdots \star_{n-1} \eta_{t_{n-1}(x)}
\end{equation}
and, for every $n$-cell~$x$ of~$X$, 
\begin{equation}
\label{E:DefinitionHomotopy}
\eta_x \:=\:
\begin{cases}
1_{x\star_0\eta_{t_0(x)}\star_1\cdots\star_{n-1}\eta_{t_{n-1}(x)}}
	&\text{if~$x\in Y$ or $x\in\im(\phi)$,} \\
\tilde{\phi}_x \star_0 \eta_{t_0(x)} \star_1 \cdots \star_{n-1} \eta_{t_{n-1}(x)} \star_n \eta_{\tilde{x}}
	&\text{if~$x\in\dom(\phi)$.}	
\end{cases}
\end{equation}
Moreover, $\iota$ is a section of the projection~$\pi:\lin{X}\fl\lin{Y}$, and~$\pi(\eta_a)=1_{\pi(a)}$ holds for every $n$-cell~$a$ of~$\lin{X}$.
\end{lemma}

\begin{proof}
Proceed by induction on the dimension.

If~$x$ is a $0$-cell of~$Y$, put $\iota(x)=x$, as required by~\eqref{E:DefinitionSection}, and extend~$\iota$ into a morphism of algebras, that must satisfy $\pi\iota=\id_{\lin{Y}}$ because~$\pi$ is the identity on $Y$.

Then, assume that~$x$ is a $0$-cell of~$X$ and define a $1$-cell~$\eta_x:x\fl\iota\pi(x)$ of~$\lin{X}$, by induction on the wellfounded order~$\trieq_{\phi}$. Fix a $0$-cell~$x$ of~$X$, and assume that~$\eta_y$ is defined by~\eqref{E:DefinitionHomotopy} for every $0$-cell~$y$ of~$X$ such that $x\tri_{\phi} y$, in such a way that $\pi(\eta_y)=1_{\pi(y)}$ holds. First, extend~$\eta$ to every $0$-cell~$a$ of~$\lin{X}$, such that $x\tri_{\phi}y$ holds for every~$y$ in~$\cell(a)$, using Lemma~\ref{L:FreeHomotopy}~\Item{ii}. We distinguish two cases.

If~$x$ is in~$Y$, then $\pi(x)=x$ holds by definition of~$\pi$, hence $\iota\pi(x)=x$ follows by definition of~$\iota$. Put $\eta_x=1_x$ as required by~\eqref{E:DefinitionHomotopy}, so that $\pi(\eta_x)=1_{\pi(x)}$ is satisfied. 

If~$x$ is in~$\dom(\phi)$, then $\pi(x)=\pi(\tilde{x})$ holds by definition of~$\pi$. Define~$\eta_x$ by
\[
\xymatrix @R=0.5em @!C {
x
	\ar@/^/ [rr] ^-{\eta_x} _-{}="s"
	\ar@/_/ [dr] _-{\tilde{\phi}_x} 
&& \pi(x)
\\
& {\tilde{x}}
	\ar@/_/ [ur] _-*+{\eta_{\tilde{x}}}
		\ar@{} "s";[] |(0.45){\sm =}
}
\]
as required by~\eqref{E:DefinitionHomotopy}, with $\pi(\eta_x)=1_{\pi(x)}$ because $\pi(\tilde{\phi}_x)=\pi(\eta_{\tilde{x}})=1_{\pi(x)}$ hold by definition of~$\pi$ and by induction hypothesis on~$\eta_{\tilde{x}}$.

Now, fix~$n\geq 1$, and assume that~$\iota$ and~$\eta$ are defined by~\eqref{E:DefinitionSection} and~\eqref{E:DefinitionHomotopy}, and that they satisfy~$\pi\iota=\id_{\lin{Y}}$ and~$\pi(\eta_a)=1_{\pi(a)}$, up to dimension~$n-1$. 

Let~$x$ be an $n$-cell of~$Y$. First, the compositions involved in~\eqref{E:DefinitionSection} are legal, because the facts that~$\eta$ is a homotopy and that~$\iota\pi$ commutes with each source map~$s_k$, for $0\leq k<n$, imply
\begin{align*}
t_k(\eta_{s_k(x)}) 
	&= \eta_{s_{k-1}(x)} \star_{k-1} \cdots \star_1 \eta_{s_0(x)} \star_0 \iota\pi(s_k(x)) \\
	&= s_k(\eta_{s_{k-1}(x)} \star_{k-1} \cdots \star_1 \eta_{s_0(x)} \star_0 \iota(x)),
\end{align*}
and similarly for the right-hand side of~\eqref{E:DefinitionSection}. Then, the two sides of~\eqref{E:DefinitionSection} are parallel, because
\begin{align*}
s(\eta_{s_{n-1}(x)} \star_{n-1} \cdots \star_1 \eta_{s_0(x)} \star_0 \iota(x))
	&= s(\eta_{s(x)}) \\
	&= s(x \star_0 \eta_{t_0(x)} \star_1\cdots\eta_{t_{n-1}(x)})
\end{align*}
hold, a similar computation giving the same equality between the targets. Thus, $\iota(x)$ can be defined by~\eqref{E:DefinitionSection}, because every other cell involved in the left-hand side is invertible. Moreover, $\pi\iota(x)=x$ holds by induction hypothesis on~$\eta$.

Next, consider an $n$-cell~$x$ of~$X$, and define~$\eta_x$ by induction on the wellfounded order~$\trieq_{\phi}$. Assume that~$\eta_y$ is defined for every $n$-cell~$y$ of $X$ such that $x\tri_{\phi} y$, and extend~$\eta$, thanks to Lemma~\ref{L:FreeHomotopy}~\Item{ii}, to every $n$-cell~$a$ of~$\lin{X}$ such that $x\tri_{\phi} y$ holds for every $y\in\cell(a)$. Define~$\eta_x$ as in~\eqref{E:DefinitionHomotopy}, and proceed by case analysis to check that it satisfies~\eqref{E:SourceHomotopy}, \eqref{E:TargetHomotopy}, and $\pi(\eta_x)=1_{\pi(x)}$.

If~$x$ is in~$Y$, then
\[
\pi(x) = x
\qquad\text{and}\qquad
\eta_x = 1_{x\star_0\eta_{t_0(x)}\star_1\cdots\star_{n-1}\eta_{t_{n-1}(x)}}
\]
hold. So~\eqref{E:SourceHomotopy} is satisfied, and~\eqref{E:TargetHomotopy} is equivalent to~\eqref{E:DefinitionSection}. Moreover, $\pi(\eta_x)=1_{\pi(x)}$ holds by induction hypothesis on~$\eta$.

If~$x$ is in~$\dom(\phi)$, then
\[
\pi(x)=\pi(\tilde{x})
\qquad\text{and}\qquad
\eta_x = \tilde{\phi}_x\star_0\eta_{t_0(x)}\star_1\cdots\star_{n-1}\eta_{t_{n-1}(x)} \star_n \eta_{\tilde{x}}
\]
are satisfied. Note that compositions involved in the definition of~$\eta_x$ are valid, because~$x$ and~$\tilde{x}$ are parallel. Now, observe that
\[
s(\eta_x) = s(\tilde{\phi}_x)\star_0\eta_{t_0(x)}\star_1\cdots\star_{n-1}\eta_{t_{n-1}(x)}
\qquad\text{and}\qquad
t(\eta_x) = t(\eta_{\tilde{x}})
\]
hold. By definition, $s(\tilde{\phi}_x)=x$, which implies~\eqref{E:SourceHomotopy}. Then~\eqref{E:TargetHomotopy} follows from the induction hypothesis applied to~$\tilde{x}$, which satisfies $x\tri_{\phi} y$ for every $y\in\cell(\tilde{x})$, and from the fact that~$\pi(x)=\pi(\tilde{x})$ holds. Finally, the definition of~$\pi$ and the induction hypothesis on~$\eta$ imply $\pi(\eta_x)=1_{\pi(x)}$.

If~$x$ is in~$\im(\phi)$, we have
\[
\pi(x)=1_{\cl{s}(x)}
\qquad\text{and}\qquad
\eta_x = 1_{x\star_0\eta_{t_0(x)}\star_1\cdots\star_{n-1}\eta_{t_{n-1}(x)}}.
\]
So~\eqref{E:SourceHomotopy} is satisfied, and
\[
\eta_{s_{n-1}(x)}\star_{n-1}\cdots\star_1\eta_{s_0(x)}\star_0\iota\pi(x) 
	= \eta_{s(x)}
\]
holds. Thus, \eqref{E:TargetHomotopy} is equivalent to 
\[
\eta_{s(x)} = x\star_0 \eta_{t_0(x)} \star_1\cdots\star_{n-1} \eta_{t_{n-1}(x)}.
\]
Recall that the right-hand side of the last equality is equal to $x - t(x) + \eta_{t(x)}$, and let~$y$ be the $(n-1)$-cell of~$\dom(\phi)$ such that~$\phi_y=x$. By induction hypothesis on~$\eta$, and using $a\star_k b = a - t_k(a) + b$, we obtain
\[
\eta_y 
	= \tilde{\phi}_y\star_0\eta_{t_0(y)}\star_1\cdots\star_{n-1}\eta_{t_{n-1}(y)}\star_n\eta_{\tilde{y}} 
	= \tilde{\phi}_y - \tilde{y} + \eta_{\tilde{y}}.
\]
Now, $\tilde{\phi}_y=\frac{1}{\lambda}(x-t(x)) + \tilde{y}$, satisfied by definition, implies
\[
\eta_y = \frac{1}{\lambda} (x - t(x)) + \eta_{\tilde{y}}.
\]
So, by definition of~$y$ and by linearity of~$\eta$, we obtain
\[
\eta_{s(x)} - \eta_{t(x)} 
	= \lambda(\eta_y - \eta_{\tilde{y}})
	= x - t(x),
\]
which, in turn, implies~\eqref{E:TargetHomotopy}. Finally, $\pi(\eta_x)=1_{\pi(x)}$ is satisfied by induction hypothesis on~$\eta$.
\end{proof}

\begin{theorem}
\label{T:Collapsing}
Let~$A$ be an algebra and~$X$ be a polygraphic resolution of~$A$. If~$Y$ is an indexed subset of~$X$, and~$\phi$ is a collapsing scheme of~$X$ onto~$Y$, then~$Y$, equipped with the structure of $\infty$-polygraph of Lemma~\ref{L:Collapsing}, is a polygraphic resolution of~$A$.
\end{theorem}

\begin{proof}
First, prove that the algebras~$\cl{X}$ and~$\cl{Y}$, respectively presented by~$X$ and~$Y$, are isomorphic. If~$a$ and~$b$ are $0$-cells of~$\lin{X}$ (resp.\ of~$\lin{Y}$) that are identified in~$\cl{X}$ (resp.\ in~$\cl{Y}$), there exists a $1$-cell $f:a\fl b$ in~$\lin{X}$ (resp.\ in~$\lin{Y}$); because~$\pi$ (resp.~$\iota$) is a morphism of $\infty$-algebras, $\pi(f):\pi(a)\fl\pi(b)$ (resp.\ $\iota(f):\iota(a)\fl\iota(b)$) is a $1$-cell of~$\lin{Y}$ (resp.\ of~$\lin{X}$), meaning that~$\pi(a)$ and~$\pi(b)$ are identified in~$\cl{Y}$ (resp.\ that~$\iota(a)$ and~$\iota(b)$ are identified in~$\cl{X}$). As a consequence, the morphisms of $\infty$-algebras $\pi:\lin{X}\fl\lin{Y}$ and $\iota:\lin{Y}\fl\lin{X}$ induce morphisms of algebras $\cl{\pi}:\cl{X}\fl\cl{Y}$ and $\cl{\iota}:\cl{Y}\fl\cl{X}$. Moreover, the equality $\pi\iota=\id_{\lin{Y}}$ induces $\cl{\pi}\cl{\iota}=\id_{\cl{Y}}$. Finally, for every $0$-cell~$a$ of~$\lin{X}$, the $1$-cell $\eta_a:a\fl\iota\pi(a)$ of~$\lin{X}$ proves that $\cl{\iota}\cl{\pi}(a)=a$.

Now, prove that, for every~$n\geq 1$ and every $n$-sphere~$(a,b)$ of~$\lin{Y}$, there exists an $(n+1)$-cell of source~$a$ and target~$b$ in~$\lin{Y}$. Apply the morphism of $\infty$-algebras~$\iota$ from Lemma~\ref{L:SectionContraction} to both~$a$ and~$b$, to obtain an $n$-sphere $(\iota(a),\iota(b))$ of~$\lin{X}$. Since~$X$ is a polygraphic resolution of~$A$, there exists an $(n+1)$-cell $f:\iota(a)\fl\iota(b)$ in~$\lin{X}$. Apply the projection~$\pi$, and use $\pi\iota=\id_{\lin{Y}}$, to obtain an $(n+1)$-cell $\pi(f):a\fl b$ in~$\lin{Y}$. 
\end{proof}

\section{Squier's polygraphic resolution of associative algebras}
\label{S:SquierResolution}

In this section, given a convergent presentation of an algebra, the coherent presentation given by Squier's theorem, Theorem~\ref{T:Squier}, is extended into a polygraphic resolution. The first step is to construct a polygraphic analogue (but with cubical cells instead of simplicial ones) of the standard resolution of an algebra, which is proved to be a resolution in Theorem~\ref{T:StandardPolygraphicResolution} by building a contraction. Using Theorem~\ref{T:Collapsing}, this explicit but very large resolution is then contracted into Squier's polygraphic resolution, with one $n$-cell for each critical $n$-branching of the initial convergent presentation. Several examples are given at the end of the section.

\subsection{The standard polygraphic resolution of an associative algebra}
\label{SS:StandardPolygraphicResolution}

\subsubsection{Cubical faces}

Let $A=\K\oplus A_+$ be an augmented algebra. For~$n\geq 1$, define the vector space
\[
A^{(n)} = 
\bigoplus_{
	\substack{
		1\leq k\leq n \\ 
		i_1+\cdots+i_k=n
		}
	}
	A^{\tens i_1+1}_+ \tens \cdots \tens A^{\tens i_k+1}_+ .
\]
Vertical bars are used to denote the innermost products of copies of~$A_+$, so that~$A^{(n)}$ is made of the linear combinations of elements~$a^1_0\vert\cdots\vert a^1_{i_1}\tens\cdots\tens a^k_0\vert\cdots\vert a^k_{i_k}$ with~$n$ vertical bars. In particular, each $a_0\vert\cdots\vert a_n$ of~$A^{\tens n+1}_+$ belongs to~$A^{(n)}$. If~$\Br$ is a linear basis of~$A_+$, the elements $u^1_0\vert\cdots\vert u^1_{i_1}\tens\cdots\tens u^k_0\vert\cdots\vert u^k_{i_k}$, with each~$u_i^j$ in~$\Br$, form a basis of the vector space~$A^{(n)}$. 

For $1\leq i\leq n$, let~$d_i^-$ (resp.~$d_i^+$) be the linear map from~$A^{(n)}$ to~$A^{(n-1)}$ that replaces the~$i^{\text{th}}$ vertical bar, counting from the left, with a tensor (resp.\ the product of~$A_+$). For example, 
\[
d_2^-(a\tens b\vert c\tens d\vert e) = a\tens b\vert c\tens d\tens e
\qquad\text{and}\qquad
d_1^+(a\tens b\vert c\tens d\vert e) = a\tens bc \tens d\vert e.
\]
It follows from the definition that these maps satisfy the cubical relations:
\begin{equation}
\label{E:CubicalRelations}
d_j^{\beta} d_i^{\alpha} = d_i^{\alpha} d_{j+1}^{\beta} 
\qquad\text{if $i\leq j$.}
\end{equation}
Define~$d_i^{(j)}$ as~$d_i^-$ if~$j$ is odd, and as~$d_i^+$ if~$j$ is even. For~$\alpha\in\ens{-,+}$,  denote by~$-\alpha$ the opposite of~$\alpha$.

\begin{theorem}
\label{T:StandardPolygraphicResolution}
Assume that $A=\K\oplus A_+$ is an augmented algebra, and that~$\Br$ is a linear basis of~$A_+$. Setting $\Std(\Br)_n=\Br^{n+1}$ and, for every $u_0\vert\cdots\vert u_n$ in~$\Br^{n+1}$, 
\begin{align*}
s(u_0\vert\cdots\vert u_n) &= 
	\sum_{
		\substack{
			1\leq k\leq n \\ 
			1\leq i_1<\cdots<i_k\leq n
		}
	}
	(-1)^{k+1} d_{i_1}^{-(n-i_1)} \cdots d_{i_k}^{-(n-i_k)} (u_0\vert\cdots\vert u_n),
\\[0.5em]
t(u_0\vert\cdots\vert u_n) &= 
	\sum_{
		\substack{
			1\leq k\leq n \\ 
			1\leq i_1<\cdots<i_k\leq n
		}
	}
	(-1)^{k+1} d_{i_1}^{(n-i_1)}\cdots d_{i_k}^{(n-i_k)} (u_0\vert\cdots\vert u_n),
\end{align*}
defines a polygraphic resolution of~$A$.
\end{theorem}

Before giving the proof of Theorem~\ref{T:StandardPolygraphicResolution}, let us state a consequence for augmented algebras given with a convergent presentation.

\begin{corollary}
\label{C:ConvergentStandardPolygraphicResolution}
Assume that~$A$ is an augmented algebra, and that~$X$ is a convergent left-monomial presentation of~$A$. Then the algebra~$A$ admits a polygraphic resolution $\Std(\red(X))$, whose $n$-cells are the $(n+1)$-uples $u_0\vert\cdots\vert u_n$ of nontrivial reduced monomials of~$\lin{X}$, with source and target 
\begin{align*}
s(u_0\vert\cdots\vert u_n) 
	&\quad= \sum_{\substack{0\leq i\leq n+1 \\ n+1-i \text{ even}}}
	d_i(u_0\vert\cdots\vert u_n) \quad+\quad s'(u_0\vert\cdots\vert u_n) \\
\text{and}\qquad
t(u_0\vert\cdots\vert u_n) 
	&\quad= \sum_{\substack{0\leq i\leq n+1 \\ n+1-i \text{ odd}}}
	d_i(u_0\vert\cdots\vert u_n) \quad+\quad t'(u_0\vert\cdots\vert u_n),
\end{align*}
where $s'(u_0\vert\cdots\vert u_n)$ and $t'(u_0\vert\cdots\vert u_n)$ are identities, and where the maps~$d_i$ are defined by
\[
d_i(u_0\vert\cdots\vert u_n) =
\begin{cases}
u_0\tens u_1\vert\cdots\vert u_n 
	&\text{if $i=0$,} 
\\
u_0\vert\cdots\vert \rep{u_{i-1}u_i} \vert\cdots\vert u_n
	&\text{if $1\leq i\leq n$,}
\\
u_0\vert\cdots\vert u_{n-1}\tens u_n
	&\text{if $i=n+1$.}
\end{cases}
\]
\end{corollary}

\begin{proof}
The algebra~$A$ being augmented, for reduced monomials~$u$ and~$v$ in~$\lin{X}$, the normal form~$\rep{uv}$ belongs to the vector space spanned by nontrivial reduced monomials of~$\lin{X}$.
By Theorem~\ref{T:StandardBasis}, the reduced monomials of~$\lin{X}$ form a linear basis $\red(X)$ of~$A$, and the product in~$A$ of two reduced monomials~$u$ and~$v$ is given by~$u\cdot v=\rep{uv}$. Apply Theorem~\ref{T:StandardPolygraphicResolution} to $\red(X)$ to obtain $\Std(\red(X))$. Then, put $d_0=d_0^-$, $d_i=d_i^+$ for $1\leq i\leq n$, and $d_{n+1}=d_n^-$, and define $s'(u_0\vert\cdots\vert u_n)$ and $t'(u_0\vert\cdots\vert u_n)$ as the rest of the terms that occur in the source and target maps of $\Std(\red(X))$. Finally, observe that~$s'(u_0\vert\cdots\vert u_n)$ and $t'(u_0\vert\cdots\vert u_n)$ are identities because they are linear combinations of identities: each $d_i^-(u_0\vert\cdots\vert u_n)$, for $1\leq i<n$, and composites of two or more face maps.
\end{proof}

Fix an augmented algebra $A=\K\oplus A_+$ and a linear basis~$\Br$ of~$A_+$ for the rest of the section, and define~$\Std(\Br)$ as in Theorem~\ref{T:StandardPolygraphicResolution}.

\begin{proof}[\bfseries Proof of Theorem~\ref{T:StandardPolygraphicResolution}]
Lemma~\ref{L:StandardResolutionDefined} states that the definition of $\Std(\Br)$ makes sense. Lemma~\ref{L:StandardResolutionPresented} shows that the $1$-polygraph underlying $\Std(\Br)$ is a presentation of the algebra~$A$. Finally, Proposition~\ref{P:StandardResolutionContraction} exhibits a right contraction of $\Std(\Br)$, so that Lemma~\ref{L:FreeRightContraction} and Theorem~\ref{T:ResolutionContraction}  conclude the proof of Theorem~\ref{T:StandardPolygraphicResolution}.
\end{proof}

\begin{lemma}
\label{L:StandardResolutionDefined}
 The $\infty$-polygraph $\Std(\Br)$ is well defined.
 \end{lemma}
 
\begin{proof} 
By definition, the source and target maps of $\Std(\Br)$ are linear. Now, check that they satisfy the globular relations. Let~$x$ be a fixed $n$-cell of $\Std(\Br)$. If~$k>1$, then $d_{i_1}^{-(n-i_1)} \cdots d_{i_k}^{-(n-i_k)} (x)$ has dimension strictly lower than~$n-2$, so that~$s_{n-2}$ is the identity on such a cell. It follows that
\[
s_{n-2}s_{n-1} (x) 
\quad = \quad 
\sum_{1\leq i\leq n} s_{n-2} (d_i^{-(n-i)} (x))
\quad + \sum_{
		\substack{
			1<k\leq n \\ 
			1\leq i_1<\cdots<i_k\leq n
		}
	}
	(-1)^{k+1} d_{i_1}^{-(n-i_1)} \cdots d_{i_k}^{-(n-i_k)} (x) .
\]
Fix~$i$ in $\ens{1,\dots,n}$. By definition, after noticing that~$d_j^{-(n-1-j)}=d_j^{(n-j)}$, we obtain
\[
s_{n-2} (d_i^{-(n-i)} (x))
\quad = \quad 
\sum_{
 			\substack{
 				1\leq k<n \\
 				1\leq j_1<\cdots<j_k<n
 			}
 		}
 		(-1)^{k+1} d_{j_1}^{(n-j_1)} \cdots d_{j_k}^{(n-j_k)} d_i^{-(n-i)} (x).
\]
Now, fix~$1\leq k<n$ and $1\leq  j_1<\cdots<j_k<n$. Let~$l$ be the smallest element of $\ens{1,\dots,k}$ such that~$i\leq j_l$, or~$k+1$ if $i>j_k$. For each $1\leq p\leq k+1$, define~$i_p$ as~$j_p$ if~$p<l$, as~$i$ if~$p=l$, and as~$j_p+1$ if~$p>l$. The cubical relations~\eqref{E:CubicalRelations} imply
\[
d_{j_1}^{(n-j_1)} \cdots d_{j_k}^{(n-j_k)} d_i^{-(n-i)} = 
	d_{i_1}^{(n-i_1)} \cdots d_{i_{l-1}}^{(n-i_{l-1})} 
	d_{i_l}^{-(n-i_l)} 
	d_{i_{l+1}}^{-(n-i_{l+1})} \cdots d_{i_{k+1}}^{-(n-i_{k+1})}.
\]
As a consequence
\[
s_{n-2} (d_i^{-(n-i)} (x))
\quad = \quad 
\sum_{
	\substack{
		1\leq l\leq k\leq n \\
		1\leq i_1<\cdots<i_k\leq n
		}
	}
	(-1)^k d_{i_1}^{(n-i_1)}\cdots d_{i_{l-1}}^{(n-i_{l-1})} d_{i_l}^{-(n-i_l)} \cdots d_{i_{k}}^{-(n-i_{k})} (x).
\]
The terms corresponding to~$l=1$ cancel out with the remaining part of $s_{n-2}s_{n-1}(x)$, leaving, if we sum over $1\leq l<k\leq n$ instead of $1<l\leq k\leq n$,
\[
s_{n-2}s_{n-1} (x) = \sum_{
	\substack{
		1\leq l<k\leq n \\
		1\leq i_1<\cdots<i_k\leq n
		}
	}
	(-1)^k d_{i_1}^{(n-i_1)}\cdots d_{i_l}^{(n-i_l)} d_{i_l+1}^{-(n-i_{l+1})} \cdots d_{i_k}^{-(n-i_k)} (x).
\]
Similar computations lead to the same expression for $s_{n-2}t_{n-1}(x)$, so that $ss=st$. A similar argument proves that both $ts(x)$ and $tt(x)$ are equal to
\[
\sum_{
	\substack{
		1\leq l<k\leq n \\
		1\leq i_1<\cdots<i_k\leq n
		}
	}
	(-1)^k 
	d_{i_1} ^{-(n-i_1)} \cdots d_{i_l} ^{-(n-i_l)} 
	d_{i_{l+1}}^{(n-i_{l+1})} \cdots d_{i_k} ^{(n-i_k)} (x).
\qedhere
\]
\end{proof}

\begin{lemma}
\label{L:StandardResolutionPresented}
The underlying $2$-polygraph of $\Std(\Br)$ is the standard coherent presentation $\Std(\Br)_2$ of~$A$, as given in Example~\ref{X:StandardCoherentPresentation}.
\end{lemma}

\begin{proof}
By definition, $\Std(\Br)$ has the same sets of $0$-cells, $1$-cells and $2$-cells as the standard coherent presentation $\Std(\Br)_2$ of Example~\ref{X:StandardCoherentPresentation}. Moreover, the source and target of $1$-cells coincide:
\[
s(u\vert v) = d_1^-(u\vert v) = u\tens v
\qquad\text{and}\qquad
t(u\vert v) = d_1^+(u\vert v) = uv.
\]
The relations between the composition~$\star_0$ and the linear structure in a $\infty$-vector space imply that the source and the target of $2$-cells also agree, because
\begin{align*}
s(u\vert v\vert w) 
	&= d_1^+ (u\vert v\vert w) + d_2^-(u\vert v\vert w) - d_1^+d_2^-(u\vert v\vert w) \\
	&= uv\vert w + u\vert v\tens w - uv\tens w \\
	&= u\vert v\tens w \star_0 uv\vert w 
\end{align*}
and
\begin{align*}
t(u\vert v\vert w) 
	&= d_1^-(u\vert v \vert w) + d_2^+(u\vert v\vert w) - d_1^-d_2^+(u\vert v\vert w) \\
	&= u\tens v \vert w + u\vert vw - u\tens vw \\
	&= u\tens v \vert w \star_0 u\vert vw.
\qedhere
\end{align*}
\end{proof}

\begin{proposition}
\label{P:StandardResolutionContraction}
Let~$\iota$ be the unital section of~$\Std(\Br)$ given by the inclusion of~$\Br$ into $\lin{\Std(\Br)}$. The assignment $\sigma_{u_0\vert\cdots\vert u_n\tens u_{n+1}}=u_0\vert\cdots\vert u_{n+1}$ defines a right $\iota$-contraction of~$\Std(\Br)$.
\end{proposition}

The proof of Proposition~\ref{P:StandardResolutionContraction} uses the following lemma.

\begin{lemma}
\label{L:StandardResolutionContractionLemma}
With the same hypotheses as in Proposition~\ref{P:StandardResolutionContraction}, fix~$n\geq 1$, and assume that the assignment $\sigma_{u_0\vert\cdots\vert u_k\tens u_{k+1}}=u_0\vert\cdots\vert u_{k+1}$ defines a right $\iota$-contraction of~$\Std(\Br)$ up to dimension~$n-1$. Then, the relations
\begin{align}
\sigma_{Dd_{n+1}^-(x)} - \sigma_{Dd_{n}^+d_{n+1}^-(x)}
	&\:= \: D(x) - Dd_{n}^+(x)
\label{E:DifferencePlus}
\\
\text{and} \qquad
\sigma_{Dd_{n+1}^-(x)} - \sigma_{Dd_{n}^-d_{n+1}^-(x)}
	&\:= \: D(x) - Dd_{n}^-(x) - Dd_{n+1}^+(x) + Dd_{n}^-d_{n+1}^+(x)
\label{E:DifferenceMinus}
\end{align}
are satisfied for every $D=d_{i_1}^{\alpha_1}\cdots d_{i_k}^{\alpha_k}$, with $1\leq i_1<\cdots<i_k<n$ and~$\alpha_1$, \dots, $\alpha_k$ in $\ens{-,+}$, and every $x=u_0\vert\cdots\vert u_{n+1}$. 
\end{lemma}

\begin{proof}
Distinguish two cases. If every~$\alpha_i$ is a~$+$, then $Dd_{n+1}^-(x)$ has shape $v_0\vert\cdots\vert v_p\vert u_n \tens u_{n+1}$, so that
\begin{align*}
\sigma_{Dd_{n+1}^-(x)} 
	&\:=\: \sigma_{v_0\vert\cdots\vert v_p\vert u_n \tens u_{n+1}} 
	\:=\: v_0\vert\cdots\vert v_p\vert u_n \vert u_{n+1}
	\:=\: D(x)
\\
\text{and} \quad
\sigma_{Dd_{n}^+d_{n+1}^-(x)}
	&\:=\: \sigma_{v_0\vert\cdots\vert v_p u_n \tens u_{n+1}}
	\:=\: v_0\vert\cdots\vert v_p u_n \vert u_{n+1}
	\:=\: Dd_{n}^+(x),
\end{align*}
and, because~$\sigma$ is a right $\iota$-contraction up to dimension~$n-1$, and since $\iota\pi(u_n\tens u_{n+1})=u_n u_{n+1}$, 
\begin{align*}
\sigma_{Dd_{n}^-d_{n+1}^-(x)}
	&\:=\: \sigma_{v_0\vert\cdots\vert v_p\tens u_n \tens u_{n+1}} \\
	&\:=\: v_0\vert\cdots\vert v_p\tens u_n\vert u_{n+1}
		\:- \: v_0\vert\cdots\vert v_p \tens u_n u_{n+1}
		\:+\: v_0\vert\cdots\vert v_p \vert u_n u_{n+1} \\
	&\:=\: Dd_{n}^-(x) - Dd_{n}^-d_{n+1}^+(x) + Dd_{n+1}^+(x).
\end{align*}

Otherwise, write $D=D_1 d_i^- D_2$, where the indices occurring in~$D_1$ are strictly smaller than~$i$, the ones that occur in~$D_2$ are strictly greater than~$i$, and all exponents that occur in~$D_2$ are~$+$s. Define $y=D_1(u_0\vert\cdots\vert u_i)$ and $v_0\vert\cdots\vert v_p=D_2(u_{i+1}\vert\cdots\vert u_{n-1})$, and write $z=y \tens v_0\cdots v_p u_n u_{n+1} - \sigma_{y \tens v_0\cdots v_p u_n u_{n+1}}$. The fact that~$\sigma$ is a right $\iota$-contraction implies, on the one hand,
\begin{align*}
\sigma_{Dd_{n+1}^-(x)} 
	&\:=\: \sigma_{y\tens v_0\vert\cdots\vert v_p \vert u_n \tens u_{n+1}} 
	\:=\: y \tens v_0\vert\cdots\vert v_p\vert u_n \vert u_{n+1} - z 
	\:=\: D(x) - z
\\
\text{and}
\quad \sigma_{Dd_{n}^+d_{n+1}^-(x)} 
	&\:=\: \sigma_{y\tens v_0\vert\cdots\vert v_p u_n \tens u_{n+1}} 
	\:=\: y \tens v_0\vert\cdots\vert v_p u_n \vert u_{n+1} - z
	\:=\: Dd_{n}^+(x) - z,
\end{align*}
and, on the other hand,
\begin{align*}
\sigma_{Dd_{n}^-d_{n+1}^-(x)}
	\:=\:& \sigma_{y\tens v_0\vert\cdots\vert v_p \tens u_n \tens u_{n+1}} \\
	\:=\:& y\tens v_0\vert\cdots\vert v_p \tens u_n\vert u_{n+1} 
		\:-\: y\tens v_0\vert\cdots\vert v_p \tens u_n u_{n+1} \\
		&\:+\: y\tens v_0\vert\cdots\vert v_p \vert u_n u_{n+1} 
		\:-\: z \\
	\:=\:& Dd_{n}^-(x) - Dd_{n}^-d_{n+1}^+(x) + Dd_{n+1}^+(x) - z. 
\qedhere
\end{align*}
\end{proof}

\begin{proof}[{\bfseries Proof of Proposition~\ref{P:StandardResolutionContraction}}]
By induction on the dimension. By definition, $\iota\pi(u_0\tens\cdots\tens u_n)=u_0\cdots u_n$, so that the $\iota$-reduced $0$-cells of~$\lin{\Std(\Br)}$ are the $0$-cells of~$\Std(\Br)$; as a consequence, the $\iota$-essential monomials of~$\lin{\Std(\Br)}$ are the monomials $u\tens v$, for~$u$ and~$v$ in~$\Br$. Define~$\sigma$ on $0$-cells as the unique right $\iota$-contraction that satisfies $\sigma_{u\tens v} = u\vert v$, which is legitimate since the equalities $s(\sigma_{u\tens v})=u\tens v$ and~$t(\sigma_{u\tens v})=\rep{u\tens v}=uv$ hold.

Fix~$n\geq 1$, and assume that the assignment $\sigma_{u_0\vert\cdots\vert u_k\tens u_{k+1}}=u_0\vert\cdots\vert u_{k+1}$ defines a right $\iota$-contraction of~$\Std(\Br)$ for every $k<n$. 
Fix an $n$-cell $u_0\vert\cdots\vert u_n$ of~$\Std(\Br)$ and an element~$u_{n+1}$ of~$\Br$. Note that the $n$-monomial $u_0\vert\cdots\vert u_n\tens u_{n+1}$ is not in the image of~$\sigma$, so that it is $\sigma$-reduced, and put $\sigma_{u_0\vert\cdots\vert u_n\tens u_{n+1}}=u_1\vert\cdots\vert u_{n+1}$. To ensure that this defines a unique right $\iota$-contraction on the $n$-cells of~$\Std(\Br)$, we must check 
\begin{align}
s(u_0\vert\cdots\vert u_{n+1}) 
	&\:=\: 
	u_0\vert\cdots\vert u_n \tens u_{n+1} 
	\:-\: t(u_0\vert\cdots\vert u_n)\tens u_{n+1} 
	\:+\: \sigma_{t(u_0\vert\cdots\vert u_n)\tens u_{n+1}}
\label{E:SourceContraction}
\\
\text{and} \quad
t(u_0\vert\cdots\vert u_{n+1}) 
	&\:=\: \sigma_{s(u_0\vert\cdots\vert u_n)\tens u_{n+1}} \,.
\label{E:TargetContraction}
\end{align}
By definition of $t(u_0\vert\cdots\vert u_n)$ and linearity of~$\sigma$,
\begin{align*}
\sigma_{t(u_0\vert\cdots\vert u_n)\tens u_{n+1}}
	&= \sum_{
		\substack{
			1\leq k\leq n \\ 
			1\leq i_1<\cdots<i_k\leq n
		}
	}
	(-1)^{k+1} \sigma_{d_{i_1}^{(n-i_1)} \cdots d_{i_k}^{(n-i_k)} d_{n+1}^- (x)} 
\\
	&= \sum_{
		\substack{
			1\leq k\leq n \\ 
			1\leq i_1<\cdots<i_k\leq n
			}
		}
	(-1)^{k+1}
	\left(
		\sigma_{d_{i_1}^{(n-i_1)} \cdots d_{i_k}^{(n-i_k)} d_{n+1}^- (x)}
		- \sigma_{d_{i_1}^{(n-i_1)} \cdots d_{i_k}^{(n-i_k)} d_n^+ d_{n+1}^- (x)}
	\right).
\end{align*}
We use~\eqref{E:DifferencePlus} for each $D=d_{i_1}^{(n-i_1)} \cdots d_{i_k}^{(n-i_k)}$, and we note that $d_{n}^+(x) = d_{n}^{-(n+1-n)}$, to obtain
\[
\sigma_{t(u_0\vert\cdots\vert u_n)\tens u_{n+1}}
	= \sum_{
		\substack{
			1\leq k\leq n \\ 
			1\leq i_1<\cdots<i_k\leq n
		}
	}
	(-1)^{k+1} d_{i_1}^{-(n+1-i_1)} \cdots d_{i_k}^{-(n+1-i_k)} (x).
\]
Since
\[
t(u_0\vert\cdots\vert u_n) \tens u_{n+1}
	= \sum_{
		\substack{
			1\leq k\leq n \\ 
			1\leq i_1<\cdots<i_k\leq n
		}
	}
	(-1)^{k+1} d_{i_1}^{-(n+1-i_1)} \cdots d_{i_k}^{-(n+1-i_k)} d_n^{-(n+1-n)} (x)
\]
and $u_0\vert\cdots\vert u_n \tens u_{n+1} = d_{n+1}^{(n+1)-n} (x)$, we conclude that~\eqref{E:SourceContraction} holds. For~\eqref{E:TargetContraction}, we start with
\begin{align*}
\sigma_{s(u_0\vert\cdots\vert u_n)\tens u_{n+1}}
	&= \sum_{
		\substack{
			1\leq k\leq n \\ 
			1\leq i_1<\cdots<i_k\leq n
		}
	}
	(-1)^{k+1} \sigma_{d_{i_1}^{-(n-i_1)} \cdots d_{i_k}^{-(n-i_k)} d_{n+1}^- (x)} 
\\
	&= \sum_{
		\substack{
			1\leq k<n \\ 
			1\leq i_1<\cdots<i_k<n
			}
		}
	(-1)^{k+1}
	\left(
		\sigma_{d_{i_1}^{-(n-i_1)} \cdots d_{i_k}^{-(n-i_k)} d_{n+1}^- (x)}
		- \sigma_{d_{i_1}^{-(n-i_1)} \cdots d_{i_k}^{-(n-i_k)} d_{n}^- d_{n+1}^- (x)}
	\right).
\end{align*}
We use~\eqref{E:DifferencePlus} for each $D=d_{i_1}^{-(n-i_1)} \cdots d_{i_k}^{-(n-i_k)}$. We note that~$D$ appears with the same sign as $Dd_{n}^-d_{n+1}^+$, which is opposite to the one of $Dd_{n}^-$ and $Dd_{n+1}^+$. Moreover, the lengths of~$D$ and $Dd_{n}^-d_{n+1}^+$ have the same parity, which is opposite to the parity of the lengths of $Dd_{n}^-$ and $Dd_{n+1}^+$. Hence, after observing that $d_{n}^-=d_{n}^{(n+1-n)}$ and $d_{n+1}^+=d_{n+1}^{(n+1-(n+1))}$, we deduce
\[
\sigma_{s(u_0\vert\cdots\vert u_n)\tens u_{n+1}}
	= \sum_{
		\substack{
			1\leq k\leq n+1 \\ 
			1\leq i_1<\cdots<i_k\leq n+1
			}
		}
	(-1)^{k+1} d_{i_1}^{(n+1-i_1)} \cdots d_{i_k}^{(n+1-i_k)}(x),
\]
which is the definition of~$t(x)$, thus concluding the proof.
\end{proof}

\subsection{Squier's polygraphic resolution}
\label{SS:SquierResolution}

\subsubsection{The rightmost rewriting step} 

Assume that~$X$ is a reduced convergent left-monomial $1$-polygraph, and let~$u$ be a monomial of~$\lin{X}$. Because~$X$ is reduced, there is at most one rewriting step of source~$u$ in~$\lin{X}$ for each submonomial of~$u$. As a consequence, the set of rewriting steps of~$\lin{X}$ of source~$u$ is finite. Moreover, this set, extended with the identity~$1_u$, is totally ordered by~$f\sqsubseteq g$ if~$f=1_u$, or if, writing $f=v\alpha v'$ and $g=w\beta w'$, the length of~$v$ is smaller than the length of~$w$. The maximal element for this order is denoted by~$\rho_u$. 
By definition, we have $\rho_u=1_u$ if, and only if, $u$ is reduced. Furthermore, if~$u$ and~$v$ are monomials of~$\lin{X}$, then the inequalities $\rho_u v \sqsubseteq \rho_{uv_1} v_2 \sqsubseteq \rho_{uv}$ hold for every factorisation $v=v_1v_2$.

\subsubsection{Critical, subcritical and supercritical cells}

Fix an augmented algebra~$A$, and a reduced convergent left-monomial presentation~$X$ of~$A$.
Let $u_0\vert\cdots\vert u_n$ be an $n$-cell of $\Std(\red(X))$. For~$i$ in $\ens{0,\dots,n}$, define $\delta_i(u_0\vert\cdots\vert u_n)$ as the pair $(v_i,w_i)$ of (possibly trivial) monomials of~$\lin{X}$ such that $u_i=v_i w_i$ and
\begin{enumerate}
\item if~$i=0$, then~$v_0$ is of length~$1$,
\item if~$i\geq 1$, then~$v_i$ is the shortest left-factor of~$u_i$ that satisfies $\rho_{u_0\cdots u_{i-1}v_i} w_i = \rho_{u_0\cdots u_i}$.
\end{enumerate}
Define $i$-critical, $i$-subcritical and $i$-supercritical $n$-cells of $\Std(\red(X))$ by induction on $0\leq i\leq n$ as follows. Fix $0\leq i\leq n$, and assume that $x=u_0\vert\cdots\vert u_n$ is $j$-critical for every $0\leq j<n$. Call~$x$ 
\begin{enumerate} 
\item \emph{$i$-subcritical} if $\delta_i(x)=(1,u_i)$,
\item \emph{$i$-critical} if $\delta_i(x)=(u_i,1)$, 
\item \emph{$i$-supercritical} otherwise.
\end{enumerate}
An $n$-cell of $\Std(\red(X))$ is \emph{critical} if it is $i$-critical for every $0\leq i\leq n$, and \emph{subcritical} (resp.\ \emph{supercritical}) if it is $i$-subcritical (resp.\ $i$-supercritical) for some $0\leq i\leq n$. 

\begin{remark}
A nontrivial monomial~$u_0$ of~$\lin{X}$ is either critical (if it is of length~$1$) or supercritical (otherwise), but never subcritical.

Let~$u_0\vert u_1$ be a $1$-cell of $\Std(\red(X))$ such that~$u_0$ is critical. By definition of~$\rho$, we have the inequality $\rho_{u_0}u_1 \sqsubseteq \rho_{u_0u_1}$, with $\rho_{u_0}u_1=1_{u_0u_1}$ by hypothesis on~$u_0$. Then $u_0\vert u_1$ is subcritical if, and only if, the inequality is an equality, i.e.~$u_0u_1$ is reduced. Otherwise, $u_1$ being reduced, $\rho_{u_0u_1}=\alpha w$ for a $1$-cell~$\alpha$ of~$X$ and a right-factor~$w$ of~$u_1$: then~$u_0\vert u_1$ is critical if~$w=1$, and supercritical otherwise.

Finally, assume that $x=u_0\vert u_1\vert u_2$ is a $2$-cell of $\Std(\red(X))$, with~$u_0\vert u_1$ critical. Denote by~$b$ the branching $(\rho_{u_0u_1}u_2,\rho_{u_0u_1u_2})$ of~$X$. Then, by definition, $x$ is subcritical if, and only if, $\rho_{u_0u_1}u_2 = \rho_{u_0u_1u_2}$, i.e.~$b$ is aspherical. Otherwise, $u_2$ being reduced, the branching~$b$ is an overlap of shape $(\alpha u_2, w\beta w')$. Then~$x$ is critical in the case~$w'=1$, i.e.\ when the branching~$b$ is critical, and supercritical otherwise.
\end{remark}

\begin{theorem}[Squier's polygraphic resolution]
\label{T:SquierResolution}
Let~$A$ be an augmented algebra and~$X$ be a reduced convergent left-monomial presentation of~$A$. 
Define~$\Sq(X)$ to be the indexed subset of $\Std(\red(X))$ consisting of its critical cells, and~$\phi$ as the indexed partial map of degree~$1$ from~$\Sq(X)$ to itself, given by
\[
\phi(u_0\vert\cdots\vert u_n) = 
	u_0\vert\cdots\vert u_{i-1}\vert 
	v_i\vert w_i\vert
	u_{i+1}\vert\cdots\vert u_n
\]
for every $i$-supercritical $n$-cell $u_0\vert\cdots\vert u_n$ with $\delta_i(u_0\vert\cdots\vert u_n)=(v_i,w_i)$. 

Then~$\phi$ is a collapsing scheme of $\Std(\red(X))$ onto~$\Sq(X)$. As a consequence, $\Sq(X)$, equipped with the structure of $\infty$-polygraph induced by~$\phi$, is a polygraphic resolution of~$A$.
\end{theorem}

\begin{proof}
First, $\phi$ induces a bijection between the $i$-supercritical $n$-cells and the $(i+1)$-subcritical $(n+1)$-cells of $\Std(\red(X))$. Indeed, if $u_0\vert\cdots\vert u_{n+1}$ is an $(i+1)$-subcritical $(n+1)$-cell, put $\psi(u_0\vert\cdots\vert u_{n+1})=u_0\vert\cdots\vert u_i u_{i+1}\vert\cdots\vert u_n$. This assignment makes sense because, if $\delta_{i+1}(u_0\vert\cdots\vert u_{n+1})=(1,u_{i+1})$, then~$u_iu_{i+1}$ is reduced, hence $u_0\vert\cdots\vert u_i u_{i+1}\vert\cdots\vert u_n$ is an $i$-supercritical $n$-cell. Moreover, $\psi$ is inverse to~$\phi$.

Second, the partition $\Std(\red(X))=\im(\phi)\amalg\Sq(X)\amalg\dom(\phi)$ is obtained by observing that~$\Sq(X)$, $\dom(\phi)$ and $\im(\phi)$ respectively consist of the critical, supercritical and subcritical cells of $\Std(\red(X))$.

Third, let $u_0\vert\cdots\vert u_n$ be an $i$-subcritical $n$-cell of $\Std(\red(X))$. Note that the definition of~$\phi$ implies $u_0\vert\cdots\vert u_n=\phi(u_0\vert\cdots\vert u_i u_{i+1}\vert\cdots\vert u_n)$. By Corollary~\ref{C:ConvergentStandardPolygraphicResolution}, the boundary of $u_0\vert\cdots\vert u_n$ satisfies
\[
\dr(u_0\vert\cdots\vert u_n) \quad= 
	\sum_{0\leq j\leq n+1} (-1)^{n+1-j} \: d_j(u_0\vert\cdots\vert u_n) 
	\quad+\quad \dr'(u_0\vert\cdots\vert u_n),
\]
where $\dr'(u_0\vert\cdots\vert u_n)$ is an identity. Observe that $u_0\vert\cdots\vert u_i u_{i+1}\vert\cdots\vert u_n=d_{i+1}(u_0\vert\cdots\vert u_n)$, and that this $(n-1)$-cell is distinct from each other $(n-1)$-cell of $\Std(\red(X))$ that appears in $\dr(u_0\vert\cdots\vert u_n)$. Hence, $\dr(u_0\vert\cdots\vert u_n)$ has the required form.

Finally, check that the relation~$\trieq_{\phi}$ induced by~$\phi$ is wellfounded, by proving that it is included into a wellfounded order~$\trieq$. Put $u_0\vert\cdots\vert u_n \tri v_0\vert\cdots\vert v_n$ if either
\begin{enumerate}
\item $v_0\cdots v_n$ is a proper submonomial of $u_0\cdots u_n$, or
\item $u_0\cdots u_n\succ_X a$, with~$a$ an $n$-cell of $\lin{\Std(\red(X))}$ such that $v_0\vert\cdots\vert v_n$ belongs to~$\supp(a)$, or
\item $u_0\cdots u_n=v_0\cdots v_n$, and there exists~$i$ in~$\ens{0,\dots,n}$ such that $u_0=v_0$, \dots, $u_{i-1}=v_{i-1}$, and $l(u_i)>l(v_i)$.
\end{enumerate}
The order relation~$\trieq$ so defined is wellfounded as a lexicographic product of wellfounded orders. Now, let $u_0\vert\cdots\vert u_n$ be an $i$-subcritical $n$-cell of $\Std(\red(X))$, and prove $u_0\vert\cdots\vert u_i u_{i+1}\vert \cdots\vert u_n \tri y$ for every other $(n-1)$-cell $y$ in~$\supp(\dr(u_0\vert\cdots\vert u_n))$, that is, in the support of a $d_j(u_0\vert\cdots\vert u_n)$, for a $j\neq i+1$. Consider the possible values of~$j$.

If~$j=0$ or~$j=n+1$, the only possible $(n-1)$-cells for~$y$ are $u_1\vert\cdots\vert u_n$ or $u_0\vert\cdots\vert u_{n-1}$. These two$(n-1)$-cells are proper submonomials of $u_0\cdots u_n$, so that Condition~\Item{i} applies.

If~$1\leq j\leq n$, then~$y$ is in the support of $u_1\vert\cdots\vert \rep{u_ju_{j+1}}\vert \cdots \vert u_n$, i.e.\ $y=u_0\vert\cdots\vert u_{j-1} \vert v\vert u_{j+2}\vert\cdots u_n$, for~$v\in\supp(\rep{u_ju_{j+1}})$. 
If~$u_ju_{j+1}$ is not reduced, then $u_0\cdots u_{j-1}v u_{j+2}\cdots u_n$ is not a proper submonomial of $u_0\cdots u_n$, because, otherwise, $X$ would not terminate, so that Condition~\Item{ii} applies. 
On the contrary, if~$u_ju_{j+1}$ is reduced, then none of Conditions~\Item{i} and~\Item{ii} is satisfied, but Condition~\Item{iii} applies, because $u_0\vert\cdots\vert u_{i-1}$ critical implies that~$u_k u_{k+1}$ is nonreduced for each~$k<i$, so that~$j>i+1$ holds.
\end{proof}

\subsection{Examples: the symmetric algebra and variations}
\label{SS:ExamplesSymmetricAlgebra}

\begin{example}
\label{X:SymmetricAlgebra3Generators}

Consider the symmetric algebra on three generators~$x$, $y$ and~$z$, presented by 
\[
\Sym(x,y,z) = \bigpres{x, y, z}{yx \: \ofl{\alpha} \: xy,\, zx \: \ofl{\beta} \: xz,\, zy \: \ofl{\gamma} \: yz}.
\]
The $1$-polygraph $\Sym(x,y,z)$ is convergent, with one critical branching of source~$zyx$, giving rise, by Theorem~\ref{T:Squier}, to an acyclic extension with exactly one $2$-cell:
\[
\xymatrix @R=0.75em {
& yzx    
  \ar[r] ^{y \beta} _{}="sF" 
& yxz
  \ar @/^/ [dr] ^{\alpha z}
\\
zyx
	\ar @/^/ [ur] ^{\gamma x}
	\ar @/_/[dr] _{z \alpha}
&&& xyz
\\
& zxy
  \ar[r] _{\beta y} ^{}="tF"
& xzy
  \ar @/_/ [ur] _{x \gamma}
\ar@2 "sF"!<0pt,-15pt>;"tF"!<0pt,15pt> ^-{\omega}
}
\]
Note that $\Sym(x,y,z)$ is also reduced, so that Theorem~\ref{T:SquierResolution} applies to it: we obtain a polygraphic resolution $\Sq(\Sym(x,y,z))$ with a unique $2$-cell and no $n$-cell for~$n\geq 3$, and the boundary of the $2$-cell is obtained by collapsing the standard polygraphic resolution, through the collapsing scheme~$\phi$ of Theorem~\ref{T:SquierResolution}. 
The cells required to understand the computation are drawn on the following figure:
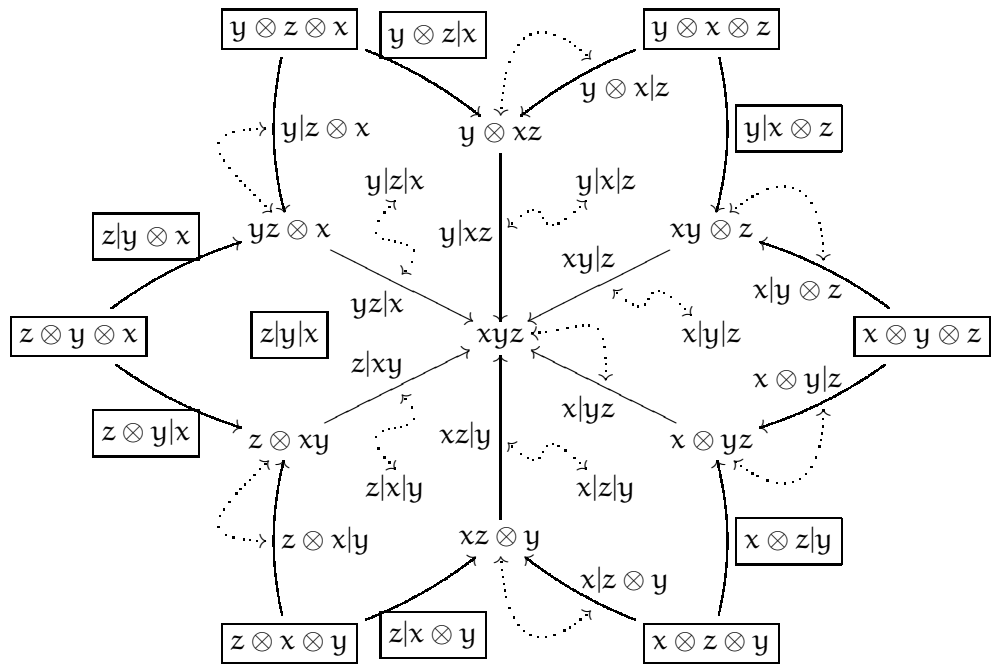
\[
\xymatrix @!C @C=2em @R=2em {
& {\boxed{y\otimes z\otimes x}}
	\ar@/_/ [dd] ^{y\vert z\otimes x} _{}="f2" 
	\ar@/^/ [dr] ^{\boxed{y \otimes z\vert x}}
&& {\boxed{y\otimes x\otimes z}}
	\ar@/_/ [dl] ^{y \otimes x\vert z} _{}="f1"  
	\ar@/^/ [dd] ^{\boxed{y\vert x \otimes z}}
\\
&& y\otimes xz 
	\ar [dd] _{y\vert xz} ^{}="f4"  
	\ar@{} [dr] |{y\vert x\vert z} |{}="f4a"
\\
& yz\otimes x 
	\ar [dr] _{yz\vert x} ^{}="f3"  
	\ar@{} [ur] |{y\vert z\vert x} |{}="f3a"
&& xy\otimes z  
	\ar [dl] _{xy\vert z} ^{}="f6" 
\\
{\boxed{z\otimes y\otimes x}}
	\ar@/^/ [ur] ^{\boxed{z\vert y\otimes x}}
	\ar@/_/ [dr] _{\boxed{z\otimes y\vert x}}
	\ar@{} [rr] |{\boxed{z\vert y\vert x}}
&& xyz 
	\ar@{} [rr] |{x\vert y\vert z}
&& {\boxed{x\otimes y\otimes z}}
	\ar@/_/ [ul] ^{x\vert y\otimes z} _{}="f5"
	\ar@/^/ [dl] _{x\otimes y\vert z} ^{}="f9" 
\\
& z\otimes xy 
	\ar [ur] ^{z\vert xy} _{}="f8"
	\ar@{} [dr] |{z\vert x\vert y} |{}="f8a"
&& x\otimes yz 
	\ar [ul] ^{x\vert yz} _{}="f7" 
\\
&& xz\otimes y 
	\ar [uu] ^{xz\vert y} _{}="f11" 
	\ar@{} [ur] |{x\vert z\vert y} |{}="f11a"
\\
& {\boxed{z\otimes x\otimes y}}
	\ar@/^/ [uu] _{z \otimes x|y} ^{}="f10"
	\ar@/_/ [ur] _{\boxed{z\vert x\otimes y}}
&& {\boxed{x\otimes z\otimes y}}
	\ar@/_/ [uu] _{\boxed{x \otimes z|y}}
	\ar@/^/ [ul] _{x|z \otimes y} ^{}="f12" 
\ar@{<.>}@(u,ul) "2,3";"f1"
\ar@{<.>}@(ul,l) "3,2";"f2"
\ar@{<.>}@(ur,dl) "f3";"f3a"!<0pt,-7.5pt>
\ar@{<.>}@(ur,dl) "f4";"f4a"!<-7.5pt,-7.5pt>
\ar@{<.>}@(ur,u) "3,4";"f5"
\ar@{<.>}@(dr,ul) "f6";"4,4"!<-7.5pt,7.5pt>
\ar@{<.>}@(r,u) "4,3";"f7"
\ar@{<.>}@(dr,ul) "f8";"f8a"!<0pt,7.5pt>
\ar@{<.>}@(dr,d) "5,4";"f9"
\ar@{<.>}@(dl,l) "5,2";"f10"
\ar@{<.>}@(dr,ul) "f11";"f11a"!<-7.5pt,7.5pt>
\ar@{<.>}@(d,dl) "6,3";"f12"
}
\]
The identifications induced by~$\phi$ are given by the dotted arrows, and the cells that belong to $\Sq(\Sym(x,y,z))$ are boxed. After the collapse, by identification of~$y\vert x$, $z\vert x$ and~$z\vert y$ to~$\alpha$, $\beta$ and~$\gamma$, respectively, and by omission of the~$\tens$ sign, we recover the former $2$-cell~$\omega$  as $z \vert y \vert x$.
\end{example}

\begin{example}
\label{X:SymmetricAlgebra}

Consider the symmetric algebra on a set~$X$ of generators: given a fixed linear order on~$X$, it admits a convergent presentation~$\Sym(X)$ whose relations are the $\alpha_{xy} : yx \fl xy$, for all $x<y$ in~$X$. The polygraphic resolution obtained by Theorem~\ref{T:SquierResolution} has exactly one $n$-cell $\omega_{x_0\cdots x_n}$ for all $x_0<\dots<x_n$ of~$X$, that corresponds to the critical $n$-cell $x_n\vert\cdots\vert x_0$ of the standard resolution. 
One can compute that Squier's resolution $\Sq(\Sym(X))$ has one $3$-cell $\omega_{xyz}$ for each triple $x<y<z$ of elements of~$X$, similar to the permutohedron $3$-cell~$\omega$ of Example~\ref{X:SymmetricAlgebra3Generators}, and one $4$-cell $\omega_{xyzt}$ for each quadruple $x<y<z<t$, whose source is 
\[
\xymatrix @!C @C=2em @R=1.5em {
&&& yztx \ar[dr] 
	\ar@{}[d] |(0.6){\sm=}  
\\
&& zytx \ar[ur] \ar[r] 
& yztx \ar[r]
& yzxt \ar[dr]
\\
& ztyx \ar[ur]
& \omega_{tzy}x 
& ytzx \ar[u] \ar[dr]
& y\omega_{tzx}
& yxzt \ar[dr]
\\
tzyx \ar[ur] \ar[dr] \ar[rr]
&& tyzx \ar[dr] \ar[ur] 
	\ar@{}[rr] |{\sm=} 
&& ytxz \ar[r]
& yxtz \ar[u] \ar[d]
	\ar@{}[r] |(0.4){\sm=}
& xyzt 
\\
& tzxy \ar[dr]
& t\omega_{zyx} 
& tyxz \ar[d] \ar[ur]
&\omega_{tyx}z 
& xytz \ar[ur]
\\
&& txzy \ar[r]
& txyz \ar[r]
& xtyz \ar[ur]
}
\]
and whose target is
\[
\xymatrix @!C @C=2em @R=1.5em {
&& zytx \ar[r]
& zyx \ar[r] \ar[d]t
& yzx \ar[dr]t 
\\
& ztyx \ar[ur] \ar[d]
& z\omega_{txy} 
& zxyt \ar[dr]
& \omega_{zyx}t 
& yxzt \ar[dr]
\\
tzyx \ar[dr] \ar[ur]
	\ar@{}[r] |(0.6){\sm=}
& ztxy \ar[r]
& zxty \ar[ur] \ar[dr]
	\ar@{}[rr] |{\sm=}
&& xzyt \ar[rr]
&& xyzt
\\
& tzxy \ar[dr] \ar[u]
& \omega_{tzx}y
& xzty \ar[ur]
& x\omega_{tzy}
& xytz \ar[ur]
\\
&& txzy \ar[dr] \ar[r]
& xtzy \ar[r] \ar[u]
	\ar@{}[d] |(0.4){\sm=}
& xty \ar[ur]z
\\
&&& txyz \ar[ur]
}
\]

More generally, fix $n\geq 1$ and $x_0<\dots<x_n$ in~$X$, and denote by $\omega_{x_0\cdots x_n}$ the $n$-cell of $\Sq(\Sym(X))$ that is the image of $x_n\vert\cdots\vert x_0$ through the projection~$\pi$ given by~\eqref{E:CollapsingProjection}. Corollary~\ref{C:ConvergentStandardPolygraphicResolution} tells us that the source and target of the $n$-cell $x_n\vert\cdots\vert x_0$ of $\Std(\red(\Sym(X)))$ are given by
\begin{align*}
s(x_n\vert\cdots\vert x_0) 
	&\quad= \sum_{\substack{0\leq i\leq n+1 \\ i \text{ even}}}
	d_i(x_n\vert\cdots\vert x_0) \quad+\quad a \\
\text{and}\qquad
t(u_0\vert\cdots\vert u_n) 
	&\quad= \sum_{\substack{0\leq i\leq n+1 \\ i \text{ odd}}}
	d_i(x_n\vert\cdots\vert x_0) \quad+\quad b \:,
\end{align*}
where~$a$ and~$b$ are identities, and  
\[
d_i(x_n\vert\cdots\vert x_0) =
\begin{cases}
x_n\tens x_{n-1} \vert\cdots\vert x_0 
	&\text{if $i=0$,} 
\\
x_n\vert\cdots\vert x_i x_{i+1} \vert\cdots\vert x_0
	&\text{if $1\leq i\leq n$,}
\\
x_n\vert\cdots\vert x_1\tens x_0
	&\text{if $i=n+1$.}
\end{cases}
\]
Observe that $d_0(x_n\vert\cdots\vert x_0)$ and $d_{n+1}(x_n\vert\cdots\vert x_0)$ belong to $\lin{\Sq(\Sym(X))}$, because $x_{n-1} \vert\cdots\vert x_0$ and~$x_n\vert\cdots\vert x_1$ are critical, whereas, for $1\leq i\leq n$, the $(n-1)$-cell $d_i(x_n\vert\cdots\vert x_0)$ is supercritical. 
Now, fix~$1\leq i\leq n$. By definition of~$\pi$, the $(n-1)$-cell $x_n\vert\cdots\vert x_i x_{i+1} \vert\cdots\vert x_0$ is mapped to the $(n-1)$-cell $c$ such that 
\[
 x_n\vert\cdots\vert x_i x_{i+1} \vert\cdots\vert x_0 = \dr(\phi(x_n\vert\cdots\vert x_i x_{i+1} \vert\cdots\vert x_0)) + c.
\]
Note that $\phi(x_n\vert\cdots\vert x_i x_{i+1} \vert\cdots\vert x_0)=x_n\vert\cdots\vert x_i \vert x_{i+1} \vert\cdots\vert x_0$. As before, the boundary of the latter $n$-cell is an alternated sum of $d_j(x_n\vert\cdots\vert x_i \vert x_{i+1} \vert\cdots\vert x_0)$, for $0\leq j\leq n+1$, plus an identity. A careful examination of the different possibilities shows that, for $j>i+1$ or $j<i-1$, the $(n-1)$-cell $d_j(x_n\vert\cdots\vert x_i \vert x_{i+1} \vert\cdots\vert x_0)$ is subcritical: so, $\pi$ maps it to an identity. Repeating the process for each of $j=i+1$ and $j=i-1$ gives, by induction respectively upwards from~$i+1$ to~$n+1$ and downwards from~$i-1$ to~$0$, that 
\begin{align*}
\pi(d_{i+1}(x_n\vert\cdots\vert x_i\vert x_{i+1} \vert\cdots\vert x_0)) 
	&= x_i\tens x_n\vert\cdots\vert\rep{x_i}\vert\cdots\vert x_0 + a
\\[1ex]
\text{and}\quad
\pi(d_{i-1}(x_n\vert\cdots\vert x_i\vert x_{i+1} \vert\cdots\vert x_0)) 
	&= x_n\vert\cdots\vert\rep{x_{i+1}}\vert\cdots\vert x_0 \tens x_{i+1} + b,
\end{align*}
where~$a$ and~$b$ are some identities, and~$\rep{x_i}$ means that~$x_i$ is omitted (and not a normal form). So, we have
\[
\pi(d_i(x_n\vert\cdots\vert x_i\vert x_{i+1} \vert\cdots\vert x_0)) 
	= x_i\tens x_n\vert\cdots\vert\rep{x_i}\vert\cdots\vert x_0
	+ x_n\vert\cdots\vert\rep{x_{i+1}}\vert\cdots\vert x_0 \tens x_{i+1}
	+ a
\]
for some identity~$a$. 
Wrapping up the results, we obtain, for some identity~$a$:
\begin{equation}
\label{E:BoundarySymmetricAlgebra}
\dr (\omega_{x_0\cdots x_n}) 
	\quad= \sum_{0\leq i\leq n+1}
	(-1)^{i+1} \big (
	x_i \omega_{x_0\cdots\rep{x_i}\cdots x_n} 
	\:-\: \omega_{x_0\cdots\rep{x_i}\cdots x_n} x_i
	\big)
	\quad+\quad a \:.
\end{equation}
\end{example}

\begin{example}
\label{X:QuanticSymmetricAlgebra}

Consider the quantum deformation of the symmetric algebra, whose generators commute up to some constant parameters, i.e.\ $yx = q_{xy} xy$, where~$q_{xy}$ is a scalar. This algebra, on a set~$X$ of generators equipped with a linear order, is presented by the convergent $1$-polygraph 
\[
\bigpres{X}{\big(\,yx \ofl{\alpha_{xy}} q_{xy} xy\,\big)_{x<y}}.
\]
Squier's resolution, given by Theorem~\ref{T:SquierResolution}, is similar to the one of the symmetric algebra of Example~\ref{X:SymmetricAlgebra}. For example, it has one $2$-cell~$\omega_{xyz}$ for each triple $x<y<z$ of generators:
\[
\xymatrix @R=0.75em {
& q_{yz} yzx    
  \ar[rr] ^-{q_{yz} y \alpha_{xz}} _{}="sF" 
&& q_{xz}q_{yz} yxz
  \ar @/^2ex/ [dr] ^-{q_{xz}q_{yz} \alpha_{xy} z}
\\
zyx
	\ar @/^/ [ur] ^-{\alpha_{yz} x}
	\ar @/_/[dr] _-{z \alpha_{xy}}
&&&& q_{xy}q_{xz}q_{yz} xyz
\\
& q_{xy} zxy
  \ar[rr] _-{q_{xy}\alpha_{xz} y} ^{}="tF"
&& q_{xy}q_{xz} xzy
  \ar @/_2ex/ [ur] _-{q_{xy}q_{xz} x \alpha_{yz}}
\ar@2 "sF"!<0pt,-15pt>;"tF"!<0pt,15pt> ^-{\omega_{xyz}}
}
\]
\end{example}

\begin{example}
Consider the exterior algebra on a set~$X$, equipped with an arbitrary linear order. As an associative algebra, it is presented by the $1$-polygraph
\[
\bigpres{X}{\big(\,yx \ofl{\alpha_{xy}} -xy\,\big)_{x<y}, \: \big(\,xx \ofl{\delta_x} 0\,\big)_x}.
\]
This $1$-polygraph is convergent, with one critical branching of source~$zyx$ for each triple $x<y<z$, two of source~$yyx$ and~$yxx$ for each pair~$x<y$, and one of source~$xxx$ for each generator~$x$. The first one gives rise to a $2$-cell with the same shape as~$\omega_{xyz}$ in Example~\ref{X:QuanticSymmetricAlgebra}, with $q_{xy}=q_{xz}=q_{yz}=-1$. The critical branching of source~$yxx$ (the one of source~$yyx$ being symmetric) induces a $2$-cell
\[
\xymatrix @!C @C=3em @R=1em {
& -xyx   
  \ar[r] ^{-x\alpha_{xy}} 
& xxy
  \ar @/^/ [dr] ^{\delta_x y}
\\
yxx
	\ar @/^/ [ur] ^{\alpha_{xy}x}
	\ar @/_2ex/ [rrr] _{y\delta_x}
&&& 0
\ar@2 "1,2"!<35pt,-12.5pt>;"2,2"!<35pt,0pt>
}
\]
which is obtained by collapsing the following cells of the standard polygraphic resolution, after identification of~$y\vert x$ with~$\alpha_{xy}$ and~$x\vert x$ with~$\delta_x$, and removal of the~$\tens$ sign:
\[
\xymatrix @!C @C=1.5em {
&& {\boxed{-x\otimes y\otimes x}}
	\ar@/_/ [dl] _-{-x\vert y\otimes x} _-{}="1"
	\ar@/^/ [dr] ^-{\boxed{-x\otimes y\vert x}}
	\ar@{} [dd] |(0.45){-x\vert y\vert x} |(0.45){}="2"
\\
& -xy\otimes x
	\ar [dr] _-{-xy\vert x} _-{}="3"
	\ar@{} [d] |(0.7){\boxed{y\vert x\vert x}}
&& x\otimes xy
	\ar [dl] _-{x\vert xy} _-{}="4"
	\ar@{} [d] |(0.7){x\vert x\vert y} |(0.7){}="5"
\\
{\boxed{y\otimes x\otimes x}}
	\ar@/^/ [ur] ^-{\boxed{y\vert x\otimes x}} 
	\ar@/_2ex/ [rr] _-{\boxed{y\otimes x\vert x}}
&& {\boxed{0}}
&& {\boxed{x\otimes x\otimes y}}
	\ar@/_/ [ul] _-{x\otimes x\vert y} ^-{}="6"
	\ar@/^2ex/ [ll] ^-{\boxed{x\vert x\otimes y}}
	\ar@{<.>}@(ul,ul) "2,2"!<-15pt,0pt>;"1"!<-25pt,17.5pt>
	\ar@{<.>}@(d,u) "2"!<0pt,-7.5pt>;"3"
	\ar@{<.>}@(dr,ul) "4";"5"!<0pt,7.5pt>
	\ar@{<.>}@(ur,u) "2,4";"6"!<0pt,15pt>
}
\]
Finally, the critical branching of source~$xxx$ produces the $2$-cell
\[
\xymatrix @!C @C=3em {
xxx
	\ar@/^3ex/ [r] ^{\delta_x x} _{}="s"
	\ar@/_3ex/ [r] _{x\delta_x} ^{}="t"
	\ar@2 "s"!<0pt,-10pt>;"t"!<0pt,10pt>
& 0
}
\]
In higher dimensions, Theorem~\ref{T:SquierResolution} produces a polygraphic resolution of the exterior algebra on~$X$, with one $n$-cell corresponding to $x_n\vert\dots\vert x_0$ for all $x_n\geq\cdots\geq x_0$ in~$X$.
\end{example}

\section{Free resolutions of associative algebras}
\label{S:FreeResolutions}

In this final section, we show that every polygraph gives rise to a chain complex of free bimodules, with the same generators, and whose differential is built from the source and target maps of the polygraph. Moreover, Theorem~\ref{T:PolRes->AbRes} proves that, starting with a polygraphic resolution of an algebra~$A$, this construction yields a resolution of~$A$ by free $A$-bimodules. The proof relies on building a contracting homotopy of the chain complex from a contraction of the polygraph. This method is applied to obtain sufficient or necessary conditions for an algebra to be Koszul.

\subsection{Free bimodules resolutions from polygraphic resolutions}
\label{SS:FreeBimodulesResolutions}

\subsubsection{Free bimodules}

Let~$X$ be an $\infty$-polygraph, and~$A$ be the algebra presented by the underlying $1$-polygraph of~$X$. Denote by $\env{A}=\op{A}\otimes A$ the enveloping algebra of~$A$, and, for~$n\geq 0$, by $\env{A}[X_n]$ the free $A$-bimodule over~$X_n$. An element $a\otimes\alpha\otimes b$ of $\env{A}[X_n]$ is written $a[\alpha]b$. 
The inclusion $\alpha\mapsto[\alpha]$ of~$X_n$ into $\env{A}[X_n]$ is extended to all the $n$-cells of~$\lin{X}$ as follows. In dimension~$0$, define~$[u]$ for every monomial~$u$ of~$\lin{X}$, by induction on the length of~$u$, by
\[
[1] \:=\: 0 
\qquad\text{and}\qquad
[uv] \:=\: [u]\cl{v} + \cl{u}[v],
\]
and, then, extend the bracket to all $0$-cells by linearity. In dimension~$n\geq 1$, define the bracket on identities and on $n$-monomials by
\[
[1_a] \:=\: 0
\qquad\text{and} \qquad
[u\alpha v] \:=\: \cl{u}[\alpha]\cl{v}.
\]
This definition, extended by linearity to linear combinations of identities and $n$-monomials is compatible with~\eqref{E:nAlgExchangeProd} in Theorem~\ref{T:nAlg}, yielding a well-defined linear map on the $n$-cells of~$\lin{X}$. Then~(\ref{E:nVectComp}) and~(\ref{E:nVectInv}) in Proposition~\ref{P:nVect} imply, for any possible $n$-cells~$a$ and~$b$ in~$\lin{X}$:
\[
[a\star_k b] \:=\: [a] + [b]
\qquad
\text{and}
\qquad
[a^-]\:=\:-[a].
\]

\subsubsection{The chain complex of a polygraph}

Assume that~$X$ is an $\infty$-polygraph, and let~$A$ be the algebra presented by~$X$. Let~$\env{A}[X]$ be the chain complex of $A$-bimodules
\[
0 \fllg 
A 
\oflg{\mu} 
\env{A}
\oflg{\delta_{0}}
\env{A}[X_0]
\oflg{\delta_{1}}
\; \cdots \;
\fllg
\env{A}[X_k]
\oflg{\delta_{k+1}}
\env{A}[X_{k+1}]
\fllg
\;\cdots
\]
with $\mu(a\otimes b) = ab$, and, for every~$k\geq 0$, 
\[
\delta_k[\alpha] =
\begin{cases}
1 \otimes \alpha  - \alpha\otimes 1
	&\text{if~$k=0$,} \\
[s_{k-1}(\alpha)] - [t_{k-1}(\alpha)]
	&\text{if~$k\geq 1$.}
\end{cases}
\]
Let us check that $\env{A}[X]$ is indeed a chain complex. First, prove by induction on the length that $\delta_{0}([u]) = 1\otimes u - u\otimes 1$ holds for every monomial~$u$ of~$\lin{X}$, yielding $\mu\delta_{0} = 0$.
Then, for every $(n+1)$-cell~$\alpha$ of~$X$, with~$n\geq 1$, the globular relations imply
\[
\delta_n\delta_{n+1}[\alpha] 
	= [s_{n-1}s_n(\alpha)] + [t_{n-1}s_n(\alpha)] - [s_{n-1}t_n(\alpha)] - [t_{n-1}t_n(\alpha)] 
	= 0.
\] 

\begin{theorem}
\label{T:PolRes->AbRes}
If~$X$ is a polygraphic resolution of an algebra~$A$, then the complex~$\env{A}[X]$ is a free resolution of the $A$-bimodule~$A$. Moreover, if~$X$ is of finite type, then so is~$\env{A}[X]$.
\end{theorem}

\begin{proof}
Fix a unital section~$\iota$ of~$X$, writing~$\rep{a}$ for~$\iota(a)$, and a right $\iota$-contraction~$\sigma$ of~$X$, thanks to Proposition~\ref{T:ResolutionContraction}. Using~$\sigma$, let us construct a contracting homotopy~$h$ of~$\env{A}[X]$. Define $h_{-1} : A \fl \env{A}$ and $h_0 : \env{A} \fl \env{A}[X_0]$ by 
\[
h_{-1}(a) = a\otimes 1
\quad
\text{and}
\quad
h_0(a\otimes b) = a[\rep{b}],
\]
for all~$a$ and~$b$ in~$A$, and, for~$n\geq 1$, define a morphism of $A$-modules $h_n : \env{A}[X_{n-1}] \fl \env{A}[X_n]$ by 
\[
h_n([\alpha]a) = [\sigma_{\alpha\rep{a}}],
\]
for every~$a$ in~$A$, and every $(n-1)$-cell $\alpha$ of~$X$. 

Now, check that $h_{n+1}([a]b) = [\sigma_{a\rep{b}}]$ holds for every $n$-cell~$a$ of~$\lin{X}$, and every~$b$ in~$A$. 
If~$a=u$ is a monomial of~$\lin{X}$, we prove, by induction on the length of~$u$, that $h_1([u]a)=[\sigma_{u\rep{a}}]$ holds. If~$u=1$, then both members of the equality are~$0$. If~$u$ and~$v$ are monomials of~$\lin{X}$, then, on the one hand, by definition of the bracket,
\[
h_1([uv]a) = \cl{u}[\sigma_{v\rep{a}}] + [\sigma_{u\rep{va}}],
\]
and, on the other hand, using the fact that~$\sigma$ is a right $\iota$-contraction,
\[
[\sigma_{uv\rep{a}}] = \cl{u}[\sigma_{v\rep{a}}] + [\sigma_{u\rep{va}}].
\]
For~$n\geq 1$, the $n$-cells of~$\lin{X}$ are linear combinations of the $n$-monomials of~$\lin{X}$ and of an identity $n$-cell of~$\lin{X}$. Suppose that~$a$ is an $n$-monomial $u\alpha v$, where~$u$ and~$v$ are monomials of~$\lin{X}$ and~$\alpha$ is an $n$-cell of~$X$. On the one hand, we have
\[
h_{n+1}([u\alpha v]a)
= h_{n+1}(\cl{u}[\alpha]\cl{v}a)
= \cl{u}[\sigma_{\alpha\rep{va}}].
\]
On the other hand, \eqref{E:RightContraction} gives, with $b=s_0(\alpha)$, 
\[
\sigma_{u\alpha v\rep{a}} 
	= ub\sigma_{v\rep{a}} 
	\star_0 u \sigma_{\alpha\rep{va}}
	\star_0 \sigma_{u\rep{bva}},
\]
from which we obtain, using the fact that $\sigma_{v\rep{a}}$ and $\sigma_{u\rep{bva}}$ are identities,
\[
[\sigma_{u\alpha v\rep{a}}] 
	= \cl{ub}[\sigma_{v\rep{a}}] 
	+ \cl{u}[\sigma_{\alpha\rep{va}}]
	+ [\sigma_{u\rep{bva}}]
	=\cl{u}[\sigma_{\alpha\rep{va}}].
\]

Finally, prove that~$h$ is a contracting homotopy for the complex~$\env{A}[X]$. We have $h_{-1}\mu(a\otimes b) = ab\otimes 1$ and $\delta_{0}h_0(a\otimes b) =  a\otimes b - ab \otimes 1$, thus $h_{-1}\mu + \delta_{0}h_0 = \id_{\env{A}}$. Then, for every~$[x]a$ in~$\env{A}[X_0]$, we have $h_0\delta_{0}([x]a) = [\rep{xa}] - \cl{x}[\rep{a}]$ and 
\[
\delta_1 h_1([x]a) = \delta_1([\sigma_{x\rep{a}}]) = [x\rep{a}] - [\rep{xa}] = \cl{x}[\rep{a}] + [x]a - [\rep{xa}],
\]
which give $h_0\delta_{0} + \delta_1h_1 = \id_{\env{A}[X_0]}$.
Now, for~$n\geq 1$ and $[\alpha]a$ be in $\env{A}[X_n]$, we have
\begin{align*}
\delta_n h_{n+1}([\alpha]a)
= \delta_n([\sigma_{\alpha\rep{a}}])
&= [\alpha\rep{a}] + [\sigma_{t_{n-1}(\alpha)\rep{a}}] - [\sigma_{s_{n-1}(\alpha)\rep{a}}],\\
&= [\alpha]a + h_{n-1}([t_{n-1}(\alpha)]a) - h_{n-1}([s_{n-1}(\alpha)]a).
\end{align*}
Thus $h_{n-1}\delta_{n-1} + \delta_n h_n = \id_{\env{A}[X_n]}$.
\end{proof}

\begin{example}
\label{X:SymmetricAlgebraAbelian}

Consider the polygraphic resolution $\Sq(\Sym(X))$ of the symmetric algebra on a set~$X$, obtained in Example~\ref{X:SymmetricAlgebra}. The associated free resolution of the symmetric algebra is generated by one element $[\omega_{x_0\cdots x_n}]$ for every tuple $x_n>\cdots>x_0$ in~$X$. From~\eqref{E:BoundarySymmetricAlgebra}, we deduce
\[
\delta_n [\omega_{x_0\cdots x_n}]
	\quad= \sum_{0\leq i\leq n+1}
	(-1)^{i+1} \big (
	x_i [\omega_{x_0\cdots\rep{x_i}\cdots x_n}]
	\:-\: [\omega_{x_0\cdots\rep{x_i}\cdots x_n}] x_i
	\big).
\]
Thus, we recover (up to the sign of the differential) the usual Koszul bimodule complex of the symmetric algebra over~$X$. In the case of a quantum deformation of the symmetric algebra of Example~\ref{X:QuanticSymmetricAlgebra}, we get the quantum version of the Koszul complex, obtained for instance by Wambst in~\cite{Wambst93}.
\end{example}

\begin{remark}
The construction of the resolution~$\env{A}[X]$ of $A$-bimodules can be adapted to obtain a resolution of~$\K$ by right $A$-modules. For~$n\geq 0$, denote by $\op{A}[X_n]$ the free right $A$-module generated by~$X_n$. The mapping inclusion $\alpha\mapsto[\alpha]$ of~$X_n$ into~$\op{A}[X_n]$ is extended to all the $n$-cells of~$\lin{X}$ by setting
\[
[1] = 0,
\quad
[xu] = [x]\cl{u} + [u],
\quad
[1_u] = 0,
\quad
[x\alpha] = 0,
\quad
[\alpha u] =[\alpha]\cl{u},
\]
for every $0$-cell~$x$ of~$X$, every $n$-cell~$\alpha$ of~$X$, with~$n\geq 1$, and every monomial~$u$ of~$\lin{X}$. The differentials are then defined in the same way as the ones of~$\env{A}[X]$, except for $\delta_0[x]=x$ for every $0$-cell~$x$ of~$X$.
A symmetric construction would give a resolution of~$\K$ by left $A$-modules, but the construction of the contracting homotopy requires the notion of a \emph{left} contraction.
\end{remark}

\subsubsection{Algebras of finite derivation type}
\label{SSS:AlgebrasFiniteDerivationType}
For~$n\geq 1$, an algebra is \emph{of finite $n$-derivation type}, $\FDT_n$ for short, if it admits a polygraphic resolution with finitely many $k$-cells for every~$k<n$, and it is \emph{of finite $\infty$-derivation type}, $\FDT_{\infty}$ for short, if it admits a polygraphic resolution of finite type. 
In particular, an algebra is $\FDT_1$ if it is finitely generated, $\FDT_2$ if it is finitely presented, and $\FDT_3$ if it admits a finite coherent presentation. The property $\FDT_3$ corresponds to the finite derivation type condition originally defined by Squier for monoids in~\cite{Squier94}. 
The property $\FDT_n$, for $n\geq 3$ for higher-dimensional categories were introduced in {\cite[2.3.6]{GuiraudMalbos12advances}}.
By definition, $\FDT_{\infty}$ implies $\FDT_n$, and $\FDT_n$ implies $\FDT_p$ for all~$n>p$. 

For~$n\geq 1$, an algebra~$A$ is of \emph{homological type bi-$\FP_n$} (over $\K$) if there is an exact sequence of~$A$-bimodules
\[
0 \fllg A \fllg F_0 \fllg F_1 \fllg \; \cdots \; \fllg F_{n-1} \fllg F_n,
\]
where each~$F_i$ is a finitely generated, free $A$-bimodule, and it is of \emph{homological type bi-$\FP_\infty$} if it is bi-$\mathrm{FP}_n$, for all $n>0$.
Theorems~\ref{T:PolRes->AbRes} and~\ref{T:SquierResolution} give the following implications.

\begin{proposition}
\label{Proposition:FDTnormalisationStrategy}
\begin{enumerate}
\item $\FDT_n$ implies bi-$\FP_n$, for every~$n\geq 1$, and $\FDT_{\infty}$ implies bi-$\FP_{\infty}$.
\item An algebra that admits a finite convergent presentation is $\FDT_{\infty}$.
\end{enumerate}
\end{proposition}

\subsection{Convergence and Koszulness}
\label{S:ConvergenceKoszulness}

\subsubsection{Koszul algebras}

Let us define the map $\ell_N:\Nb\setminus\ens{0}\fl\Nb\setminus\ens{0}$, for a fixed~$N\geq 2$, by
\[
\ell_N(2p) = Np
\qquad\text{and}\qquad
\ell_N(2p+1) = Np+1.
\]
So, if~$X$ is an $\ell_N$-concentrated graded $\infty$-polygraph, then its $0$-cells are in degree $\ell_N(1)=1$, its $1$-cells in degree $\ell_N(2)=N$, its $2$-cells in degree $\ell_N(3)=N+1$, and so on; in particular, the $1$-polygraph underlying~$X$ is $N$-homogeneous.
An $N$-homogeneous algebra~$A$ is said to be \emph{Koszul} if there exists a resolution
\[
0 \fllg A \fllg M_0 \fllg M_1 \fllg M_2 \fllg \cdots
\]
by projective graded $A$-bimodules such that each~$M_n$ is generated by~$M_n^{(\ell_N(n))}$. Note that~\cite[Proposition~4.4]{BergerMarconnet06} implies that one obtains the same notion of Koszulness by replacing bimodules by left or right modules. 
As a consequence of the definition, if~$A$ is Koszul, the vector spaces $\textrm{Tor}_{n,(i)}^{A} (\K,\K)$  vanish for~$i \neq \ell_N(n+1)$, where the first grading in the~$\textrm{Tor}$ refers to the homological degree and the second one to the internal grading of the algebra~$A$. This property of the~$\textrm{Tor}$ groups is an equivalent definition of Koszul algebras, as proved by Berger~{\cite[Theorem 2.11]{Berger01}}.
From Theorem~\ref{T:PolRes->AbRes}, we deduce the following sufficient condition for an algebra to be Koszul.

\begin{proposition}
\label{P:Koszul}
Assume that~$X$ is an $\omega$-concentrated polygraphic resolution of a graded algebra $A$, for some map $\omega:\Nb\setminus\ens{0}\fl\Nb\setminus\ens{0}$. Then each graded $A$-bimodule~$\env{A}[X_n]$ is generated by its component of degree~$\omega(n+1)$. As a consequence, if $\omega=\ell_N$ for some~$N\geq 2$, then~$A$ is Koszul. In particular, if an algebra admits a quadratic convergent presentation, then it is Koszul.
\end{proposition}

\begin{proof}
The $A$-bimodule~$\env{A}[X_n]$ is generated by the $n$-cells of~$X$ which, by hypothesis, are homogeneous of degree~$\omega(n+1)$. Thus, when $\omega=\ell_N$, the resolution~$\env{A}[X]$ satisfies the required properties to prove that~$A$ is Koszul. 
If a graded algebra~$A$ admits a quadratic convergent presentation~$X$, then the critical $n$-cells of $\Std(\red(X))$ are all of the form $x_0\vert\cdots\vert x_n$, where each~$x_i$ belongs to~$X$. As a consequence, the generators of~$\env{A}[X_n]$ lie in degree $n+1=\ell_2(n+1)$.
\end{proof}

\begin{example}
The Koszulness of the (quantum) symmetric algebra is a consequence of Example~\ref{X:SymmetricAlgebraAbelian}.
\end{example}

\begin{proposition}
\label{P:coherentEmpty=>Koszul}
If a graded algebra~$A$ has an $N$-homogeneous presentation that admits an empty acyclic cellular extension, then~$A$ is Koszul. In particular, if~$A$ admits an $N$-homogeneous terminating presentation with no critical branching, then~$A$ is Koszul.
\end{proposition}

\subsubsection{Example} 

Consider the cubical algebra~$A$ presented by the convergent left-monomial $1$-polygraph~$X$ of Example~\ref{X:xyz=x3+y3+z3}.
Since~$X$ has no critical branching, $A$ is Koszul, and admits the resolution
\[
0 \fllg A \fllg \env{A} \fllg \env{A}[x,y,z] \fllg \env{A}[\gamma] \fllg 0 
\]
Thus $\Tor_{0,(0)}^{A}(\K,\K) \simeq \K$, $\Tor_{1,(1)}^{A}(\K,\K) \simeq \K^3$, $\Tor_{2,(3)}^{A}(\K,\K) \simeq \K$, while $\Tor_{n,(i)}^{A}(\K,\K)$ vanishes for other values of~$n$ and~$i$.

Let us compare with the resolution obtained when starting with another presentation.
Using the deglex order induced by $x<y<z$, the leading monomial of $z^3 + y^3+x^3 - xyz$ is~$z^3$. The corresponding terminating presentation~$Y$ of~$A$ has~$x$, $y$ and~$z$ as $0$-cells, and
 $\alpha:z^3\fl xyz - x^3 - y^3$ as only $1$-cell. This presentation is not confluent, because neither of its two critical branchings is:
\[
\vcenter{\xymatrix @R=0.5em {
& xyz^2 - x^3z- y^3z
\\
z^4 
	\ar@/^3ex/ [ur] ^-{\alpha_f z}
	\ar@/_3ex/ [dr] _-{z\alpha_f}
\\
& zxyz - zx^3- zy^3
}}
\qquad\quad
\vcenter{\xymatrix @R=0.5em {
& *+\txt{$xyz^3 - x^3z^2$\\$- y^3z^2$}
	\ar [rr] _-*+\txt{$xy\alpha_f - x^3z^2$\\$- y^3z^2$}
&& *+\txt{$xyxyz - xy^4 - xyx^3$\\$- x^3z^2 - y^3z^2$} 
\\
z^5 
	\ar@/^3ex/ [ur] ^-{\alpha_f z^2}
	\ar@/_3ex/ [dr] _-{z^2\alpha_f}
\\
& z^2xyz - z^2x^3- z^2y^3
}}
\]
The adjunction of the $1$-cell $\beta:zy^3\fl zxyz - zx^3 + y^3z + x^3z -xyz^2$ to~$Y$ yields convergent left-monomial presentation~$Z$ of~$A$, with three critical branchings. 
In that case, Squier's resolution~\ref{T:SquierResolution} obtained from~$Z$ is way larger than~$\Sq(X)$, with cells in every dimension. This induces a resolution~$\env{A}[Z]$ of~$A$ by $A$-bimodules of infinite length, with a non-trivial differential, making homological computations harder than with~$\env{A}[X]$.

\begin{example}
\label{X:PP05Suite}
In Example~\ref{X:PP05}, we have seen that the quadratic $1$-polygraph  
\[
X = 
	\bigpres
		{x, y, z}
		{yz \ofl{\alpha} -x^2 ,\; zy \ofl{\beta} - \lambda^{-1} x^2}.
\]
can be extended into a coherent presentation without any $2$-cell.
By Proposition~\ref{P:coherentEmpty=>Koszul}, it follows that the presented algebra is Koszul.
Moreover we obtain that $\Tor_{0,(0)}^{A}(\K,\K) \simeq \K$, $\Tor_{1,(1)}^{A}(\K,\K) \simeq \K^3$, $\Tor_{2,(2)}^{A}(\K,\K) \simeq \K^2$, and $\Tor_{n,(i)}^{A}(\K,\K)$ vanishes for other values of~$n$ and~$i$.
\end{example}

\begin{proposition}
\label{Proposition:pasKoszul}
Let~$X$ be a polygraphic resolution of an algebra~$A$, whose underlying $n$-polygraph is $\ell_N$-concentrated, for~$N\geq 2$.
If~$X_n^{(i)}$ has strictly more elements than~$X_{n+1}^{(i)}$ for any $i>\ell_N(n+1)$, then~$A$ is not Koszul.
\end{proposition}

\subsubsection{Example}\label{Example:pasKoszul} 

Consider the quadratic algebra~$A$ presented by the terminating $1$-polygraph
\[
X \:=\: \bigpres{x,y}{xy\ofl{\alpha}x^2,\; y^2\ofl{\beta}x^2},
\]
which has two critical branchings, only one of them being confluent:
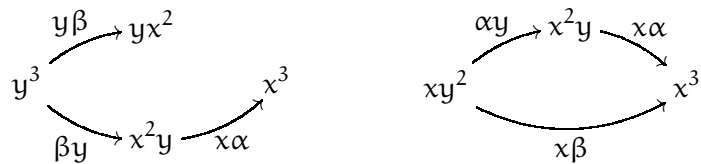
\[
\vcenter{\xymatrix @!C @R=0.25em {
& yx^2
\\
y^3
	\ar @/^/ [ur] ^{y\beta}
	\ar @/_/ [dr] _{\beta y}
&& x^3
\\
& x^2y
	\ar @/_/ [ur] _{x\alpha}
}}
\qquad\qquad
\vcenter{\xymatrix @R=0.5em {
& x^2y
	\ar @/^/ [dr] ^-{x\alpha}
\\
xy^2
	\ar @/^/ [ur] ^-{\alpha y}
	\ar @/_3ex/ [rr] _-{x\beta} ^-{}="tgt"
&& x^3
}}
\]
Adding the $1$-cell $ \gamma : yx^2 \fl x^3$ to~$X$ gives a convergent polygraph~$Y$, so that, by Theorem~\ref{T:Squier}, the following $2$-cells extend the $1$-polygraph~$Y$ into a coherent presentation of~$A$: 
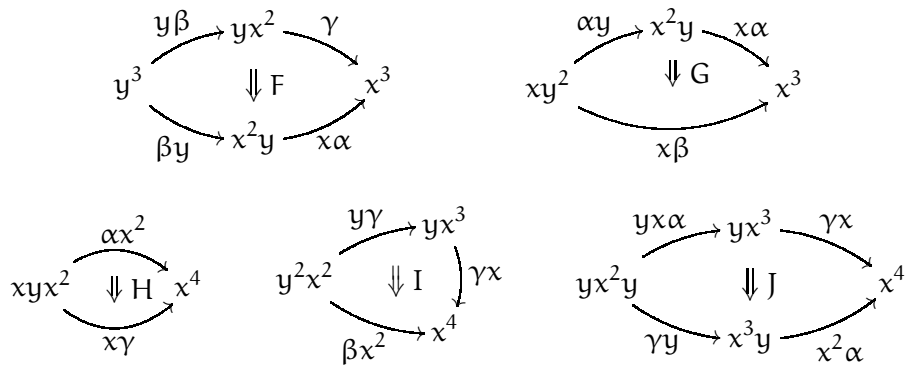
\[
\vcenter{\xymatrix @!C @R=0.25em {
& yx^2
            \ar @/^/ [dr] ^{\gamma}
\\
y^3
	\ar @/^/ [ur] ^{y\beta}
	\ar @/_/ [dr] _{\beta y}
&& x^3
\\
& x^2y
	\ar @/_/ [ur] _{x\alpha}
\ar@2 "1,2"!<0pt,-15pt>;"3,2"!<0pt,15pt> ^-*+{F}
}}
\qquad\qquad
\vcenter{\xymatrix @R=0.5em {
& x^2y
	\ar @/^/ [dr] ^-{x\alpha}
\\
xy^2
	\ar @/^/ [ur] ^-{\alpha y}
	\ar @/_3ex/ [rr] _-{x\beta} ^{}="tgt"
&& x^3
\ar@2 "1,2"!<0pt,-12.5pt>;"2,2"!<0pt,-20pt> ^-*+{G}
}}
\]
\[
\vcenter{\xymatrix @C=3em {
xyx^2
	\ar @/^3ex/ [r] ^-{\alpha x^2} _{}="src"
	\ar @/_3ex/ [r] _-{x\gamma} ^{}="tgt"
& x^4
\ar@2 "src"!<0pt,-10pt>;"tgt"!<0pt,10pt> ^-*+{H}
}}
\qquad
\vcenter{\xymatrix @C=2.5em @R=0.25em {
& yx^3
	\ar @/^/ [dd] ^-{\gamma x}
\\
y^2x^2
	\ar @/^/ [ur] ^-{y\gamma} _{}="src"
	\ar @/_/ [dr] _-{\beta x^2} ^{}="tgt"
\\
&x^4
\ar@2 "src"!<10pt,-10pt>;"tgt"!<10pt,10pt> ^-*+{I}
}}
\qquad
\vcenter{\xymatrix @!C @R=0.25em {
& yx^3
            \ar @/^/ [dr] ^{\gamma x}
\\
yx^2y
	\ar @/^/ [ur] ^{yx\alpha}
	\ar @/_/ [dr] _{\gamma y}
&& x^4
\\
& x^3y
	\ar @/_/ [ur] _{x^2 \alpha}
\ar@2 "1,2"!<0pt,-15pt>;"3,2"!<0pt,15pt> ^-*+{J}
}}
\]
The standard polygraphic resolution has seven critical $3$-cells, inducing the same amount of $3$-cells in~$\Sq(Y)$, with the following $0$-sources: 
$xyx^2y$, $xy^2x^2$, $xy^3$, $yx^2yy$, $y^2x^2y$, $y^3x^2$ and~$y^4$.
Thus, $\Sq(Y)$ has three $2$-cells and two $3$-cells in degree~$4$, preventing~$A$ to be Koszul.

\begin{small}
\bibliographystyle{amsplain}
\bibliography{biblio}
\end{small}

\vfill

\begin{footnotesize}
\auteur{Yves Guiraud}{yves.guiraud@irif.fr}
{INRIA~$\pi r^2$ \\
IRIF, CNRS UMR~8243 \\
Université Paris 7, Case 7014 \\
75205 Paris Cedex 13, France}

\bigskip
\auteur{Eric Hoffbeck}{hoffbeck@math.univ-paris13.fr}
{Université Paris 13, Sorbonne Paris Cité \\
LAGA, CNRS UMR 7539 \\
99 avenue Jean-Baptiste Clément \\
93430 Villetaneuse, France}

\bigskip
\auteur{Philippe Malbos}{malbos@math.univ-lyon1.fr}
{Université de Lyon,\\
Institut Camille Jordan, CNRS UMR 5208\\
Université Claude Bernard Lyon 1\\
43, boulevard du 11 novembre 1918,\\
69622 Villeurbanne Cedex, France}
\end{footnotesize}

\end{document}